
\def\numadroite{oui}

  
  %
  %

  %
  %
  \let\footnotea=\footnote
  \def\anote#1#2{\footnotea{\hbox{$^{#1}$}}{\eightpoint#2}}  
  \catcode`@=12 

 \def\defrefnote#1{\definexref{#1}{{\the\footnotenumber}}{refnotes}}

  %
  %


\ifx\eplain\undefined
  \let\next\relax
\else
  \expandafter\let\expandafter\next\csname endinput\endcsname
\fi
\next
\begingroup
  \catcode44 12 
  \catcode45 12 
  \catcode46 12 
  \catcode58 12 
  \catcode64 11 
  \expandafter\let\expandafter\x\csname ver@ifpdf.sty\endcsname
  \ifcase 0%
    \ifx\x\relax 
    \else
      \ifx\x\empty 
      \else
        1%
      \fi
    \fi
  \else
    \catcode35 6 
    \catcode123 1 
    \catcode125 2 
    \expandafter\ifx\csname PackageInfo\endcsname\relax
      \def\x#1#2{%
        \immediate\write-1{Package #1 Info: #2.}%
      }%
    \else
      \def\x#1#2{\PackageInfo{#1}{#2, stopped}}%
    \fi
    \x{ifpdf}{The package is already loaded}%
    \endgroup
  \fi
\endgroup
\begingroup
  \catcode35 6 
  \catcode40 12 
  \catcode41 12 
  \catcode44 12 
  \catcode45 12 
  \catcode46 12 
  \catcode47 12 
  \catcode58 12 
  \catcode64 11 
  \catcode123 1 
  \catcode125 2 
  \expandafter\ifx\csname ProvidesPackage\endcsname\relax
    \def\x#1#2#3[#4]{\endgroup
      \xdef#1{#4}%
    }%
  \else
    \def\x#1#2[#3]{\endgroup
      #2[{#3}]%
      \ifx#1\@undefined
        \xdef#1{#3}%
      \fi
      \ifx#1\relax
        \xdef#1{#3}%
      \fi
    }%
  \fi
\expandafter\x\csname ver@ifpdf.sty\endcsname
\ProvidesPackage{ifpdf}%
  [2009/04/10 v2.0 Provides the ifpdf switch (HO)]
\begingroup
  \catcode123 1 
  \catcode125 2 
  \def\x{\endgroup
    \expandafter\edef\csname ifpdf@AtEnd\endcsname{%
      \catcode35 \the\catcode35\relax
      \catcode64 \the\catcode64\relax
      \catcode123 \the\catcode123\relax
      \catcode125 \the\catcode125\relax
    }%
  }%
\x
\catcode35 6 
\catcode64 11 
\catcode123 1 
\catcode125 2 
\def\TMP@EnsureCode#1#2{%
  \edef\ifpdf@AtEnd{%
    \ifpdf@AtEnd
    \catcode#1 \the\catcode#1\relax
  }%
  \catcode#1 #2\relax
}
\TMP@EnsureCode{10}{12}
\TMP@EnsureCode{39}{12}
\TMP@EnsureCode{40}{12}
\TMP@EnsureCode{41}{12}
\TMP@EnsureCode{44}{12}
\TMP@EnsureCode{45}{12}
\TMP@EnsureCode{46}{12}
\TMP@EnsureCode{47}{12}
\TMP@EnsureCode{58}{12}
\TMP@EnsureCode{60}{12}
\TMP@EnsureCode{61}{12}
\TMP@EnsureCode{94}{7}
\TMP@EnsureCode{96}{12}
\begingroup
  \expandafter\ifx\csname ifpdf\endcsname\relax
  \else
    \edef\i/{\expandafter\string\csname ifpdf\endcsname}%
    \expandafter\ifx\csname PackageError\endcsname\relax
      \def\x#1#2{%
        \edef\z{#2}%
        \expandafter\errhelp\expandafter{\z}%
        \errmessage{Package ifpdf Error: #1}%
      }%
      \def\y{^^J}%
      \newlinechar=10 %
    \else
      \def\x#1#2{%
        \PackageError{ifpdf}{#1}{#2}%
      }%
      \def\y{\MessageBreak}%
    \fi
    \x{Name clash, \i/ is already defined}{%
      Incompatible versions of \i/ can cause problems,\y
      therefore package loading is aborted.%
    }%
    \endgroup
    \ifpdf@AtEnd
  \fi
\endgroup
\begingroup
  \expandafter\ifx\csname pdfoutput\endcsname\relax
  \else
    \def\skip#1\relax\endgroup{\csname fi\endcsname\endgroup}%
    \skip
  \fi
  \expandafter\ifx\csname directlua\endcsname\relax
    \def\skip#1\endgroup{\csname fi\endcsname\endgroup}%
    \skip
  \fi
  \expandafter\ifx\csname RequirePackage\endcsname\relax
    \input ifluatex.sty\relax
  \else
    \RequirePackage{ifluatex}[2009/04/10]%
  \fi
  \ifluatex
    \ifnum\luatexversion<36 %
    \else
      \directlua{tex.enableprimitives('ifpdf', {'pdfoutput'})}%
      \global\let\pdfoutput\ifpdfpdfoutput
    \fi
  \fi
  \relax
\endgroup
\newif\ifpdf
\begingroup\expandafter\expandafter\expandafter\endgroup
\expandafter\ifx\csname pdfoutput\endcsname\relax
\else
  \ifnum\pdfoutput<1 %
  \else
    \pdftrue
  \fi
\fi
\begingroup
  \expandafter\ifx\csname pdfoutput\endcsname\relax
  \else
    \escapechar=-1 %
    \edef\m{\meaning\pdfoutput}%
    \edef\p{%
      \string p\string d\string f%
      \string o\string u\string t\string p\string u\string t%
    }%
    \ifx\m\p
    \else
      \expandafter\ifx\csname PackageWarningNoLine\endcsname\relax
        \def\PackageWarningNoLine#1#2{%
          \immediate\write16{%
            Package `#1' Warning: #2.%
          }%
        }%
      \fi
      \PackageWarningNoLine{ifpdf}{%
        Someone has redefined \string\\pdfoutput%
      }%
    \fi
  \fi
\endgroup
\begingroup
  \expandafter\ifx\csname PackageInfo\endcsname\relax
    \def\x#1#2{%
    }%
  \else
    \let\x\PackageInfo
    \expandafter\let\csname on@line\endcsname\empty
  \fi
  \x{ifpdf}{pdfTeX in pdf mode \ifpdf\else not \fi detected}%
\endgroup
\ifpdf@AtEnd
\def\makeactive#1{\catcode`#1 = \active \ignorespaces}%
\chardef\letter = 11
\chardef\other = 12
\def\makeatletter{%
  \edef\resetatcatcode{\catcode`\noexpand\@\the\catcode`\@\relax}%
  \catcode`\@11\relax
}%
\def\makeatother{%
  \edef\resetatcatcode{\catcode`\noexpand\@\the\catcode`\@\relax}%
  \catcode`\@12\relax
}%
\edef\leftdisplays{\the\catcode`@}%
\catcode`@ = \letter
\let\@eplainoldatcode = \leftdisplays
\toksdef\toks@ii = 2
\def\uncatcodespecials{%
   \def\do##1{\catcode`##1 = \other}%
   \dospecials
}%
{%
   \makeactive\^^M %
   \long\gdef\letreturn#1{\let^^M = #1}%
}%
\let\@eattoken = \relax  
\def\eattoken{\let\@eattoken = }%
\def\gobble#1{}%
\def\gobbletwo#1#2{}%
\def\gobblethree#1#2#3{}%
\def\identity#1{#1}%
\def\@emptymarkA{\@emptymarkB} 
\def\ifempty#1{\@@ifempty #1\@emptymarkA\@emptymarkB}%
\def\@@ifempty#1#2\@emptymarkB{\ifx #1\@emptymarkA}%
\def\@gobbleminus#1{\ifx-#1\else#1\fi}%
\def\ifinteger#1{\ifcat_\ifnum9<1\@gobbleminus#1 _\else A\fi}%
\def\isinteger{TT\fi\ifinteger}%
\def\@gobblemeaning#1:->{}%
\def\sanitize{\expandafter\@gobblemeaning\meaning}%
\def\ifundefined#1{\expandafter\ifx\csname#1\endcsname\relax}%
\def\csn#1{\csname#1\endcsname}%
\def\ece#1#2{\expandafter#1\csname#2\endcsname}%
\def\expandonce{\expandafter\noexpand}%
\let\@plainwlog = \wlog
\let\wlog = \gobble
\newlinechar = `^^J
\def\loggingall{\tracingcommands\tw@\tracingstats\tw@
   \tracingpages\@ne\tracingoutput\@ne\tracinglostchars\@ne
   \tracingmacros\tw@\tracingparagraphs\@ne\tracingrestores\@ne
   \showboxbreadth\maxdimen\showboxdepth\maxdimen
}%
\def\tracingoff{\tracingonline\z@\tracingcommands\z@\tracingstats\z@
  \tracingpages\z@\tracingoutput\z@\tracinglostchars\z@
  \tracingmacros\z@\tracingparagraphs\z@\tracingrestores\z@
  \showboxbreadth5 \showboxdepth3
}%
\begingroup
  \catcode`\{ = 12 \catcode`\} = 12
  \catcode`\[ = 1 \catcode`\] = 2
  \gdef\lbracechar[{]%
  \gdef\rbracechar[}]%
  \catcode`\% = \other
  \gdef\percentchar[
\def\vpenalty{\ifhmode\par\fi \penalty}%
\def\hpenalty{\ifvmode\leavevmode\fi \penalty}%
\def\iterate{%
  \let\loop@next\relax
  \body
  \let\loop@next\iterate
  \fi
  \loop@next
}%
\def\edefappend#1#2{%
  \toks@ = \expandafter{#1}%
  \edef#1{\the\toks@ #2}%
}%
\def\allowhyphens{\nobreak\hskip\z@skip}%
\long\def\hookprepend{\@hookassign{\the\toks@ii \the\toks@}}%
\long\def\hookappend{\@hookassign{\the\toks@ \the\toks@ii}}%
\let\hookaction = \hookappend 
\long\def\@hookassign#1#2#3{%
  \expandafter\ifx\csname @#2hook\endcsname \relax
    \toks@ = {}%
  \else
    \expandafter\let\expandafter\temp \csname @#2hook\endcsname
    \toks@ = \expandafter{\temp}%
  \fi
  \toks2 = {#3}
  \ece\edef{@#2hook}{#1}%
}%
\long\def\hookactiononce#1#2{%
  \edefappend#2{\global\let\noexpand#2\relax}
  \hookaction{#1}#2%
}%
\def\hookrun#1{%
  \expandafter\ifx\csname @#1hook\endcsname \relax \else
    \def\temp{\csname @#1hook\endcsname}%
    \expandafter\temp
  \fi
}%
\def\setpropertyglobal#1#2#3{\ece\xdef{#1@p#2}{#3}}%
\def\getproperty#1#2{%
  \expandafter\ifx\csname#1@p#2\endcsname\relax
  \else \csname#1@p#2\endcsname
  \fi
}%
\ifx\@undefinedmessage\@undefined
  \def\@undefinedmessage
    {No .aux file; I won't warn you about undefined labels.}%
\fi
\edef\cite{\the\catcode`@}%
\catcode`@ = 11
\let\@oldatcatcode = \cite
\chardef\@letter = 11
\chardef\@other = 12
\def\@innerdef#1#2{\edef#1{\expandafter\noexpand\csname #2\endcsname}}%
\@innerdef\@innernewcount{newcount}%
\@innerdef\@innernewdimen{newdimen}%
\@innerdef\@innernewif{newif}%
\@innerdef\@innernewwrite{newwrite}%
\def\@gobble#1{}%
\ifx\inputlineno\@undefined
   \let\@linenumber = \empty 
\else
   \def\@linenumber{\the\inputlineno:\space}%
\fi
\def\@futurenonspacelet#1{\def\cs{#1}%
   \afterassignment\@stepone\let\@nexttoken=
}%
\begingroup 
\def\\{\global\let\@stoken= }%
\\ 
\endgroup
\def\@stepone{\expandafter\futurelet\cs\@steptwo}%
\def\@steptwo{\expandafter\ifx\cs\@stoken\let\@@next=\@stepthree
   \else\let\@@next=\@nexttoken\fi \@@next}%
\def\@stepthree{\afterassignment\@stepone\let\@@next= }%
\def\@getoptionalarg#1{%
   \let\@optionaltemp = #1%
   \let\@optionalnext = \relax
   \@futurenonspacelet\@optionalnext\@bracketcheck
}%
\def\@bracketcheck{%
   \ifx [\@optionalnext
      \expandafter\@@getoptionalarg
   \else
      \let\@optionalarg = \empty
      \expandafter\@optionaltemp
   \fi
}%
\def\@@getoptionalarg[#1]{%
   \def\@optionalarg{#1}%
   \@optionaltemp
}%
\def\@nnil{\@nil}%
\def\@fornoop#1\@@#2#3{}%
\def\@for#1:=#2\do#3{%
   \edef\@fortmp{#2}%
   \ifx\@fortmp\empty \else
      \expandafter\@forloop#2,\@nil,\@nil\@@#1{#3}%
   \fi
}%
\def\@forloop#1,#2,#3\@@#4#5{\def#4{#1}\ifx #4\@nnil \else
       #5\def#4{#2}\ifx #4\@nnil \else#5\@iforloop #3\@@#4{#5}\fi\fi
}%
\def\@iforloop#1,#2\@@#3#4{\def#3{#1}\ifx #3\@nnil
       \let\@nextwhile=\@fornoop \else
      #4\relax\let\@nextwhile=\@iforloop\fi\@nextwhile#2\@@#3{#4}%
}%
\@innernewif\if@fileexists
\def\@testfileexistence{\@getoptionalarg\@finishtestfileexistence}%
\def\@finishtestfileexistence#1{%
   \begingroup
      \def\extension{#1}%
      \immediate\openin0 =
         \ifx\@optionalarg\empty\jobname\else\@optionalarg\fi
         \ifx\extension\empty \else .#1\fi
         \space
      \ifeof 0
         \global\@fileexistsfalse
      \else
         \global\@fileexiststrue
      \fi
      \immediate\closein0
   \endgroup
}%
\toks0 = {%
\def\bibliographystyle#1{%
   \@readauxfile
   \@writeaux{\string\bibstyle{#1}}%
}%
\let\bibstyle = \@gobble
\let\bblfilebasename = \jobname
\def\bibliography#1{%
   \@readauxfile
   \@writeaux{\string\bibdata{#1}}%
   \@testfileexistence[\bblfilebasename]{bbl}%
   \if@fileexists
      \nobreak
      \@readbblfile
   \fi
}%
\let\bibdata = \@gobble
\def\nocite#1{%
   \@readauxfile
   \@writeaux{\string\citation{#1}}%
}%
\@innernewif\if@notfirstcitation
\def\cite{\@getoptionalarg\@cite}%
\def\@cite#1{%
   \let\@citenotetext = \@optionalarg
   \printcitestart
   \nocite{#1}%
   \@notfirstcitationfalse
   \@for \@citation :=#1\do
   {%
      \expandafter\@onecitation\@citation\@@
   }%
   \ifx\empty\@citenotetext\else
      \printcitenote{\@citenotetext}%
   \fi
   \printcitefinish
}%
\def\@onecitation#1\@@{%
   \if@notfirstcitation
      \printbetweencitations
   \fi
   \expandafter \ifx \csname\@citelabel{#1}\endcsname \relax
      \if@citewarning
         \message{\@linenumber Undefined citation `#1'.}%
      \fi
      \expandafter\gdef\csname\@citelabel{#1}\endcsname{%
         {\tt
            \escapechar = -1
            \nobreak\hskip0pt
            \expandafter\string\csname#1\endcsname
            \nobreak\hskip0pt
         }%
      }%
   \fi
   \printcitepreitem{#1}%
   \csname\@citelabel{#1}\endcsname
   \printcitepostitem
   \@notfirstcitationtrue
}%
\def\@citelabel#1{b@#1}%
\def\@citedef#1#2{{\expandafter\gdef\csname\@citelabel{#1}\endcsname{#2}}}%
\def\@readbblfile{%
   \ifx\@itemnum\@undefined
      \@innernewcount\@itemnum
   \fi
   \begingroup
      \ifx\begin\@undefined
         \def\begin##1##2{%
            \setbox0 = \hbox{\biblabelcontents{##2}}%
            \biblabelwidth = \wd0
         }%
         \let\end = \@gobble 
      \fi
      \@itemnum = 0
      \def\bibitem{\@getoptionalarg\@bibitem}%
      \def\@bibitem{%
         \ifx\@optionalarg\empty
            \expandafter\@numberedbibitem
         \else
            \expandafter\@alphabibitem
         \fi
      }%
      \def\@alphabibitem##1{%
         \expandafter \xdef\csname\@citelabel{##1}\endcsname {\@optionalarg}%
         \ifx\biblabelprecontents\@undefined
            \let\biblabelprecontents = \relax
         \fi
         \ifx\biblabelpostcontents\@undefined
            \let\biblabelpostcontents = \hss
         \fi
         \@finishbibitem{##1}%
      }%
      \def\@numberedbibitem##1{%
         \advance\@itemnum by 1
         \expandafter \xdef\csname\@citelabel{##1}\endcsname{\number\@itemnum}%
         \ifx\biblabelprecontents\@undefined
            \let\biblabelprecontents = \hss
         \fi
         \ifx\biblabelpostcontents\@undefined
            \let\biblabelpostcontents = \relax
         \fi
         \@finishbibitem{##1}%
      }%
      \def\@finishbibitem##1{%
         \bblitemhook{##1}%
         \biblabelprint{\csname\@citelabel{##1}\endcsname}%
         \@writeaux{\string\@citedef{##1}{\csname\@citelabel{##1}\endcsname}}%
         \ignorespaces
      }%
      \ifx\undefined\em \let\em=\bblem \fi
      \ifx\undefined\emph \let\emph=\bblemph \fi
      \ifx\undefined\mbox \let\mbox=\bblmbox \fi
      \ifx\undefined\newblock \let\newblock=\bblnewblock \fi
      \ifx\undefined\sc \let\sc=\bblsc \fi
      \ifx\undefined\textbf \let\textbf=\bbltextbf \fi
      \frenchspacing
      \clubpenalty = 4000 \widowpenalty = 4000
      \tolerance = 10000 \hfuzz = .5pt
      \everypar = {\hangindent = \biblabelwidth
                      \advance\hangindent by \biblabelextraspace}%
      \bblrm
      \parskip = 1.5ex plus .5ex minus .5ex
      \biblabelextraspace = .5em
      \bblhook
      \input \bblfilebasename.bbl
   \endgroup
}%
\@innernewdimen\biblabelwidth
\@innernewdimen\biblabelextraspace
\def\biblabelprint#1{%
   \noindent
   \hbox to \biblabelwidth{%
      \biblabelprecontents
      \biblabelcontents{#1}%
      \biblabelpostcontents
   }%
   \kern\biblabelextraspace
}%
\def\biblabelcontents#1{{\bblrm [#1]}}%
\def\bblrm{\rm}%
\def\bblem{\it}%
\def\bblemph#1{{\bblem #1\/}}
\def\textbf#1{{\bf #1}}
\def\bblmbox{\leavevmode\hbox}
\def\bblsc{\ifx\@scfont\@undefined
              \font\@scfont = cmcsc10
           \fi
           \@scfont
}%
\def\bblnewblock{\hskip .11em plus .33em minus .07em }%
\let\bblhook = \empty
\let\bblitemhook = \@gobble
\def\printcitestart{[}
\def\printcitefinish{]}
\def\printbetweencitations{, }
\let\printcitepreitem\@gobble 
\let\printcitepostitem\empty
\def\printcitenote#1{, #1}
\let\citation = \@gobble
\@innernewcount\@numparams
\ifx\newcommand\undefined
\def\newcommand#1{%
   \def\@commandname{#1}%
   \@getoptionalarg\@continuenewcommand
}%
\fi
\ifx\renewcommand\undefined
\let\renewcommand = \newcommand
\fi
\ifx\providecommand\undefined
\def\providecommand#1{%
   \def\@commandname{#1}%
   \expandafter\ifx\@commandname \@undefined
     \let\cs=\@continuenewcommand  
   \else
     \let\cs=\@gobble              
   \fi
   \@getoptionalarg\cs
}%
\fi
\def\@continuenewcommand{%
   \@numparams = \ifx\@optionalarg\empty 0\else\@optionalarg \fi \relax
   \@newcommand
}%
\def\@newcommand#1{%
   \def\@startdef{\expandafter\edef\@commandname}%
   \ifnum\@numparams=0
      \let\@paramdef = \empty
   \else
      \ifnum\@numparams>9
         \errmessage{\the\@numparams\space is too many parameters}%
      \else
         \ifnum\@numparams<0
            \errmessage{\the\@numparams\space is too few parameters}%
         \else
            \edef\@paramdef{%
               \ifcase\@numparams
                  \empty  No arguments.
               \or ####1%
               \or ####1####2%
               \or ####1####2####3%
               \or ####1####2####3####4%
               \or ####1####2####3####4####5%
               \or ####1####2####3####4####5####6%
               \or ####1####2####3####4####5####6####7%
               \or ####1####2####3####4####5####6####7####8%
               \or ####1####2####3####4####5####6####7####8####9%
               \fi
            }%
         \fi
      \fi
   \fi
   \expandafter\@startdef\@paramdef{#1}%
}%
}%
\ifx\nobibtex\@undefined \the\toks0 \fi
\def\@readauxfile{%
   \if@auxfiledone \else 
      \global\@auxfiledonetrue
      \@testfileexistence{aux}%
      \if@fileexists
         \begingroup
            \endlinechar = -1
            \catcode`@ = 11
            \input \jobname.aux
         \endgroup
      \else
         \message{\@undefinedmessage}%
         \global\@citewarningfalse
      \fi
      \immediate\openout\@auxfile = \jobname.aux
   \fi
}%
\newif\if@auxfiledone
\ifx\noauxfile\@undefined \else \@auxfiledonetrue\fi
\@innernewwrite\@auxfile
\def\@writeaux#1{\ifx\noauxfile\@undefined \write\@auxfile{#1}\fi}%
\ifx\@undefinedmessage\@undefined
   \def\@undefinedmessage{No .aux file; I won't give you warnings about
                          undefined citations.}%
\fi
\@innernewif\if@citewarning
\ifx\noauxfile\@undefined \@citewarningtrue\fi
\catcode`@ = \@oldatcatcode
\let\auxfile = \@auxfile
\let\for = \@for
\let\futurenonspacelet = \@futurenonspacelet
\def\iffileexists{\if@fileexists}%
\let\innerdef = \@innerdef
\let\innernewcount = \@innernewcount
\let\innernewdimen = \@innernewdimen
\let\innernewif = \@innernewif
\let\innernewwrite = \@innernewwrite
\let\linenumber = \@linenumber
\let\readauxfile = \@readauxfile
\let\spacesub = \@spacesub
\let\testfileexistence = \@testfileexistence
\let\writeaux = \@writeaux
\def\innerinnerdef#1{\expandafter\innerdef\csname inner#1\endcsname{#1}}%
\innerinnerdef{newbox}%
\innerinnerdef{newfam}%
\innerinnerdef{newhelp}%
\innerinnerdef{newinsert}%
\innerinnerdef{newlanguage}%
\innerinnerdef{newmuskip}%
\innerinnerdef{newread}%
\innerinnerdef{newskip}%
\innerinnerdef{newtoks}%
\def\immediatewriteaux#1{%
  \ifx\noauxfile\@undefined
    \immediate\write\@auxfile{#1}%
  \fi
}%
\def\bblitemhook#1{\gdef\@hlbblitemlabel{#1}}%
\def\biblabelprint#1{%
   \noindent
   \hbox to \biblabelwidth{%
      \hldest@impl{bib}{\@hlbblitemlabel}%
      \biblabelprecontents
      \biblabelcontents{#1}%
      \biblabelpostcontents
   }%
   \kern\biblabelextraspace
}%
\def\eplainprintcitepreitem#1{\hlstart@impl{cite}{#1}}%
\def\eplainprintcitepostitem{\hlend@impl{cite}}%
\def\printcitepreitem#1{%
  \testfileexistence[\bblfilebasename]{bbl}%
  \iffileexists
    \global\let\printcitepreitem\eplainprintcitepreitem
    \global\let\printcitepostitem\eplainprintcitepostitem
  \else
    \global\let\printcitepreitem\gobble
    \global\let\printcitepostitem\relax
  \fi
  \printcitepreitem{#1}%
}%
\def\@Nnil{\@Nil}%
\def\@Fornoop#1\@@#2#3{}%
\def\For#1:=#2\do#3{%
   \edef\@Fortmp{#2}%
   \ifx\@Fortmp\empty \else
      \expandafter\@Forloop#2,\@Nil,\@Nil\@@#1{#3}%
   \fi
}%
\def\@Forloop#1,#2,#3\@@#4#5{\@Fordef#1\@@#4\ifx #4\@Nnil \else
       #5\@Fordef#2\@@#4\ifx #4\@Nnil \else#5\@iForloop #3\@@#4{#5}\fi\fi
}%
\def\@iForloop#1,#2\@@#3#4{\@Fordef#1\@@#3\ifx #3\@Nnil
       \let\@Nextwhile=\@Fornoop \else
      #4\relax\let\@Nextwhile=\@iForloop\fi\@Nextwhile#2\@@#3{#4}%
}%
\def\@Forspc{ }%
\def\@Fordef{\futurelet\@Fortmp\@@Fordef}
\def\@@Fordef{%
  \expandafter\ifx\@Forspc\@Fortmp 
    \expandafter\@Fortrim
  \else
    \expandafter\@@@Fordef
  \fi
}%
\expandafter\def\expandafter\@Fortrim\@Forspc#1\@@{\@Fordef#1\@@}%
\def\@@@Fordef#1\@@#2{\def#2{#1}}%
\def\tmpfileextension{.tmp}%
\let\tmpfilebasename = \jobname
\ifx\eTeXversion\undefined
  \innernewwrite\eplain@tmpfile
  \def\scantokens#1{%
    \toks@={#1}%
    \immediate\openout\eplain@tmpfile=\tmpfilebasename\tmpfileextension
    \immediate\write\eplain@tmpfile{\the\toks@}%
    \immediate\closeout\eplain@tmpfile
    \input \tmpfilebasename\tmpfileextension\relax
  }%
\fi
\begingroup
   \makeactive\^^M \makeactive\ 
\gdef\obeywhitespace{%
\makeactive\^^M\def^^M{\par\futurelet\next\@finishobeyedreturn}%
\makeactive\ \let =\ %
\aftergroup\@removebox%
\futurelet\next\@finishobeywhitespace%
}%
\gdef\@finishobeywhitespace{{%
\ifx\next %
\aftergroup\@obeywhitespaceloop%
\else\ifx\next^^M%
\aftergroup\gobble%
\fi\fi}}%
\gdef\@finishobeyedreturn{%
\ifx\next^^M\vskip\blanklineskipamount\fi%
\indent%
}%
\endgroup
\def\@obeywhitespaceloop#1{\futurelet\next\@finishobeywhitespace}%
\def\@removebox{%
  \ifhmode
    \setbox0 = \lastbox
    \ifdim\wd0=\parindent
      \setbox2 = \hbox{\unhcopy0}
      \ifdim\wd2=0pt
        \ignorespaces
      \else
        \box0 
      \fi
    \else
       \box0 
    \fi
  \fi
}%
\newskip\blanklineskipamount
\blanklineskipamount = 0pt
\def\frac#1/#2{\leavevmode
   \kern.1em \raise .5ex \hbox{\the\scriptfont0 #1}%
   \kern-.1em $/$%
   \kern-.15em \lower .25ex \hbox{\the\scriptfont0 #2}%
}%
\newdimen\hruledefaultheight  \hruledefaultheight = 0.4pt
\newdimen\hruledefaultdepth   \hruledefaultdepth = 0.0pt
\newdimen\vruledefaultwidth   \vruledefaultwidth = 0.4pt
\def\ehrule{\hrule height\hruledefaultheight depth\hruledefaultdepth}%
\def\evrule{\vrule width\vruledefaultwidth}%
\ifx\sc\undefined
    \def\sc{%
      \expandafter\ifx\the\scriptfont\fam\nullfont
        \font\temp = cmr7 \temp
      \else
        \the\scriptfont\fam
      \fi
      \def\uppercasesc{\char\uccode`}%
    }%
\fi
\ifx\uppercasesc\undefined
  \let\uppercasesc = \relax
\fi
\def\TeX{T\kern-.1667em\lower.5ex\hbox{E}\kern-.125emX\spacefactor1000 }%
\ifx\AmS\undefined
    \def\AmS{{\the\textfont2 A}\kern-.1667em\lower.5ex\hbox
        {\the\textfont2 M}\kern-.125em{\the\textfont2 S}}
\fi
\ifx\AMS\undefined \let\AMS=\AmS \fi
\ifx\AmSLaTeX\undefined
    \def\AmSLaTeX{\AmS-\LaTeX}
\fi
\ifx\AMSLaTeX\undefined \let\AMSLaTeX=\AmSLaTeX \fi
\ifx\AmSTeX\undefined
    \def\AmSTeX{$\cal A$\kern-.1667em\lower.5ex\hbox{$\cal M$}%
            \kern-.125em$\cal S$-\TeX}%
\fi
\ifx\AMSTEX\undefined \let\AMSTEX=\AmSTeX \fi
\ifx\AMSTeX\undefined \let\AMSTeX=\AmSTeX \fi
\ifx\BibTeX\undefined
    \def\BibTeX{B{\sc \uppercasesc i\kern-.025em \uppercasesc b}\kern-.08em
                \TeX}%
\fi
\ifx\BIBTeX\undefined \let\BIBTeX=\BibTeX \fi
\ifx\BIBTEX\undefined \let\BIBTEX=\BibTeX \fi
\ifx\LAMSTeX\undefined
    \def\LAMSTeX{L\raise.42ex\hbox{\kern-.3em\the\scriptfont2 A}%
                 \kern-.2em\lower.376ex\hbox{\the\textfont2 M}%
                 \kern-.125em {\the\textfont2 S}-\TeX}%
\fi
\ifx\LamSTeX\undefined \let\LamSTeX=\LAMSTeX \fi
\ifx\LAmSTeX\undefined \let\LAmSTeX=\LAMSTeX \fi
\ifx\LaTeX\undefined
    \def\LaTeX{L\kern-.36em\raise.3ex\hbox{\sc \uppercasesc a}\kern-.15em\TeX}%
\fi
\ifx\LATEX\undefined \let\LATEX=\LaTeX \fi
\ifx\LaTeXe\undefined
    \def\LaTeXe{\LaTeX{}\kern.05em2$_{\textstyle\varepsilon}$}
\fi
\ifx\MF\undefined
    \ifx\manfnt\undefined
            \font\manfnt=logo10
    \fi
    \ifx\manfntsl\undefined
            \font\manfntsl=logosl10
    \fi
    \def\MF{{\ifdim\fontdimen1\font>0pt \let\manfnt = \manfntsl \fi
      {\manfnt META}\-{\manfnt FONT}}\spacefactor1000 }%
\fi
\ifx\METAFONT\undefined \let\METAFONT=\MF \fi
\ifx\SLITEX\undefined
    \def\SLITEX{S\kern-.065em L\kern-.18em\raise.32ex\hbox{i}\kern-.03em\TeX}%
\fi
\ifx\SLiTeX\undefined \let\SLiTeX=\SLITEX \fi
\ifx\SliTeX\undefined \let\SliTeX=\SLITEX \fi
\ifx\SLITeX\undefined \let\SLITeX=\SLITEX \fi
\edef\path{\the\catcode`@}%
\catcode`@ = 11
\let\@oldatcatcode = \path
\newcount \c@tcode
\newcount \c@unter
\newif \ifspecialpathdelimiters
\begingroup
\catcode `\ = 10
\gdef \passivesp@ce { }%
\catcode `\ = 13\relax%
\gdef\activesp@ce{ }%
\endgroup
\def \discretionaries 
    {\begingroup
        \c@tcodes = 13
        \discr@tionaries
    }%
\def \discr@tionaries #1
    {\def \discr@ti@naries ##1#1
         {\endgroup
          \def \discr@ti@n@ries ####1
              {\if   \noexpand ####1\noexpand #1%
                     \let \n@xt = \relax
               \else
                     \catcode `####1 = 13
                     \def ####1{\discretionary
                                  {\char `####1}{}{\char `####1}}%
                     \let \n@xt = \discr@ti@n@ries
               \fi
               \n@xt
              }%
          \def \discr@ti@n@ri@s {\discr@ti@n@ries ##1#1}%
         }%
     \discr@ti@naries
    }%
\let\pathafterhook = \relax
\def \path
    {\ifspecialpathdelimiters
        \begingroup
        \c@tcodes = 12
        \def \endp@th {\endgroup \endgroup \pathafterhook}%
     \else
        \def \endp@th {\endgroup \pathafterhook}%
     \fi
     \p@th
    }%
\def \p@th #1
    {\begingroup
        \tt
        \c@tcode = \catcode `#1
        \discr@ti@n@ri@s
        \catcode `\ = \active
        \expandafter \edef \activesp@ce {\passivesp@ce \hbox {}}%
        \catcode `#1 = \c@tcode
        \def \p@@th ##1#1
            {\leavevmode \hbox {}##1%
             \endp@th
            }%
     \p@@th
    }%
\def \c@tcodes {\afterassignment \c@tc@des \c@tcode}%
\def \c@tc@des
    {\c@unter = 0
     \loop
            \ifnum \catcode \c@unter = \c@tcode
            \else
                \catcode \c@unter = \c@tcode
            \fi
     \ifnum \c@unter < 255
            \advance \c@unter by 1
     \repeat
     \catcode `\ = 10
    }%
\catcode `\@ = \@oldatcatcode
\discretionaries |~!@$
\ifx\eTeX\undefined
  \def\eTeX{\hbox{\mathsurround=0pt $\varepsilon$-\kern-.125em\TeX}}%
\fi
\ifx\ExTeX\undefined
  \def\ExTeX{\hbox{\mathsurround=0pt
    $\textstyle\varepsilon_{\kern-0.15em\cal{X}}$\kern-.2em\TeX}}%
\fi
\def\eplain@Xe@reflect#1{%
  \ifx\reflectbox\undefined
    \errmessage{A graphics package must be loaded for \string\XeTeX}%
  \else
    \ifdim \fontdimen1\font>0pt
      \raise 1.75ex \hbox{\kern.1em\rotatebox{180}{#1}}\kern-.1em
    \else
      \reflectbox{#1}%
    \fi
  \fi
}%
\def\eplain@Xe#1{\leavevmode
  \smash{\hbox{X%
    \setbox0=\hbox{\TeX}\setbox2=\hbox{E}%
    \lower\dp0\hbox{\raise\dp2\hbox{\kern-.125em\eplain@Xe@reflect{E}}}%
    \kern-.1667em #1}}}%
\ifx\XeTeX\undefined
  \def\XeTeX{\eplain@Xe\TeX}%
\fi
\ifx\XeLaTeX\undefined
  \def\XeLaTeX{\eplain@Xe{\thinspace\LaTeX}}%
\fi
\def\blackbox{\vrule height .8ex width .6ex depth -.2ex \relax}
\def\makeblankbox#1#2{%
  \ifvoid0
    \errhelp = \@makeblankboxhelp
    \errmessage{Box 0 is void}%
  \fi
  \hbox{\lower\dp0
    \vbox{\hidehrule{#1}{#2}%
      \kern -#1
      \hbox to \wd0{\hidevrule{#1}{#2}%
        \raise\ht0\vbox to #1{}
        \lower\dp0\vtop to #1{}
        \hfil\hidevrule{#2}{#1}%
      }%
      \kern-#1\hidehrule{#2}{#1}%
    }%
  }%
}%
\newhelp\@makeblankboxhelp{Assigning to the dimensions of a void^^J%
  box has no effect.  Do `\string\setbox0=\string\null' before you^^J%
  define its dimensions.}%
\def\hidehrule#1#2{\kern-#1\hrule height#1 depth#2 \kern-#2}%
\def\hidevrule#1#2{%
  \kern-#1%
  \dimen@=#1\advance\dimen@ by #2%
  \vrule width\dimen@
  \kern-#2%
}%
\newdimen\boxitspace \boxitspace = 3pt
\long\def\boxit#1{%
  \vbox{%
    \ehrule
    \hbox{%
      \evrule
      \kern\boxitspace
      \vbox{\kern\boxitspace \parindent = 0pt #1\kern\boxitspace}%
      \kern\boxitspace
      \evrule
    }%
    \ehrule
  }%
}%
\def\numbername#1{\ifcase#1%
   zero%
   \or one%
   \or two%
   \or three%
   \or four%
   \or five%
   \or six%
   \or seven%
   \or eight%
   \or nine%
   \or ten%
   \or #1%
   \fi
}%
\let\@plainnewif = \newif
\let\@plainnewdimen = \newdimen
\let\newif = \innernewif
\let\newdimen = \innernewdimen
\edef\@eplainoldandcode{\the\catcode`& }%
\catcode`& = 11
\toks0 = {%
\edef\thinlines{\the\catcode`@ }%
\catcode`@ = 11
\let\@oldatcatcode = \thinlines
\def\smash@@{\relax 
  \ifmmode\def\next{\mathpalette\mathsm@sh}\else\let\next\makesm@sh
  \fi\next}
\def\makesm@sh#1{\setbox\z@\hbox{#1}\finsm@sh}
\def\mathsm@sh#1#2{\setbox\z@\hbox{$\m@th#1{#2}$}\finsm@sh}
\def\finsm@sh{\ht\z@\z@ \dp\z@\z@ \box\z@}
\edef\@oldandcatcode{\the\catcode`& }%
\catcode`& = 11
\def\&whilenoop#1{}%
\def\&whiledim#1\do #2{\ifdim #1\relax#2\&iwhiledim{#1\relax#2}\fi}%
\def\&iwhiledim#1{\ifdim #1\let\&nextwhile=\&iwhiledim 
        \else\let\&nextwhile=\&whilenoop\fi\&nextwhile{#1}}%
\newif\if&negarg
\newdimen\&wholewidth
\newdimen\&halfwidth
\font\tenln=line10
\def\thinlines{\let\&linefnt\tenln \let\&circlefnt\tencirc
  \&wholewidth\fontdimen8\tenln \&halfwidth .5\&wholewidth}%
\def\thicklines{\let\&linefnt\tenlnw \let\&circlefnt\tencircw
  \&wholewidth\fontdimen8\tenlnw \&halfwidth .5\&wholewidth}%
\def\drawline(#1,#2)#3{\&xarg #1\relax \&yarg #2\relax \&linelen=#3\relax
  \ifnum\&xarg =0 \&vline \else \ifnum\&yarg =0 \&hline \else \&sline\fi\fi}%
\def\&sline{\leavevmode
  \ifnum\&xarg< 0 \&negargtrue \&xarg -\&xarg \&yyarg -\&yarg
  \else \&negargfalse \&yyarg \&yarg \fi
  \ifnum \&yyarg >0 \&tempcnta\&yyarg \else \&tempcnta -\&yyarg \fi
  \ifnum\&tempcnta>6 \&badlinearg \&yyarg0 \fi
  \ifnum\&xarg>6 \&badlinearg \&xarg1 \fi
  \setbox\&linechar\hbox{\&linefnt\&getlinechar(\&xarg,\&yyarg)}%
  \ifnum \&yyarg >0 \let\&upordown\raise \&clnht\z@
  \else\let\&upordown\lower \&clnht \ht\&linechar\fi
  \&clnwd=\wd\&linechar
  \&whiledim \&clnwd <\&linelen \do {%
    \&upordown\&clnht\copy\&linechar
    \advance\&clnht \ht\&linechar
    \advance\&clnwd \wd\&linechar
  }%
  \advance\&clnht -\ht\&linechar
  \advance\&clnwd -\wd\&linechar
  \&tempdima\&linelen\advance\&tempdima -\&clnwd
  \&tempdimb\&tempdima\advance\&tempdimb -\wd\&linechar
  \hskip\&tempdimb \multiply\&tempdima \@m
  \&tempcnta \&tempdima \&tempdima \wd\&linechar \divide\&tempcnta \&tempdima
  \&tempdima \ht\&linechar \multiply\&tempdima \&tempcnta
  \divide\&tempdima \@m
  \advance\&clnht \&tempdima
  \ifdim \&linelen <\wd\&linechar \hskip \wd\&linechar
  \else\&upordown\&clnht\copy\&linechar\fi}%
\def\&hline{\vrule height \&halfwidth depth \&halfwidth width \&linelen}%
\def\&getlinechar(#1,#2){\&tempcnta#1\relax\multiply\&tempcnta 8
  \advance\&tempcnta -9 \ifnum #2>0 \advance\&tempcnta #2\relax\else
  \advance\&tempcnta -#2\relax\advance\&tempcnta 64 \fi
  \char\&tempcnta}%
\def\drawvector(#1,#2)#3{\&xarg #1\relax \&yarg #2\relax
  \&tempcnta \ifnum\&xarg<0 -\&xarg\else\&xarg\fi
  \ifnum\&tempcnta<5\relax \&linelen=#3\relax
    \ifnum\&xarg =0 \&vvector \else \ifnum\&yarg =0 \&hvector
    \else \&svector\fi\fi\else\&badlinearg\fi}%
\def\&hvector{\ifnum\&xarg<0 \rlap{\&linefnt\&getlarrow(1,0)}\fi \&hline
  \ifnum\&xarg>0 \llap{\&linefnt\&getrarrow(1,0)}\fi}%
\def\&vvector{\ifnum \&yarg <0 \&downvector \else \&upvector \fi}%
\def\&svector{\&sline
  \&tempcnta\&yarg \ifnum\&tempcnta <0 \&tempcnta=-\&tempcnta\fi
  \ifnum\&tempcnta <5 
    \if&negarg\ifnum\&yarg>0                   
      \llap{\lower\ht\&linechar\hbox to\&linelen{\&linefnt
        \&getlarrow(\&xarg,\&yyarg)\hss}}\else 
      \llap{\hbox to\&linelen{\&linefnt\&getlarrow(\&xarg,\&yyarg)\hss}}\fi
    \else\ifnum\&yarg>0                        
      \&tempdima\&linelen \multiply\&tempdima\&yarg
      \divide\&tempdima\&xarg \advance\&tempdima-\ht\&linechar
      \raise\&tempdima\llap{\&linefnt\&getrarrow(\&xarg,\&yyarg)}\else
      \&tempdima\&linelen \multiply\&tempdima-\&yarg 
      \divide\&tempdima\&xarg
      \lower\&tempdima\llap{\&linefnt\&getrarrow(\&xarg,\&yyarg)}\fi\fi
  \else\&badlinearg\fi}%
\def\&getlarrow(#1,#2){\ifnum #2 =\z@ \&tempcnta='33\else
\&tempcnta=#1\relax\multiply\&tempcnta \sixt@@n \advance\&tempcnta
-9 \&tempcntb=#2\relax\multiply\&tempcntb \tw@
\ifnum \&tempcntb >0 \advance\&tempcnta \&tempcntb\relax
\else\advance\&tempcnta -\&tempcntb\advance\&tempcnta 64
\fi\fi\char\&tempcnta}%
\def\&getrarrow(#1,#2){\&tempcntb=#2\relax
\ifnum\&tempcntb < 0 \&tempcntb=-\&tempcntb\relax\fi
\ifcase \&tempcntb\relax \&tempcnta='55 \or 
\ifnum #1<3 \&tempcnta=#1\relax\multiply\&tempcnta
24 \advance\&tempcnta -6 \else \ifnum #1=3 \&tempcnta=49
\else\&tempcnta=58 \fi\fi\or 
\ifnum #1<3 \&tempcnta=#1\relax\multiply\&tempcnta
24 \advance\&tempcnta -3 \else \&tempcnta=51\fi\or 
\&tempcnta=#1\relax\multiply\&tempcnta
\sixt@@n \advance\&tempcnta -\tw@ \else
\&tempcnta=#1\relax\multiply\&tempcnta
\sixt@@n \advance\&tempcnta 7 \fi\ifnum #2<0 \advance\&tempcnta 64 \fi
\char\&tempcnta}%
\def\&vline{\ifnum \&yarg <0 \&downline \else \&upline\fi}%
\def\&upline{\hbox to \z@{\hskip -\&halfwidth \vrule width \&wholewidth
   height \&linelen depth \z@\hss}}%
\def\&downline{\hbox to \z@{\hskip -\&halfwidth \vrule width \&wholewidth
   height \z@ depth \&linelen \hss}}%
\def\&upvector{\&upline\setbox\&tempboxa\hbox{\&linefnt\char'66}\raise 
     \&linelen \hbox to\z@{\lower \ht\&tempboxa\box\&tempboxa\hss}}%
\def\&downvector{\&downline\lower \&linelen
      \hbox to \z@{\&linefnt\char'77\hss}}%
\def\&badlinearg{\errmessage{Bad \string\arrow\space argument.}}%
\thinlines
\countdef\&xarg     0
\countdef\&yarg     2
\countdef\&yyarg    4
\countdef\&tempcnta 6
\countdef\&tempcntb 8
\dimendef\&linelen  0
\dimendef\&clnwd    2
\dimendef\&clnht    4
\dimendef\&tempdima 6
\dimendef\&tempdimb 8
\chardef\@arrbox    0
\chardef\&linechar  2
\chardef\&tempboxa  2           
\let\lft^%
\let\rt_
\newif\if@pslope 
\def\@findslope(#1,#2){\ifnum#1>0
  \ifnum#2>0 \@pslopetrue \else\@pslopefalse\fi \else
  \ifnum#2>0 \@pslopefalse \else\@pslopetrue\fi\fi}%
\def\generalsmap(#1,#2){\getm@rphposn(#1,#2)\plnmorph\futurelet\next\addm@rph}%
\def\sline(#1,#2){\setbox\@arrbox=\hbox{\drawline(#1,#2){\sarrowlength}}%
  \@findslope(#1,#2)\d@@blearrfalse\generalsmap(#1,#2)}%
\def\arrow(#1,#2){\setbox\@arrbox=\hbox{\drawvector(#1,#2){\sarrowlength}}%
  \@findslope(#1,#2)\d@@blearrfalse\generalsmap(#1,#2)}%
\newif\ifd@@blearr
\def\bisline(#1,#2){\@findslope(#1,#2)%
  \if@pslope \let\@upordown\raise \else \let\@upordown\lower\fi
  \getch@nnel(#1,#2)\setbox\@arrbox=\hbox{\@upordown\@vchannel
    \rlap{\drawline(#1,#2){\sarrowlength}}%
      \hskip\@hchannel\hbox{\drawline(#1,#2){\sarrowlength}}}%
  \d@@blearrtrue\generalsmap(#1,#2)}%
\def\biarrow(#1,#2){\@findslope(#1,#2)%
  \if@pslope \let\@upordown\raise \else \let\@upordown\lower\fi
  \getch@nnel(#1,#2)\setbox\@arrbox=\hbox{\@upordown\@vchannel
    \rlap{\drawvector(#1,#2){\sarrowlength}}%
      \hskip\@hchannel\hbox{\drawvector(#1,#2){\sarrowlength}}}%
  \d@@blearrtrue\generalsmap(#1,#2)}%
\def\adjarrow(#1,#2){\@findslope(#1,#2)%
  \if@pslope \let\@upordown\raise \else \let\@upordown\lower\fi
  \getch@nnel(#1,#2)\setbox\@arrbox=\hbox{\@upordown\@vchannel
    \rlap{\drawvector(#1,#2){\sarrowlength}}%
      \hskip\@hchannel\hbox{\drawvector(-#1,-#2){\sarrowlength}}}%
  \d@@blearrtrue\generalsmap(#1,#2)}%
\newif\ifrtm@rph
\def\@shiftmorph#1{\hbox{\setbox0=\hbox{$\scriptstyle#1$}%
  \setbox1=\hbox{\hskip\@hm@rphshift\raise\@vm@rphshift\copy0}%
  \wd1=\wd0 \ht1=\ht0 \dp1=\dp0 \box1}}%
\def\@hm@rphshift{\ifrtm@rph
  \ifdim\hmorphposnrt=\z@\hmorphposn\else\hmorphposnrt\fi \else
  \ifdim\hmorphposnlft=\z@\hmorphposn\else\hmorphposnlft\fi \fi}%
\def\@vm@rphshift{\ifrtm@rph
  \ifdim\vmorphposnrt=\z@\vmorphposn\else\vmorphposnrt\fi \else
  \ifdim\vmorphposnlft=\z@\vmorphposn\else\vmorphposnlft\fi \fi}%
\def\addm@rph{\ifx\next\lft\let\temp=\lftmorph\else
  \ifx\next\rt\let\temp=\rtmorph\else\let\temp\relax\fi\fi \temp}%
\def\plnmorph{\dimen1\wd\@arrbox \ifdim\dimen1<\z@ \dimen1-\dimen1\fi
  \vcenter{\box\@arrbox}}%
\def\lftmorph\lft#1{\rtm@rphfalse \setbox0=\@shiftmorph{#1}%
  \if@pslope \let\@upordown\raise \else \let\@upordown\lower\fi
  \llap{\@upordown\@vmorphdflt\hbox to\dimen1{\hss 
    \llap{\box0}\hss}\hskip\@hmorphdflt}\futurelet\next\addm@rph}%
\def\rtmorph\rt#1{\rtm@rphtrue \setbox0=\@shiftmorph{#1}%
  \if@pslope \let\@upordown\lower \else \let\@upordown\raise\fi
  \llap{\@upordown\@vmorphdflt\hbox to\dimen1{\hss
    \rlap{\box0}\hss}\hskip-\@hmorphdflt}\futurelet\next\addm@rph}%
\def\getm@rphposn(#1,#2){\ifd@@blearr \dimen@\morphdist \advance\dimen@ by
  .5\channelwidth \@getshift(#1,#2){\@hmorphdflt}{\@vmorphdflt}{\dimen@}\else
  \@getshift(#1,#2){\@hmorphdflt}{\@vmorphdflt}{\morphdist}\fi}%
\def\getch@nnel(#1,#2){\ifdim\hchannel=\z@ \ifdim\vchannel=\z@
    \@getshift(#1,#2){\@hchannel}{\@vchannel}{\channelwidth}%
    \else \@hchannel\hchannel \@vchannel\vchannel \fi
  \else \@hchannel\hchannel \@vchannel\vchannel \fi}%
\def\@getshift(#1,#2)#3#4#5{\dimen@ #5\relax
  \&xarg #1\relax \&yarg #2\relax
  \ifnum\&xarg<0 \&xarg -\&xarg \fi
  \ifnum\&yarg<0 \&yarg -\&yarg \fi
  \ifnum\&xarg<\&yarg \&negargtrue \&yyarg\&xarg \&xarg\&yarg \&yarg\&yyarg\fi
  \ifcase\&xarg \or  
    \ifcase\&yarg    
      \dimen@i \z@ \dimen@ii \dimen@ \or 
      \dimen@i .7071\dimen@ \dimen@ii .7071\dimen@ \fi \or
    \ifcase\&yarg    
      \or 
      \dimen@i .4472\dimen@ \dimen@ii .8944\dimen@ \fi \or
    \ifcase\&yarg    
      \or 
      \dimen@i .3162\dimen@ \dimen@ii .9486\dimen@ \or
      \dimen@i .5547\dimen@ \dimen@ii .8321\dimen@ \fi \or
    \ifcase\&yarg    
      \or 
      \dimen@i .2425\dimen@ \dimen@ii .9701\dimen@ \or\or
      \dimen@i .6\dimen@ \dimen@ii .8\dimen@ \fi \or
    \ifcase\&yarg    
      \or 
      \dimen@i .1961\dimen@ \dimen@ii .9801\dimen@ \or
      \dimen@i .3714\dimen@ \dimen@ii .9284\dimen@ \or
      \dimen@i .5144\dimen@ \dimen@ii .8575\dimen@ \or
      \dimen@i .6247\dimen@ \dimen@ii .7801\dimen@ \fi \or
    \ifcase\&yarg    
      \or 
      \dimen@i .1645\dimen@ \dimen@ii .9864\dimen@ \or\or\or\or
      \dimen@i .6402\dimen@ \dimen@ii .7682\dimen@ \fi \fi
  \if&negarg \&tempdima\dimen@i \dimen@i\dimen@ii \dimen@ii\&tempdima\fi
  #3\dimen@i\relax #4\dimen@ii\relax }%
\catcode`\&=4  
}%
\catcode`& = 4
\toks2 = {%
\catcode`\&=4  
\def\generalhmap{\futurelet\next\@generalhmap}%
\def\@generalhmap{\ifx\next^ \let\temp\generalhm@rph\else
  \ifx\next_ \let\temp\generalhm@rph\else \let\temp\m@kehmap\fi\fi \temp}%
\def\generalhm@rph#1#2{\ifx#1^
    \toks@=\expandafter{\the\toks@#1{\rtm@rphtrue\@shiftmorph{#2}}}\else
    \toks@=\expandafter{\the\toks@#1{\rtm@rphfalse\@shiftmorph{#2}}}\fi
  \generalhmap}%
\def\m@kehmap{\mathrel{\smash@@{\the\toks@}}}%
\def\mapright{\toks@={\mathop{\vcenter{\smash@@{\drawrightarrow}}}\limits}%
  \generalhmap}%
\def\mapleft{\toks@={\mathop{\vcenter{\smash@@{\drawleftarrow}}}\limits}%
  \generalhmap}%
\def\bimapright{\toks@={\mathop{\vcenter{\smash@@{\drawbirightarrow}}}\limits}%
  \generalhmap}%
\def\bimapleft{\toks@={\mathop{\vcenter{\smash@@{\drawbileftarrow}}}\limits}%
  \generalhmap}%
\def\adjmapright{\toks@={\mathop{\vcenter{\smash@@{\drawadjrightarrow}}}\limits}%
  \generalhmap}%
\def\adjmapleft{\toks@={\mathop{\vcenter{\smash@@{\drawadjleftarrow}}}\limits}%
  \generalhmap}%
\def\hline{\toks@={\mathop{\vcenter{\smash@@{\drawhline}}}\limits}%
  \generalhmap}%
\def\bihline{\toks@={\mathop{\vcenter{\smash@@{\drawbihline}}}\limits}%
  \generalhmap}%
\def\drawrightarrow{\hbox{\drawvector(1,0){\harrowlength}}}%
\def\drawleftarrow{\hbox{\drawvector(-1,0){\harrowlength}}}%
\def\drawbirightarrow{\hbox{\raise.5\channelwidth
  \hbox{\drawvector(1,0){\harrowlength}}\lower.5\channelwidth
  \llap{\drawvector(1,0){\harrowlength}}}}%
\def\drawbileftarrow{\hbox{\raise.5\channelwidth
  \hbox{\drawvector(-1,0){\harrowlength}}\lower.5\channelwidth
  \llap{\drawvector(-1,0){\harrowlength}}}}%
\def\drawadjrightarrow{\hbox{\raise.5\channelwidth
  \hbox{\drawvector(-1,0){\harrowlength}}\lower.5\channelwidth
  \llap{\drawvector(1,0){\harrowlength}}}}%
\def\drawadjleftarrow{\hbox{\raise.5\channelwidth
  \hbox{\drawvector(1,0){\harrowlength}}\lower.5\channelwidth
  \llap{\drawvector(-1,0){\harrowlength}}}}%
\def\drawhline{\hbox{\drawline(1,0){\harrowlength}}}%
\def\drawbihline{\hbox{\raise.5\channelwidth
  \hbox{\drawline(1,0){\harrowlength}}\lower.5\channelwidth
  \llap{\drawline(1,0){\harrowlength}}}}%
\def\generalvmap{\futurelet\next\@generalvmap}%
\def\@generalvmap{\ifx\next\lft \let\temp\generalvm@rph\else
  \ifx\next\rt \let\temp\generalvm@rph\else \let\temp\m@kevmap\fi\fi \temp}%
\toksdef\toks@@=1
\def\generalvm@rph#1#2{\ifx#1\rt 
    \toks@=\expandafter{\the\toks@
      \rlap{$\vcenter{\rtm@rphtrue\@shiftmorph{#2}}$}}\else 
    \toks@@={\llap{$\vcenter{\rtm@rphfalse\@shiftmorph{#2}}$}}%
    \toks@=\expandafter\expandafter\expandafter{\expandafter\the\expandafter
      \toks@@ \the\toks@}\fi \generalvmap}%
\def\m@kevmap{\the\toks@}%
\def\mapdown{\toks@={\vcenter{\drawdownarrow}}\generalvmap}%
\def\mapup{\toks@={\vcenter{\drawuparrow}}\generalvmap}%
\def\bimapdown{\toks@={\vcenter{\drawbidownarrow}}\generalvmap}%
\def\bimapup{\toks@={\vcenter{\drawbiuparrow}}\generalvmap}%
\def\adjmapdown{\toks@={\vcenter{\drawadjdownarrow}}\generalvmap}%
\def\adjmapup{\toks@={\vcenter{\drawadjuparrow}}\generalvmap}%
\def\vline{\toks@={\vcenter{\drawvline}}\generalvmap}%
\def\bivline{\toks@={\vcenter{\drawbivline}}\generalvmap}%
\def\drawdownarrow{\hbox to5pt{\hss\drawvector(0,-1){\varrowlength}\hss}}%
\def\drawuparrow{\hbox to5pt{\hss\drawvector(0,1){\varrowlength}\hss}}%
\def\drawbidownarrow{\hbox to5pt{\hss\hbox{\drawvector(0,-1){\varrowlength}}%
  \hskip\channelwidth\hbox{\drawvector(0,-1){\varrowlength}}\hss}}%
\def\drawbiuparrow{\hbox to5pt{\hss\hbox{\drawvector(0,1){\varrowlength}}%
  \hskip\channelwidth\hbox{\drawvector(0,1){\varrowlength}}\hss}}%
\def\drawadjdownarrow{\hbox to5pt{\hss\hbox{\drawvector(0,-1){\varrowlength}}%
  \hskip\channelwidth\lower\varrowlength
  \hbox{\drawvector(0,1){\varrowlength}}\hss}}%
\def\drawadjuparrow{\hbox to5pt{\hss\hbox{\drawvector(0,1){\varrowlength}}%
  \hskip\channelwidth\raise\varrowlength
  \hbox{\drawvector(0,-1){\varrowlength}}\hss}}%
\def\drawvline{\hbox to5pt{\hss\drawline(0,1){\varrowlength}\hss}}%
\def\drawbivline{\hbox to5pt{\hss\hbox{\drawline(0,1){\varrowlength}}%
  \hskip\channelwidth\hbox{\drawline(0,1){\varrowlength}}\hss}}%
\def\commdiag#1{\null\,
  \vcenter{\commdiagbaselines
  \m@th\ialign{\hfil$##$\hfil&&\hfil$\mkern4mu ##$\hfil\crcr
      \mathstrut\crcr\noalign{\kern-\baselineskip}
      #1\crcr\mathstrut\crcr\noalign{\kern-\baselineskip}}}\,}%
\def\commdiagbaselines{\baselineskip15pt \lineskip3pt \lineskiplimit3pt }%
\def\gridcommdiag#1{\null\,
  \vcenter{\offinterlineskip
  \m@th\ialign{&\vbox to\vgrid{\vss
    \hbox to\hgrid{\hss\smash@@{$##$}\hss}}\crcr
      \mathstrut\crcr\noalign{\kern-\vgrid}
      #1\crcr\mathstrut\crcr\noalign{\kern-.5\vgrid}}}\,}%
\newdimen\harrowlength \harrowlength=60pt
\newdimen\varrowlength \varrowlength=.618\harrowlength
\newdimen\sarrowlength \sarrowlength=\harrowlength
\newdimen\hmorphposn \hmorphposn=\z@
\newdimen\vmorphposn \vmorphposn=\z@
\newdimen\morphdist  \morphdist=4pt
\dimendef\@hmorphdflt 0       
\dimendef\@vmorphdflt 2       
\newdimen\hmorphposnrt  \hmorphposnrt=\z@
\newdimen\hmorphposnlft \hmorphposnlft=\z@
\newdimen\vmorphposnrt  \vmorphposnrt=\z@
\newdimen\vmorphposnlft \vmorphposnlft=\z@

\newdimen\hgrid \hgrid=15pt
\newdimen\vgrid \vgrid=15pt
\newdimen\hchannel  \hchannel=0pt
\newdimen\vchannel  \vchannel=0pt
\newdimen\channelwidth \channelwidth=3pt
\dimendef\@hchannel 0         
\dimendef\@vchannel 2         
\catcode`& = \@oldandcatcode
\catcode`@ = \@oldatcatcode
}%
\let\newif = \@plainnewif
\let\newdimen = \@plainnewdimen
\ifx\noarrow\@undefined \the\toks0 \the\toks2 \fi
\catcode`& = \@eplainoldandcode
\def\environment#1{%
   \ifx\@groupname\@undefined\else
      \errhelp = \@unnamedendgrouphelp
      \errmessage{`\@groupname' was not closed by \string\endenvironment}%
   \fi
   \edef\@groupname{#1}%
   \begingroup
      \let\@groupname = \@undefined
}%
\def\endenvironment#1{%
   \endgroup
   \edef\@thearg{#1}%
   \ifx\@groupname\@thearg
   \else
      \ifx\@groupname\@undefined
         \errhelp = \@isolatedendenvironmenthelp
         \errmessage{Isolated \string\endenvironment\space for `#1'}%
      \else
         \errhelp = \@mismatchedenvironmenthelp
         \errmessage{Environment `#1' ended, but `\@groupname' started}%
         \endgroup 
      \fi
   \fi
   \let\@groupname = \@undefined
}%
\newhelp\@unnamedendgrouphelp{Most likely, you just forgot an^^J%
   \string\endenvironment.  Maybe you should try inserting another^^J%
   \string\endgroup to recover.}%
\newhelp\@isolatedendenvironmenthelp{You ended an environment X, but^^J%
   no \string\environment{X} to start it is anywhere in sight.^^J%
   You might also be at an \string\endenvironment\space that would match^^J%
   a \string\begingroup, i.e., you forgot an \string\endgroup.}%
\newhelp\@mismatchedenvironmenthelp{You started an environment named X, but^^J%
   you ended one named Y.  Maybe you made a typo in one^^J%
   or the other of the names?}%
\newif\ifenvironment
\def\checkenv{\ifenvironment \errhelp = \@interwovenenvhelp
   \errmessage{Interwoven environments}%
   \egroup \fi
}%
\newhelp\@interwovenenvhelp{Perhaps you forgot to end the previous^^J%
   environment? I'm finishing off the current group,^^J%
   hoping that will fix it.}%
\newtoks\previouseverydisplay
\let\@leftleftfill\relax 
\newdimen\leftdisplayindent \leftdisplayindent=\parindent
\newif\if@leftdisplays
\def\leftdisplays{%
  \if@leftdisplays\else
    \previouseverydisplay = \everydisplay
    \everydisplay = {\the\previouseverydisplay \leftdisplaysetup}%
    \let\@save@maybedisableeqno = \@maybedisableeqno
    \let\@saveeqno = \eqno
    \let\@saveleqno = \leqno
    \let\@saveeqalignno = \eqalignno
    \let\@saveleqalignno = \leqalignno
    \let\@maybedisableeqno = \relax
    \def\eqno{\hfill\textstyle\enspace}%
    \def\leqno{%
      \hfill
      \hbox to0pt\bgroup
        \kern-\displaywidth
        \kern-\leftdisplayindent    
        $\aftergroup\@leftleqnoend  
    }%
    \@redefinealignmentdisplays
    \@leftdisplaystrue
  \fi
}%
\newbox\@lignbox
\newdimen\disprevdepth
\def\centereddisplays{%
  \if@leftdisplays
    \everydisplay = \previouseverydisplay
    \let\@maybedisableeqno = \@save@maybedisableeqno
    \let\eqno = \@saveeqno
    \let\leqno = \@saveleqno
    \let\eqalignno = \@saveeqalignno
    \let\leqalignno = \@saveleqalignno
    \@leftdisplaysfalse
  \fi
}%
\def\leftdisplaysetup{%
   \dimen@ = \leftdisplayindent
   \advance\dimen@ by \leftskip
   \advance\displayindent by \dimen@
   \advance\displaywidth by -\dimen@
   \halign\bgroup##\cr \noalign\bgroup
      \disprevdepth = \prevdepth
      \setbox\z@ = \hbox to\displaywidth\bgroup
      $\displaystyle
      \aftergroup\@lefteqend 
}
\def\@lefteqend{
   \hfil\egroup
   \@putdisplay}
\def\@leftleqnoend{\hss \egroup $}
\def\@putdisplay{%
   \ifvoid\@lignbox 
     \moveright\displayindent\box\z@ 
   \else 
     \prevdepth = \dp\@lignbox 
     \unvbox\@lignbox
   \fi
   \egroup\egroup 
   $
}
\def\@redefinealignmentdisplays{%
  \def\displaylines##1{
    \global\setbox\@lignbox\vbox{%
      \prevdepth = \disprevdepth
      \displ@y
      \tabskip\displayindent
      \halign{\hbox to\displaywidth
        {$\@lign\displaystyle####\hfil$\hfil}\crcr
              ##1\crcr}}}%
  \def\eqalignno##1{%
    \def\eqno{&}%
    \def\leqno{&}%
    \global\setbox\@lignbox\vbox{%
      \prevdepth = \disprevdepth
      \displ@y
      \advance\displaywidth by \displayindent
      \tabskip\displayindent
      \halign to\displaywidth{%
         \hfil $\@lign\displaystyle{####}$\@leftleftfill\tabskip\z@skip
        &$\@lign\displaystyle{{}####}$\hfil\tabskip\centering
        &\llap{$\@lign####$}\tabskip\z@skip\crcr
        ##1\crcr}}}%
  \def\leqalignno##1{%
    \def\eqno{&}%
    \def\leqno{&}%
    \global\setbox\@lignbox\vbox{%
      \prevdepth = \disprevdepth
      \displ@y
      \advance\displaywidth by \displayindent
      \tabskip\displayindent
      \halign to\displaywidth{%
         \hfil $\@lign\displaystyle{####}$\@leftleftfill\tabskip\z@skip
        &$\@lign\displaystyle{{}####}$\hfil\tabskip\centering
        &\kern-\displaywidth 
         \rlap{\kern\displayindent \kern-\leftdisplayindent$\@lign####$}%
         \tabskip\displaywidth\crcr
        ##1\crcr}}}%
}%
\let\@primitivenoalign = \noalign
\newtoks\@everynoalign
\def\@lefteqalignonoalign#1{%
  \@primitivenoalign{%
    \advance\leftskip by -\parindent
    \advance\leftskip by -\leftdisplayindent
    \parskip = 0pt
    \parindent = 0pt
    \the\@everynoalign
    #1%
  }%
}%
\def\monthname{%
   \ifcase\month
      \or Jan\or Feb\or Mar\or Apr\or May\or Jun%
      \or Jul\or Aug\or Sep\or Oct\or Nov\or Dec%
   \fi
}%
\def\fullmonthname{%
   \ifcase\month
      \or January\or February\or March\or April\or May\or June%
      \or July\or August\or September\or October\or November\or December%
   \fi
}%
\def\timestring{\begingroup
   \count0 = \time
   \divide\count0 by 60
   \count2 = \count0   
   \count4 = \time
   \multiply\count0 by 60
   \advance\count4 by -\count0   
   \ifnum\count4<10
      \toks1 = {0}%
   \else
      \toks1 = {}%
   \fi
   \ifnum\count2<12
      \toks0 = {a.m.}%
   \else
      \toks0 = {p.m.}%
      \advance\count2 by -12
   \fi
   \ifnum\count2=0
      \count2 = 12
   \fi
   \number\count2:\the\toks1 \number\count4 \thinspace \the\toks0
\endgroup}%
\def\today{\the\day\ \fullmonthname\ \the\year}%
\newskip\abovelistskipamount      \abovelistskipamount = .5\baselineskip
  \newcount\abovelistpenalty      \abovelistpenalty    = 10000
  \def\abovelistskip{\vpenalty\abovelistpenalty \vskip\abovelistskipamount}%
\newskip\interitemskipamount      \interitemskipamount = 0pt
  \newcount\belowlistpenalty      \belowlistpenalty    = -50
\newskip\belowlistskipamount      \belowlistskipamount = .5\baselineskip
  \newcount\interitempenalty      \interitempenalty    = 0
  \def\interitemskip{\vpenalty\interitempenalty \vskip\interitemskipamount}%
\newdimen\listleftindent    \listleftindent = 0pt
\newdimen\listrightindent   \listrightindent = 0pt        
\let\listmarkerspace = \enspace
\newtoks\everylist
\newdimen\@listindent
\def\beginlist{%
  \abovelistskip
  \@listindent = \parindent
  \advance\@listindent by \listleftindent
  \advance\leftskip by \@listindent
  \advance\rightskip by \listrightindent
  \itemnumber = 1
  \the\everylist
}%
\def\li{\@getoptionalarg\@finli}%
\def\@finli{%
  \let\@lioptarg\@optionalarg
  \ifx\@lioptarg\empty \else
    \begingroup
      \@@hldestoff
      \expandafter\writeitemxref\expandafter{\@lioptarg}%
    \endgroup
  \fi
  \ifnum\itemnumber=1 \else \interitemskip \fi
  \begingroup
    \ifx\@lioptarg\empty \else
      \toks@ = \expandafter{\@lioptarg}%
      \let\li@nohldest@marker\marker
      \edef\marker{\noexpand\hldest@impl{li}{\the\toks@}\noexpand\li@nohldest@marker}%
    \fi
    \printitem
  \endgroup
  \advance\itemnumber by 1
  \advance\itemletter by 1
  \advance\itemromannumeral by 1
  \ignorespaces
}%
\def\writeitemxref#1{\definexref{#1}\marker{item}}%
\def\printitem{%
  \par
  \nobreak
  \vskip-\parskip
  \noindent
  \printmarker\marker
}%
\def\printmarker#1{\llap{\marker \enspace}}%
\newcount\numberedlistdepth
\newcount\itemnumber
\newcount\itemletter
\newcount\itemromannumeral
\def\numberedmarker{%
  \ifcase\numberedlistdepth
      (impossible)%
  \or \printitemnumber
  \or \printitemletter
  \or \printitemromannumeral
  \else *%
  \fi
}%
\def\printitemnumber{\number\itemnumber}%
\def\printitemletter{\char\the\itemletter}%
\def\printitemromannumeral{\romannumeral\itemromannumeral}%
\def\numberedprintmarker#1{\llap{#1) \listmarkerspace}}%
\def\numberedlist{\environment{@numbered-list}%
  \advance\numberedlistdepth by 1
  \itemletter = `a
  \itemromannumeral = 1
  \beginlist
  \let\marker = \numberedmarker
  \let\printmarker = \numberedprintmarker
}%

\newcount\unorderedlistdepth
\def\unorderedmarker{%
  \ifcase\unorderedlistdepth
      (impossible)%
  \or \blackbox
  \or ---%
  \else *%
  \fi
}%
\def\unorderedprintmarker#1{\llap{#1\listmarkerspace}}%
\def\unorderedlist{\environment{@unordered-list}%
  \advance\unorderedlistdepth by 1
  \beginlist
  \let\marker = \unorderedmarker
  \let\printmarker = \unorderedprintmarker
}%
\def\listing#1{%
   \par \begingroup
   \@setuplisting
   \setuplistinghook
   \input #1
   \endgroup
}%
\let\setuplistinghook = \relax
\def\linenumberedlisting{%
  \ifx\lineno\undefined \innernewcount\lineno \fi
  \lineno = 0
  \everypar = {\advance\lineno by 1 \printlistinglineno}%
}%
\def\printlistinglineno{\llap{[\the\lineno]\quad}}%
\def\nolastlinelisting{\aftergroup\@removeboxes}%
\def\@removeboxes{%
  \setbox0 = \lastbox
  \ifvoid0
    \ignorespaces 
  \else
    \expandafter\@removeboxes
  \fi
}%
{%
  \makeactive\^^L
  \let^^L = \relax
  \gdef\@setuplisting{%
     \uncatcodespecials
     \obeywhitespace
     \makeactive\`
     \makeactive\^^I
     \makeactive\^^L
     \def^^L{\vfill\break}%
     \parskip = 0pt
     \listingfont
  }%
}%
\def\listingfont{\tt}%
{%
   \makeactive\`
   \gdef`{\relax\lq}
}%
{%
   \makeactive\^^I
   \gdef^^I{\hskip8\fontdimen2}%
}%
\def\verbatimescapechar#1{%
  \gdef\@makeverbatimescapechar{%
    \@makeverbatimdoubleescape #1%
    \catcode`#1 = 0
  }%
}%
\def\@makeverbatimdoubleescape#1{%
  \catcode`#1 = \other
  \begingroup
    \lccode`\* = `#1%
    \lowercase{\endgroup \ece\def*{*}}%
}%
\verbatimescapechar\|  
\def\verbatim{\begingroup
  \uncatcodespecials
  \makeactive\` 
  \@makeverbatimescapechar
  \tt\obeywhitespace}

\def\definecontentsfile#1{%
  \ece\innernewwrite{#1file}%
  \ece\innernewif{if@#1fileopened}%
  \ece\let{#1filebasename} = \jobname
  \ece\def{open#1file}{\opencontentsfile{#1}}%
  \ece\def{write#1entry}{\writecontentsentry{#1}}%
  \ece\def{writenumbered#1entry}{\writenumberedcontentsentry{#1}}%
  \ece\def{writenumbered#1line}{\writenumberedcontentsline{#1}}%
  \ece\innernewif{ifrewrite#1file} \csname rewrite#1filetrue\endcsname
  \ece\def{read#1file}{\readcontentsfile{#1}}%
}%
\definecontentsfile{toc}%
\def\opencontentsfile#1{%
  \csname if@#1fileopened\endcsname \else
     \ece{\immediate\openout}{#1file} = \csname #1filebasename\endcsname.#1
     \ece\global{@#1fileopenedtrue}%
  \fi
}%
\def\writecontentsentry#1#2#3{\writenumberedcontentsentry{#1}{#2}{#3}{}}%
\def\writenumberedcontentsentry#1#2#3#4{%
  \csname ifrewrite#1file\endcsname
    \writenumberedcontents@cmdname{#1}{#2}%
    \def\temp{#3}
    \toks2 = \expandafter{#4}%
    \edef\cs{\the\toks2}%
    \edef\@wr{%
      \write\csname #1file\endcsname{%
        \the\toks0 
        {\sanitize\temp}
        \ifx\empty\cs\else {\sanitize\cs}\fi 
        {\noexpand\folio}
      }%
    }%
    \@wr
  \fi
  \ignorespaces
}%
\def\writenumberedcontentsline#1#2#3#4{%
  \csname ifrewrite#1file\endcsname
    \writenumberedcontents@cmdname{#1}{#2}%
    \def\temp{#4}
    \toks2 = \expandafter{#3}%
    \edef\cs{\the\toks2}%
    \edef\@wr{%
      \write\csname #1file\endcsname{%
        \the\toks0 
        \ifx\empty\cs\else {\sanitize\cs}\fi 
        {\sanitize\temp}
        {\noexpand\folio}
      }%
    }%
    \@wr
  \fi
  \ignorespaces
}%
\def\writenumberedcontents@cmdname#1#2{%
  \csname open#1file\endcsname
  \edef\temp{#2}
  \expandafter\if\expandafter\isinteger\expandafter{\temp}%
    \toks0 = {\expandafter\noexpand \csname #1entry\endcsname}%
    \edef\temp{\the\toks0{\temp}}%
    \toks0 = \expandafter{\temp}%
  \else
    \toks0 = {\expandafter\noexpand \csname #1#2entry\endcsname}%
  \fi
}%
\def\readcontentsfile#1{%
   \edef\temp{%
     \noexpand\testfileexistence[\csname #1filebasename\endcsname]{#1}%
   }\temp
   \if@fileexists
      \input \csname #1filebasename\endcsname.#1\relax
   \fi
}%
\let\ifxrefwarning = \iftrue
\def\xrefwarningtrue{\@citewarningtrue \let\ifxrefwarning = \iftrue}%
\def\xrefwarningfalse{\@citewarningfalse \let\ifxrefwarning = \iffalse}%
\begingroup
  \catcode`\_ = 8
  \gdef\xrlabel#1{#1_x}%
\endgroup
\def\xrdef#1{%
  \begingroup
    \hldest@impl{xrdef}{#1}%
    \begingroup
      \@@hldestoff
      \definexref{#1}{\noexpand\folio}{page}%
    \endgroup
  \endgroup
  \ignorespaces
}%
\def\definexref#1#2#3{%
  \hldest@impl{definexref}{#1}%
  \edef\temp{#1}%
  \readauxfile
  \edef\@wr{\noexpand\writeaux{\string\@definelabel{\temp}{#2}{#3}}}%
  \@wr
  \ignorespaces
}%
\def\@definelabel#1{
  \begingroup 
    \expandafter\ifx\csname\xrlabel{#1}\endcsname \relax
      \expandafter\@definelabel@nocheck
    \else
      \expandafter\@definelabel@warn
    \fi
    {#1}%
}%
\def\@definelabel@nocheck#1#2#3{%
    \expandafter\gdef\csname\xrlabel{#1}\endcsname{#2}%
    \setpropertyglobal{\xrlabel{#1}}{class}{#3}%
  \endgroup 
}%
\def\@definelabel@warn#1#2#3{%
  \message{^^J\linenumber Label `#1' multiply defined,
           value `#2' of class `#3' overwriting value
           `\csname\xrlabel{#1}\endcsname' of class
           `\getproperty{\xrlabel{#1}}{class}'.}%
  \@definelabel@nocheck{#1}{#2}{#3}%
}%
\def\reftie{\penalty\@M \ }
\let\refspace\ 
\def\xrefn{\@getoptionalarg\@finxrefn}%
\def\@finxrefn#1{%
  \hlstart@impl{ref}{#1}%
  \ifx\@optionalarg\empty \else
    \let\@xrefnoptarg\@optionalarg
    \readauxfile
    {\@@hloff\@xrefnoptarg}\reftie
  \fi
  \plain@xrefn{#1}%
  \hlend@impl{ref}%
}%
\def\plain@xrefn#1{%
  \readauxfile
  \expandafter \ifx\csname\xrlabel{#1}\endcsname\relax
    \if@citewarning
       \message{\linenumber Undefined label `#1'.}%
    \fi
    \expandafter\def\csname\xrlabel{#1}\endcsname{%
      `{\tt
        \escapechar = -1
        \expandafter\string\csname#1\endcsname
      }'%
    }%
  \fi
  \csname\xrlabel{#1}\endcsname 
}%

\def\xrefpageword{p.\thinspace}%
\def\xref{\@getoptionalarg\@finxref}%
\def\@finxref#1{%
  \hlstart@impl{xref}{#1}%
  \ifx\@optionalarg\empty \else
    {\@@hloff\@optionalarg}\refspace
  \fi
  \xrefpageword\plain@xrefn{#1}%
  \hlend@impl{xref}%
}%
\def\@maybewarnref{%
  \ifundefined{amsppt.sty}%
  \else
    \message{Warning: amsppt.sty and Eplain both define \string\ref. See
             the Eplain manual.}%
    \let\amsref = \ref
  \fi
  \let\ref = \eplainref
  \ref
}
\let\ref = \@maybewarnref
\def\eplainref{\@getoptionalarg\@fineplainref}%
\def\@fineplainref{\@generalref{1}{}}%
\def\refs{\let\@optionalarg\empty \@generalref{0}s}%
\def\@generalref#1#2#3{%
  \let\@generalrefoptarg\@optionalarg
  \readauxfile
  \ifcase#1 \else \hlstart@impl{ref}{#3}\fi
  \edef\@generalref@class{\getproperty{\xrlabel{#3}}{class}}%
  \expandafter\ifx\csname \@generalref@class word\endcsname\relax
    \ifx\@generalrefoptarg\empty \else {\@@hloff\@generalrefoptarg\reftie}\fi
  \else
    \begingroup
      \@@hloff
      \ifx\@generalrefoptarg\empty \else \@generalrefoptarg \refspace \fi
      \csname \@generalref@class word\endcsname
      #2\reftie
    \endgroup
  \fi
  \ifcase#1 \hlstart@impl{ref}{#3}\fi
  \plain@xrefn{#3}%
  \hlend@impl{ref}%
}%
\newcount\eqnumber
\newcount\subeqnumber
\def\eqdefn{\@getoptionalarg\@fineqdefn}%
\def\@fineqdefn#1{%
  \ifx\@optionalarg\empty
    \global\advance\eqnumber by 1
    \def\temp{\eqconstruct{\number\eqnumber}}%
  \else
    \def\temp{\@optionalarg}%
  \fi
  \global\subeqnumber = 0
  \gdef\@currenteqlabel{#1}%
  \toks0 = \expandafter{\@currenteqlabel}%
  \begingroup
    \def\eqrefn{\noexpand\plain@xrefn}%
    \def\xrefn{\noexpand\plain@xrefn}%
    \edef\temp{\noexpand\@eqdefn{\the\toks0}{\temp}}%
    \temp
  \endgroup
}%
\def\eqsubdefn#1{%
  \global\advance\subeqnumber by 1
  \toks0 = {#1}%
  \toks2 = \expandafter{\@currenteqlabel}%
  \begingroup
    \def\eqrefn{\noexpand\plain@xrefn}%
    \def\xrefn{\noexpand\plain@xrefn}%
    \def\eqsubreftext{\noexpand\eqsubreftext}%
    \edef\temp{%
      \noexpand\@eqdefn
        {\the\toks0}%
        {\eqsubreftext{\eqrefn{\the\toks2}}{\the\subeqnumber}}%
    }%
    \temp           
  \endgroup
}%
\newcount\phantomeqnumber
\def\phantomeqlabel{PHEQ\the\phantomeqnumber}%
\def\@eqdefn#1#2{%
  \ifempty{#1}%
    \global\advance\phantomeqnumber by 1
    \edef\hl@eqlabel{\phantomeqlabel}%
    \readauxfile
  \else
    \def\hl@eqlabel{#1}%
    {\@@hldestoff \definexref{#1}{#2}{eq}}%
  \fi
  \hldest@impl{eq}{\hl@eqlabel}%
  \begingroup 
    \@definelabel@nocheck{#1}{#2}{eq}%
}%
\def\eqdef{\@getoptionalarg\@fineqdef}%
\def\@fineqdef{%
  \toks0 = \expandafter{\@optionalarg}%
  \edef\temp{\noexpand\@eqdef{\noexpand\eqdefn[\the\toks0]}}%
  \temp
}%
\def\eqsubdef{\@eqdef\eqsubdefn}%
\def\@eqdef#1#2{%
  \@maybedisableeqno
  \eqnum #1{#2}
        \let\@optionalarg\empty 
        {\@@hloff\@fineqref{#2}}
  \@mayberestoreeqno
  \ignorespaces
}%
\let\@mayberestoreeqno = \relax
\def\@maybedisableeqno{%
  \ifinner
    \global\let\eqno = \relax
    \global\let\leqno = \relax
    \global\let\@mayberestoreeqno = \@restoreeqno
  \fi
}%
\let\@primitiveeqno = \eqno
\let\@primitiveleqno = \leqno
\def\@restoreeqno{%
  \global\let\eqno = \@primitiveeqno
  \global\let\leqno = \@primitiveleqno
  \global\let\@mayberestoreeqno = \empty
}%
\def\righteqnumbers{%
  \def\eqnum{\eqno}%
  \def\eqalignnum{\eqalignno}%
}%
\def\lefteqnumbers{%
  \def\eqnum{\leqno}%
  \def\eqalignnum{\leqalignno}%
}%
\righteqnumbers
\def\eqrefn{\@getoptionalarg\@fineqrefn}%
\def\@fineqrefn#1{%
  \eqref@start{#1}%
  \plain@xrefn{#1}%
  \hlend@impl{eq}%
}%
\def\eqref{\@getoptionalarg\@fineqref}%
\def\@fineqref#1{%
  \eqref@start{#1}%
  \eqprint{\plain@xrefn{#1}}%
  \hlend@impl{eq}%
}%
\def\eqref@start#1{%
  \let\@eqrefoptarg\@optionalarg
  \ifempty{#1}%
    \hlstart@impl{eq}{\phantomeqlabel}%
  \else
    \hlstart@impl{eq}{#1}%
  \fi
  \ifx\@eqrefoptarg\empty \else
    {\@@hloff\@eqrefoptarg\reftie}%
  \fi
}%
\let\eqconstruct = \identity
\def\eqprint#1{(#1)}%
\def\eqsubreftext#1#2{#1.#2}%
\let\extraidxcmdsuffixes = \empty
\outer\def\defineindex#1{%
  \def\@idxprefix{#1}%
  \expandafter\innernewif\csname if\@idxprefix dx\endcsname
  \csname \@idxprefix dxtrue\endcsname
  \for\@idxcmd:=,marked,submarked,name%
                \extraidxcmdsuffixes\do
  {%
    \@defineindexcmd\@idxcmd
  }%
  \ece\innernewwrite{@#1indexfile}%
  \ece\innernewif{if@#1indexfileopened}%
}%
\newif\ifsilentindexentry
\def\@defineindexcmd#1{%
  \@defineoneindexcmd{s}{#1}\silentindexentrytrue
  \@defineoneindexcmd{}{#1}\silentindexentryfalse
}%
\def\@defineoneindexcmd#1#2#3{%
  \toks@ = {#3}%
  \edef\temp{%
    \def
      \expandonce\csname#1\@idxprefix dx#2\endcsname 
      {\def\noexpand\@idxprefix{\@idxprefix}
       \expandonce\csname @@#1idx#2\endcsname
      }%
    \def
      \expandonce\csname @@#1idx#2\endcsname{
        \the\toks@
        \noexpand\@idxgetrange\expandonce\csname @#1idx#2\endcsname
      }%
  }%
  \temp
}%
\let\indexfilebasename = \jobname
\def\@idxwrite#1#2{%
  \csname if\@idxprefix dx\endcsname
    \@openidxfile
    \def\temp{#1}%
    \edef\@wr{%
      \expandafter\write\csname @\@idxprefix indexfile\endcsname{%
        \string\indexentry
        {\sanitize\temp}%
        {\noexpand#2}%
      }%
    }%
    \@wr
  \else
    \write-1{}%
  \fi
  \ifindexproofing
    \def\temp{#1}%
    \edef\temp{%
      \insert\@indexproof{\noexpand\indexproofterm{\sanitize\temp}}%
    }%
    \temp
    \ifhmode\allowhyphens\fi
  \fi
  \hookrun{afterindexterm}%
  \ifsilentindexentry \expandafter\ignorespaces\fi
}%
\def\@openidxfile{%
  \csname if@\@idxprefix indexfileopened\endcsname \else
    \expandafter\immediate\openout\csname @\@idxprefix indexfile\endcsname =
      \indexfilebasename.\@idxprefix dx
    \expandafter\global\csname @\@idxprefix indexfileopenedtrue\endcsname
  \fi
}%
\newif\ifindexproofing
\newinsert\@indexproof
\dimen\@indexproof = \maxdimen                  
\count\@indexproof = 0  \skip\@indexproof = 0pt 
\font\indexprooffont = cmtt8
\def\indexproofterm#1{\hbox{\strut \indexprooffont #1}}%
\let\@plainmakeheadline = \makeheadline
\def\makeheadline{%
  \expandafter\ifx\csname\idxpageanchor{\folio}\endcsname\relax \else
    {\@@hldeston \hldest@impl{idx}{\hlidxpagelabel{\folio}}}%
  \fi
  \indexproofunbox
  \@plainmakeheadline
}%
\def\indexsetmargins{%
  \ifx\undefined\outsidemargin
    \dimen@ = 1truein
    \advance\dimen@ by \hoffset
    \edef\outsidemargin{\the\dimen@}%
    \let\insidemargin = \outsidemargin
  \fi
}%
\def\indexproofunbox{%
  \ifvoid\@indexproof\else
    \indexsetmargins
    \rlap{%
      \kern\hsize
      \ifodd\pageno \kern\outsidemargin \else \kern\insidemargin \fi
      \vbox to 0pt{\unvbox\@indexproof\vss}%
    }\nointerlineskip
  \fi
}%
\def\idxrangebeginword{begin}%
\def\idxbeginrangemark{(}
\def\idxrangeendword{end}%
\def\idxendrangemark{)}%
\def\idxseecmdword{see}%
\def\idxseealsocmdword{seealso}%
\newif\if@idxsee
\newif\if@hlidxencap
\let\@idxseenterm = \relax
\def\idxpagemarkupcmdword{pagemarkup}%
\let\@idxpagemarkup = \relax
\def\@idxgetrange#1{%
  \let\@idxrangestr = \empty
  \let\@afteridxgetrange = #1%
  \begingroup
    \catcode\idxargopen=1
    \@getoptionalarg\@finidxgetopt
}%
\def\@finidxgetopt{%
    \global\let\@idxgetrange@arg\@optionalarg
  \endgroup
  \@hlidxencaptrue
  \for\@idxarg:=\@idxgetrange@arg\do{%
    \expandafter\@idxcheckpagemarkup\@idxarg=,%
    \ifx\@idxarg\idxrangebeginword
      \def\@idxrangestr{\idxencapoperator\idxbeginrangemark}%
    \else
      \ifx\@idxarg\idxrangeendword
        \def\@idxrangestr{\idxencapoperator\idxendrangemark}%
        \@hlidxencapfalse
      \else
        \ifx\@idxarg\idxseecmdword
          \def\@idxpagemarkup{indexsee}%
          \@idxseetrue
          \@hlidxencapfalse
        \else
          \ifx\@idxarg\idxseealsocmdword
            \def\@idxpagemarkup{indexseealso}%
            \@idxseetrue
            \@hlidxencapfalse
          \else
             \ifx\@idxpagemarkup\relax
               \errmessage{Unrecognized index option `\@idxarg'}%
             \fi
          \fi
        \fi
      \fi
    \fi
  }%
  \ifnum\hldest@place@idx < 0 \else
    \if@hlidxencap
      \ifx\@idxpagemarkup\relax
        \let\@idxpagemarkup\empty
      \fi
      \ifcase\hldest@place@idx \relax
        \edef\@idxpagemarkup{hlidxpage{\@idxpagemarkup}}%
        \definepageanchor{\noexpand\folio}%
      \else
        \global\advance\hlidxlabelnumber by 1
        \edef\@idxpagemarkup{hlidx{\hlidxlabel}{\@idxpagemarkup}}%
        \hldest@impl{idx}{\hlidxlabel}%
      \fi
    \fi
  \fi
  \@afteridxgetrange
}%
\def\@idxcheckpagemarkup#1=#2,{%
  \def\temp{#1}%
  \ifx\temp\idxpagemarkupcmdword
    \if ,#2, 
      \errmessage{Missing markup command to `pagemarkup'}%
    \else
      \def\temp##1={##1}%
      \edef\@idxpagemarkup{\temp\string#2}%
    \fi
  \fi
}%
\def\hldest@place@idx{-1}%
\begingroup
  \catcode`\_ = 8
  \gdef\idxpageanchor#1{#1_p}%
\endgroup
\def\definepageanchor#1{%
  \readauxfile
  \edef\@wr{\noexpand\writeaux{\string\@definepageanchor{#1}}}%
  \@wr
  \ignorespaces
}%
\def\@definepageanchor#1{%
  \expandafter\gdef\csname\idxpageanchor{#1}\endcsname{}%
}%
\newcount\hlidxlabelnumber
\def\hlidxlabel{IDX\the\hlidxlabelnumber}%
\def\hlidxpagelabel#1{IDXPG#1}%
\def\hlidx#1#2#3{%
  \ifempty{#2}%
    \let\@idxpageencap\relax
  \else
    \expandafter\let\expandafter\@idxpageencap\csname #2\endcsname
  \fi
  \hlstart@impl{idx}{#1}%
  \@idxpageencap{#3}%
  \hlend@impl{idx}%
}%
\def\hlidxpage#1#2{%
  \@hlidxgetpages{#2}%
  \ifempty{#1}%
    \let\@idxpageencap\relax
  \else
    \expandafter\let\expandafter\@idxpageencap\csname #1\endcsname
  \fi
  \hlstart@impl{idx}{\hlidxpagelabel{\@idxpageiref}}%
  \expandafter\@idxpageencap\expandafter{\@idxpagei}%
  \hlend@impl{idx}%
  \ifx\@idxpageii\empty \else
    \@idxpagesep
    \hlstart@impl{idx}{\hlidxpagelabel{\@idxpageiiref}}%
    \expandafter\@idxpageencap\expandafter{\@idxpageii}%
    \hlend@impl{idx}%
  \fi
}%
\def\@hlidxgetpages#1{%
  \idxparselist{#1}%
  \ifx\idxpagei\empty
    \idxparserange{#1}%
    \ifx\idxpagei\empty
      \def\@idxpageiref{#1}
    \else
      \let\@idxpageiref\idxpagei 
    \fi
    \def\@idxpagei{#1}%
    \let\@idxpageii\empty
  \else
    \let\@idxpagei\idxpagei
    \let\@idxpageii\idxpageii
    \let\@idxpageiref\idxpagei 
    \let\@idxpageiiref\idxpageii 
    \let\@idxpagesep\idxpagelistdelimiter
  \fi
}%
\def\setidxpagelistdelimiter#1{%
  \gdef\idxpagelistdelimiter{#1}%
  \gdef\@removepagelistdelimiter##1#1{##1}%
  \gdef\@idxparselist##1#1##2\@@{%
    \ifempty{##2}%
      \let\idxpagei\empty
    \else
      \def\idxpagei{##1}%
      \edef\idxpageii{\@removepagelistdelimiter##2}%
    \fi
  }%
  \gdef\idxparselist##1{\@idxparselist##1#1\@@}%
}%
\def\setidxpagerangedelimiter#1{%
  \gdef\idxpagerangedelimiter{#1}%
  \gdef\@removepagerangedelimiter##1#1{##1}%
  \gdef\@idxparserange##1#1##2\@@{%
    \ifempty{##2}%
      \let\idxpagei\empty
    \else
      \def\idxpagei{##1}%
      \edef\idxpageii{\@removepagerangedelimiter##2}%
    \fi
  }%
  \gdef\idxparserange##1{\@idxparserange##1#1\@@}%
}%
\setidxpagelistdelimiter{, }%
\setidxpagerangedelimiter{--}%
\def\idxsubentryseparator{!}%
\def\idxencapoperator{|}%
\def\idxmaxpagenum{99999}%
\newtoks\@idxmaintoks
\newtoks\@idxsubtoks
\def\@idxtokscollect{%
  \edef\temp{\the\@idxsubtoks}%
  \edef\@indexentry{%
    \the\@idxmaintoks
    \ifx\temp\empty\else \idxsubentryseparator\the\@idxsubtoks \fi
    \@idxrangestr
  }%
  \if@idxsee
    \@idxseefalse 
    \edef\temp{\noexpand\idx@getverbatimarg
      {\noexpand\@finidxtokscollect{\idxmaxpagenum}}}%
  \else
    \def\temp{\@finfinidxtokscollect\folio}%
  \fi
  \temp
}%
\def\@finidxtokscollect#1#2{%
  \def\@idxseenterm{#2}%
  \@finfinidxtokscollect{#1}%
}%
\def\@finfinidxtokscollect#1{%
  \ifx\@idxpagemarkup\relax \else
    \toks@ = \expandafter{\@indexentry}%
    \edef\@indexentry{%
      \the\toks@
      \ifx\@idxrangestr\empty \idxencapoperator \fi
      \@idxpagemarkup
    }%
    \let\@idxpagemarkup = \relax
  \fi
  \ifx\@idxseenterm\relax \else
    \toks@ = \expandafter{\@indexentry}%
    \edef\@indexentry{\the\toks@{\sanitize\@idxseenterm}}%
    \let\@idxseenterm = \relax
  \fi
  \expandafter\@idxwrite\expandafter{\@indexentry}{#1}%
}%
\def\@idxcollect#1#2{%
  \@idxmaintoks = {#1}%
  \@idxsubtoks = {#2}%
  \@idxtokscollect
}%
\def\idxargopen{`\{}%
\def\idxargclose{`\}}%
\def\idx@getverbatimarg#1{%
  \def\idx@getverbatimarg@cmd{#1}
  \begingroup
    \uncatcodespecials
    \catcode\idxargopen=1
    \catcode\idxargclose=2
    \catcode`\ =10   
    \catcode`\^^I=10 
    \catcode`\^^M=5  
    \idx@fingetverbatimarg
}%
\def\idx@fingetverbatimarg#1{\endgroup\idx@getverbatimarg@cmd{#1}}%
\def\idx@getverboptarg#1{%
  \def\idx@optionaltemp{#1}
  \let\idx@optionalnext = \relax
  \begingroup
    \if@idxsee \catcode\idxargopen=1 \fi
    \@futurenonspacelet\idx@optionalnext\idx@bracketcheck
}%
\def\idx@bracketcheck{%
  \ifx [\idx@optionalnext
    \endgroup 
    \expandafter\idx@finbracketcheck
  \else
    \endgroup 
    \let\@optionalarg = \empty
    \expandafter\idx@optionaltemp
  \fi
}%
\def\idx@finbracketcheck{%
  \begingroup
    \uncatcodespecials
    \catcode`\ =10   
    \catcode`\^^I=10 
    \catcode`\^^M=5  
    \idx@@getoptionalarg
}%
\def\idx@@getoptionalarg[#1]{%
  \endgroup
  \def\@optionalarg{#1}%
  \idx@optionaltemp
}%
{\catcode`\&=0
\gdef\idx@scanterm#1{%
  \edef\idx@scanterm@nl@catcode{\the\catcode`\^^M\relax}%
  \catcode`\^^M=5
  \scantokens{#1&endinput}%
  \catcode`\^^M=\idx@scanterm@nl@catcode
}}%
\def\@idx{\idx@getverbatimarg\@finidx}%
\def\@finidx#1{%
  \idx@scanterm{#1}
  \@idxcollect{#1}{}%
}%
\def\@sidx{\idx@getverbatimarg\@finsidx}%
\def\@finsidx#1{\@idxmaintoks = {#1}\idx@getverboptarg\@finfinsidx}%
\def\@finfinsidx{%
  \@idxsubtoks = \expandafter{\@optionalarg}%
  \@idxtokscollect
}%
\def\idxsortkeysep{@}
\def\@idxconstructmarked#1#2#3{%
  \toks@ = {#2}
  \toks2 = {#3}
  \edef\temp{\the\toks2 \idxsortkeysep \the\toks@{\the\toks2}}%
  #1 = \expandafter{\temp}%
}%
\def\@idxmarked#1{\idx@getverbatimarg{\@finidxmarked{#1}}}%
\def\@finidxmarked#1#2{%
  \idx@scanterm{#1{#2}}
  \@idxconstructmarked\@idxmaintoks{#1}{#2}%
  \@idxsubtoks = {}%
  \@idxtokscollect
}%
\def\@sidxmarked#1{\idx@getverbatimarg{\@finsidxmarked{#1}}}%
\def\@finsidxmarked#1#2{%
  \@idxconstructmarked\toks@{#1}{#2}%
  \edef\temp{{\the\toks@}}%
  \expandafter\@finsidx\temp
}%
\def\@idxsubmarked{%
  \let\sidxsubmarked@print\idxsubmarked@print
  \idx@getverbatimarg\@finsidxsubmarked
}%
\def\idxsubmarked@print{\expandafter\idx@scanterm\expandafter}%
\def\@sidxsubmarked{%
  \let\sidxsubmarked@print\gobble
  \idx@getverbatimarg\@finsidxsubmarked
}%
\def\@finsidxsubmarked#1{\@idxmaintoks = {#1}\@@finsidxsubmarked}
\def\@@finsidxsubmarked#1{\idx@getverbatimarg{\@@@finsidxsubmarked{#1}}}
\def\@@@finsidxsubmarked#1#2{
  \sidxsubmarked@print 
    {\the\@idxmaintoks\space #1{#2}}
  \@idxconstructmarked\@idxsubtoks{#1}{#2}%
  \@idxtokscollect
}%
\def\idxnameseparator{, }
\def\@idxcollectname#1#2{%
  \def\temp{#1}%
  \ifx\temp\empty
    \toks@ = {}%
  \else
    \toks@ = \expandafter{\idxnameseparator #1}%
  \fi
  \toks2 = {#2}%
  \edef\temp{\the\toks2 \the\toks@}%
}%
\def\@idxname{\idx@getverbatimarg\@finidxname}%
\def\@finidxname#1{\idx@getverbatimarg{\@finfinidxname{#1}}}%
\def\@finfinidxname#1#2{%
  \idx@scanterm{#1 #2}
  \@idxcollectname{#1}{#2}%
  \expandafter\@idxcollect\expandafter{\temp}{}%
}%
\def\@sidxname{\idx@getverbatimarg\@finsidxname}%
\def\@finsidxname#1{\idx@getverbatimarg{\@finfinsidxname{#1}}}%
\def\@finfinsidxname#1#2{%
  \@idxcollectname{#1}{#2}%
  \expandafter\@finsidx\expandafter{\temp}%
}%
\let\indexfonts = \relax
\def\readindexfile#1{%
  \edef\@idxprefix{#1}%
  \testfileexistence[\indexfilebasename]{\@idxprefix nd}%
  \iffileexists \begingroup
    \ifx\begin\undefined
      \def\begin##1{\@beginindex}%
      \let\end = \@gobble
    \fi
    \input \indexfilebasename.\@idxprefix nd
    \singlecolumn
  \endgroup
  \else
    \message{No index file \indexfilebasename.\@idxprefix nd.}%
  \fi
}%
\def\@beginindex{%
  \let\item = \@indexitem
  \let\subitem = \@indexsubitem
  \let\subsubitem = \@indexsubsubitem
  \indexfonts
  \doublecolumns
  \parindent = 0pt
  \hookrun{beginindex}%
}%

\newskip\aboveindexitemskipamount  \aboveindexitemskipamount = 0pt plus2pt
\def\aboveindexitemskip{\vskip\aboveindexitemskipamount}%
\def\@indexitem{\begingroup
  \@indexitemsetup
  \leftskip = 0pt
  \aboveindexitemskip
  \penalty-100 
  \def\par{\endgraf\endgroup\nobreak}%
}%
\def\@indexsubitem{%
  \@indexitemsetup
  \leftskip = 1em
}%
\def\@indexsubsubitem{%
  \@indexitemsetup
  \leftskip = 2em
}%
\def\@indexitemsetup{%
  \par
  \hangindent = 1em
  \raggedright
  \hyphenpenalty = 10000
  \hookrun{indexitem}%
}%
\def\seevariant{\it}%
\def\indexseeword{see}%
\def\indexsee{\idx@getverbatimarg\@finindexsee}%
\def\@finindexsee#1#2{{\seevariant \indexseeword\/ }\idx@scanterm{#1}}%
\def\indexseealsowords{see also}%
\def\indexseealso{\idx@getverbatimarg\@finindexseealso}%
\def\@finindexseealso#1#2{{\seevariant \indexseealsowords\/ }\idx@scanterm{#1}}%
\defineindex{i}%
\begingroup
  \catcode `\^^M = \active %
  \gdef\flushleft{%
    \def\@endjustifycmd{\@endflushleft}%
    \def\@eoljustifyaction{\null\hfil\break}%
    \let\@firstlinejustifyaction = \relax
    \@startjustify %
  }%
  \gdef\flushright{%
    \def\@endjustifycmd{\@endflushright}%
    \def\@eoljustifyaction{\break\null\hfil}%
    \def\@firstlinejustifyaction{\hfil\null}%
    \@startjustify %
  }%
  \gdef\center{%
    \def\@endjustifycmd{\@endcenter}%
    \def\@eoljustifyaction{\hfil\break\null\hfil}%
    \def\@firstlinejustifyaction{\hfil\null}%
    \@startjustify %
  }%
  \gdef\@startjustify{%
    \parskip = 0pt
    \catcode`\^^M = \active %
    \def^^M{\futurelet\next\@finjustifyreturn}%
    \def\@eateol##1^^M{%
      \def\temp{##1}%
      \@firstlinejustifyaction %
      \ifx\temp\empty\else \temp^^M\fi %
    }%
    \expandafter\aftergroup\@endjustifycmd %
    \checkenv \environmenttrue %
    \par\noindent %
    \@eateol %
  }%
  \gdef\@finjustifyreturn{%
    \@eoljustifyaction %
    \ifx\next^^M%
      \def\par{\endgraf\vskip\blanklineskipamount \global\let\par = \endgraf}%
      \@endjustifycmd %
      \noindent %
      \@firstlinejustifyaction %
    \fi %
  }%
\endgroup
\def\@endflushleft{\unpenalty{\parfillskip = 0pt plus1fil\par}\ignorespaces}%
\def\@endflushright{
   \unskip \setbox0=\lastbox \unpenalty
   {\parfillskip = 0pt \par}\ignorespaces
}%
\def\@endcenter{
   \unskip \setbox0=\lastbox \unpenalty
   {\parfillskip = 0pt plus1fil \par}\ignorespaces
}%
\newcount\abovecolumnspenalty   \abovecolumnspenalty = 10000
\newcount\@linestogo         
\newcount\@linestogoincolumn 
\newcount\@columndepth       
\newdimen\@columnwidth       
\newtoks\crtok  \crtok = {\cr}%
\newcount\currentcolumn
\def\makecolumns#1/#2 {\par \begingroup  
   \@columndepth = #1
   \advance\@columndepth by -1
   \divide \@columndepth by #2
   \advance\@columndepth by 1
   \@linestogoincolumn = \@columndepth
   \@linestogo = #1
   \currentcolumn = 1
   \def\@endcolumnactions{%
      \ifnum \@linestogo<2 
         \the\crtok \egroup \endgroup \par 
      \else
         \global\advance\@linestogo by -1
         \ifnum\@linestogoincolumn<2
            \global\advance\currentcolumn by 1
            \global\@linestogoincolumn = \@columndepth
            \the\crtok
         \else
            &\global\advance\@linestogoincolumn by -1
         \fi
      \fi
   }%
   \makeactive\^^M
   \letreturn \@endcolumnactions
   \@columnwidth = \hsize
     \advance\@columnwidth by -\parindent
     \divide\@columnwidth by #2
   \penalty\abovecolumnspenalty
   \noindent 
   \valign\bgroup
     &\hbox to \@columnwidth{\strut \hsize = \@columnwidth ##\hfil}\cr
}%
\newcount\footnotenumber
\newcount\hlfootlabelnumber
\newdimen\footnotemarkseparation \footnotemarkseparation = .5em
\newskip\interfootnoteskip \interfootnoteskip = 0pt
\newtoks\everyfootnote
\newdimen\footnoterulewidth \footnoterulewidth = 2in
\newdimen\footnoteruleheight \footnoteruleheight = 0.4pt
\newdimen\belowfootnoterulespace \belowfootnoterulespace = 2.6pt
\let\@plainfootnote = \footnote
\let\@plainvfootnote = \vfootnote
\def\hlfootlabel{FOOT\the\hlfootlabelnumber}%
\def\hlfootbacklabel{FOOTB\the\hlfootlabelnumber}%
\def\@eplainfootnote#1{\let\@sf\empty 
  \ifhmode\edef\@sf{\spacefactor\the\spacefactor}\/\fi
  \global\advance\hlfootlabelnumber by 1
  \hlstart@impl{foot}{\hlfootlabel}%
  \hldest@impl{footback}{\hlfootbacklabel}%
  #1%
  \hlend@impl{foot}%
  \@sf\vfootnote{#1}%
}%
\let\footnote\@eplainfootnote
\def\vfootnote#1{\insert\footins\bgroup
  \interlinepenalty\interfootnotelinepenalty
  \splittopskip\ht\strutbox 
  \advance\splittopskip by \interfootnoteskip
  \splitmaxdepth\dp\strutbox
  \floatingpenalty\@MM
  \leftskip\z@skip \rightskip\z@skip \spaceskip\z@skip \xspaceskip\z@skip
  \everypar = {}%
  \parskip = 0pt 
  \ifnum\@numcolumns > 1 \hsize = \@normalhsize \fi
  \the\everyfootnote
  \vskip\interfootnoteskip
  \indent\llap{%
    \hlstart@impl{footback}{\hlfootbacklabel}%
    \hldest@impl{foot}{\hlfootlabel}%
    #1%
    \hlend@impl{footback}%
    \kern\footnotemarkseparation
  }%
  \footstrut\futurelet\next\fo@t
}%
\def\footnoterule{\dimen@ = \footnoteruleheight
  \advance\dimen@ by \belowfootnoterulespace
  \kern-\dimen@
  \hrule width\footnoterulewidth height\footnoteruleheight depth0pt
  \kern\belowfootnoterulespace
  \vskip-\interfootnoteskip
}%
\def\numberedfootnote{%
  \global\advance\footnotenumber by 1
  \@eplainfootnote{{\number\footnotenumber}}
}%
\newdimen\paperheight 
\ifnum\mag=1000
  \paperheight = 11in 
\else
  \paperheight = 11truein 
\fi
\def\topmargin{\afterassignment\@finishtopmargin \dimen@}%
\def\@finishtopmargin{%
  \dimen2 = \voffset		
  \voffset = \dimen@ \advance\voffset by -1truein
  \advance\dimen2 by -\voffset	
  \advance\vsize by \dimen2	
}%
\def\advancetopmargin{%
  \dimen@ = 0pt \afterassignment\@finishadvancetopmargin \advance\dimen@
}%
\def\@finishadvancetopmargin{%
  \advance\voffset by \dimen@
  \advance\vsize by -\dimen@
}%
\def\bottommargin{\afterassignment\@finishbottommargin \dimen@}%
\def\@finishbottommargin{%
  \@computebottommargin		
  \advance\dimen2 by -\dimen@	
  \advance\vsize by \dimen2	
}%
\def\advancebottommargin{%
  \dimen@ = 0pt \afterassignment\@finishadvancebottommargin \advance\dimen@
}%
\def\@finishadvancebottommargin{%
  \advance\vsize by -\dimen@
}%
\def\@computebottommargin{%
  \dimen2 = \paperheight	
  \advance\dimen2 by -\vsize	
  \advance\dimen2 by -\voffset	
  \advance\dimen2 by -1truein	
}%
\newdimen\paperwidth
\ifnum\mag=1000
  \paperwidth = 8.5in 
\else
  \paperwidth = 8.5truein 
\fi
\def\leftmargin{\afterassignment\@finishleftmargin \dimen@}%
\def\@finishleftmargin{%
  \dimen2 = \hoffset		
  \hoffset = \dimen@ \advance\hoffset by -1truein
  \advance\dimen2 by -\hoffset	
  \advance\hsize by \dimen2	
}%
\def\advanceleftmargin{%
  \dimen@ = 0pt \afterassignment\@finishadvanceleftmargin \advance\dimen@
}%
\def\@finishadvanceleftmargin{%
  \advance\hoffset by \dimen@
  \advance\hsize by -\dimen@
}%
\def\rightmargin{\afterassignment\@finishrightmargin \dimen@}%
\def\@finishrightmargin{%
  \@computerightmargin		
  \advance\dimen2 by -\dimen@	
  \advance\hsize by \dimen2	
}%
\def\advancerightmargin{%
  \dimen@ = 0pt \afterassignment\@finishadvancerightmargin \advance\dimen@
}%
\def\@finishadvancerightmargin{%
  \advance\hsize by -\dimen@
}%
\def\@computerightmargin{%
  \dimen2 = \paperwidth		
  \advance\dimen2 by -\hsize	
  \advance\dimen2 by -\hoffset	
  \advance\dimen2 by -1truein	
}%
\let\@plainm@g = \m@g
\def\m@g{\@plainm@g \paperwidth = 8.5 true in \paperheight = 11 true in}%
\newskip\abovecolumnskip \abovecolumnskip = \bigskipamount
\newskip\belowcolumnskip \belowcolumnskip = \bigskipamount
\newdimen\gutter \gutter = 2pc
\newbox\@partialpage
\newdimen\@normalhsize
\newdimen\@normalvsize  
\newtoks\previousoutput
\def\quadcolumns{\@columns4}%
\def\triplecolumns{\@columns3}%
\def\doublecolumns{\@columns2}%
\def\begincolumns#1{\ifcase#1\relax \or \singlecolumn \or \@columns2 \or
                            \@columns3 \or \@columns4 \else \relax \fi}%
\let\@ndcolumns = \relax
\chardef\@numcolumns = 1
\mathchardef\@ejectpartialpenalty = 10141
\chardef\@col@minlines = 3
\chardef\@col@extralines = 3
\newdimen\@col@extraheight
\def\@columns#1{%
  \@ndcolumns
  \global\let\@ndcolumns = \@endcolumns
  \global\chardef\@numcolumns = #1
  \global\previousoutput = \expandafter{\the\output}%
  \global\output = {%
    \ifnum\outputpenalty = -\@ejectpartialpenalty
      \dimen@ = \vsize
      \advance\dimen@ by \@col@minlines\baselineskip
      \global\setbox\@partialpage =
        \vbox  \ifdim \pagetotal > \vsize  to \dimen@  \fi  {%
	  \unvbox255 \unskip
	}%
    \else
      \the\previousoutput
    \fi
  }%
  \vskip \abovecolumnskip
  \vskip \@col@minlines\baselineskip
  \penalty -\@ejectpartialpenalty
  \global\output = {\@columnoutput}%
  \global\@normalhsize = \hsize
  \global\@normalvsize = \vsize
  \count@ = \@numcolumns
  \advance\count@ by -1
  \global\advance\hsize by -\count@\gutter
  \global\divide\hsize by \@numcolumns
  \advance\vsize by -\ht\@partialpage
  \advance\vsize by -\ht\footins
  \ifvoid\footins\else \advance\vsize by -\skip\footins \fi
  \multiply\count\footins by \@numcolumns
  \advance\vsize by -\ht\topins
  \ifvoid\topins\else \advance\vsize by -\skip\topins \fi
  \multiply\count\topins by \@numcolumns
  \global\vsize = \@numcolumns\vsize
  \@col@extraheight=\@col@extralines\baselineskip
  \multiply\@col@extraheight by \@numcolumns
  \global\advance\vsize by \@col@extraheight
}%
\def\gutterbox{\vbox to \dimen0{\vfil\hbox{\hfil}\vfil}}%
\def\@columnsplit{%
  \splittopskip = \topskip
  \splitmaxdepth = \baselineskip
  \begingroup
    \vbadness = 10000
    \global\setbox1 = \vsplit255 to \dimen@  \global\wd1 = \hsize
    \global\setbox3 = \vsplit255 to \dimen@  \global\wd3 = \hsize
    \ifnum\@numcolumns > 2
      \global\setbox5 = \vsplit255 to \dimen@ \global\wd5 = \hsize
    \fi
    \ifnum\@numcolumns > 3
      \global\setbox7 = \vsplit255 to \dimen@ \global\wd7 = \hsize
    \fi
  \endgroup
  \setbox0 = \box255
  \global\setbox255 = \vbox{%
    \unvbox\@partialpage
    \ifcase\@numcolumns \relax\or\relax
      \or \hbox to \@normalhsize{\box1\hfil\gutterbox\hfil\box3}%
      \or \hbox to \@normalhsize{\box1\hfil\gutterbox\hfil\box3%
                                      \hfil\gutterbox\hfil\box5}%
      \or \hbox to \@normalhsize{\box1\hfil\gutterbox\hfil\box3%
                                      \hfil\gutterbox\hfil\box5%
                                      \hfil\gutterbox\hfil\box7}%
    \fi
  }%
  \setbox\@partialpage = \box0
}%
\def\@columnoutput{%
  \dimen@ = \ht255
    \advance\dimen@ by -\@col@extraheight
    \divide\dimen@ by \@numcolumns
  \@columnsplit
  \@recoverclubpenalty 
  \hsize = \@normalhsize 
  \vsize = \@normalvsize
  \the\previousoutput
  \unvbox\@partialpage
  \penalty\outputpenalty
  \global\vsize = \@numcolumns\@normalvsize
  \global\advance\vsize by \@col@extraheight
}%
\def\singlecolumn{%
  \@ndcolumns
  \chardef\@numcolumns = 1
  \vskip\belowcolumnskip
  \nointerlineskip
}%
\newbox\@singlecolumnbox 
\newdimen\column@pagegoal
\newdimen\column@vsize
\def\@endcolumns{%
  \global\let\@ndcolumns = \relax
  \par 
  \column@pagegoal = \pagegoal
  \advance\column@pagegoal by-\@col@extraheight
  \ifdim \pagetotal > \column@pagegoal
    \column@vsize = \column@pagegoal
  \else
    \column@vsize = \pagetotal
  \fi
  \global\output = {\global\setbox1 = \box255}%
  \pagegoal = \pagetotal
  \break                     
  \setbox2 = \box1           
  \global\output = \expandafter{\the\previousoutput}%
  \setbox\@singlecolumnbox = \box\@partialpage
  \@balancecolumns
}%
\def\@balancecolumns{%
  \global\setbox255 = \copy2  
  \dimen@ = \column@vsize
    \divide\dimen@ by \@numcolumns
  \@columnsplit
  \ifvoid\@partialpage
    \global\vsize = \@normalvsize
    \global\hsize = \@normalhsize
    \dump@balanced@columns
    \let\next\relax
  \else
    \advance \column@vsize by \@numcolumns pt
    \test@spill@columns
    \ifspill@columns
      \begingroup
        \vsize = \@normalvsize
        \hsize = \@normalhsize
        \dump@balanced@columns
        \break
        \@recoverclubpenalty
      \endgroup
      \unvbox\@partialpage
      \let\next\@endcolumns
    \else
      \setbox0=\box\@partialpage 
      \let\next\@balancecolumns
    \fi
  \fi
  \next
}%
\def\dump@balanced@columns{%
  \ifvoid\topins\else\topinsert\unvbox\topins\endinsert\fi
  \unvbox\@singlecolumnbox
  \nointerlineskip
  \box255
}%
\newif\ifspill@columns
\def\test@spill@columns{%
  \spill@columnsfalse
  \ifdim \column@vsize > \column@pagegoal
    \ifvoid\footins
      \ifvoid\topins
        \spill@columnstrue
      \fi
    \fi
  \fi
}%
\def\@saveclubpenalty{
  \edef\@recoverclubpenalty{%
     \global\clubpenalty=\the\clubpenalty\relax%
     \global\let\noexpand\@recoverclubpenalty\relax
  }
}%
\let\@recoverclubpenalty\relax
\newdimen\temp@dimen
\def\columnfill{%
  \par
  \dimen@=\pagetotal   
  \temp@dimen = \vsize 
  \divide\temp@dimen by \@numcolumns 
  \loop
    \ifdim \dimen@ > \temp@dimen
      \advance \dimen@ by -\temp@dimen
      \advance \dimen@ by \topskip 
  \repeat
  \advance \temp@dimen by -\dimen@
  \advance \temp@dimen by -\prevdepth
  \@saveclubpenalty 
  \clubpenalty=10000\relax
  \hrule height\temp@dimen width0pt depth0pt\relax
  \nointerlineskip
  \par
  \nointerlineskip
  \penalty0\vfil 
  \relax
}%
\def\@hldest{%
  \def\hl@prefix{hldest}%
  \let\after@hl@getparam\hldest@aftergetparam
  \begingroup
    \hl@getparam
}%
\def\hldest@aftergetparam{%
  \ifvmode
    \hldest@driver
  \else
    \allowhyphens
    \smash{\ifx\hldest@opt@raise\empty \else \raise\hldest@opt@raise\fi
             \hbox{\hldest@driver}}%
    \allowhyphens
  \fi
  \endgroup
}%
\def\@hlstart{%
  \leavevmode
  \def\hl@prefix{hl}%
  \let\after@hl@getparam\hlstart@aftergetparam
  \begingroup
    \hl@getparam
}%
\def\hlstart@aftergetparam{%
  \ifx\color\undefined \else
    \ifx\hl@opt@color\empty \else
      \ifx\hl@opt@colormodel\empty
        \edef\temp{\noexpand\color{\hl@opt@color}}%
      \else
        \edef\temp{\noexpand\color[\hl@opt@colormodel]{\hl@opt@color}}%
      \fi
      \temp
    \fi
  \fi
  \hl@driver
}%
\def\hl@getparam#1#2{
  \edef\@hltype{#1}%
  \ifx\@hltype\empty
    \expandafter\let\expandafter\@hltype
      \csname \hl@prefix @type\endcsname
  \fi
  \expandafter\ifx\csname \hl@prefix @typeh@\@hltype\endcsname \relax
    \errmessage{Invalid hyperlink type `\@hltype'}%
  \fi
  \For\hl@arg:=#2\do{%
    \ifx\hl@arg\empty \else
      \expandafter\hl@set@opt\hl@arg=,%
    \fi
  }%
  \bgroup
    \uncatcodespecials
    \catcode`\{=1 \catcode`\}=2
    \@hl@getparam
}%
\def\@hl@getparam#1{%
  \egroup
  \edef\@hllabel{#1}%
  \after@hl@getparam
  \ignorespaces
}%
\def\hl@set@opt#1=#2,{%
  \expandafter\ifx\csname \hl@prefix @opt@#1\endcsname \relax
    \errmessage{Invalid hyperlink option `#1'}%
  \fi
  \if ,#2, 
    \errmessage{Missing value for option `#1'}%
  \else
    \def\temp##1={##1}%
    \expandafter\edef\csname \hl@prefix @opt@#1\endcsname{\temp#2}%
  \fi
}%
\def\hldest@impl#1{%
  \expandafter\ifcase\csname hldest@on@#1\endcsname
    \relax\expandafter\gobble
  \else
    \toks@=\expandafter{\csname hldest@type@#1\endcsname}%
    \toks@ii=\expandafter{\csname hldest@opts@#1\endcsname}%
    \edef\temp{\noexpand\hldest{\the\toks@}{\the\toks@ii}}%
    \expandafter\temp
  \fi
}%
\def\hlstart@impl#1{%
  \expandafter\ifcase\csname hl@on@#1\endcsname
    \leavevmode\expandafter\gobble
  \else
    \toks@=\expandafter{\csname hl@type@#1\endcsname}%
    \toks@ii=\expandafter{\csname hl@opts@#1\endcsname}%
    \edef\temp{\noexpand\hlstart{\the\toks@}{\the\toks@ii}}%
    \expandafter\temp
  \fi
}%
\def\hlend@impl#1{%
  \expandafter\ifcase\csname hl@on@#1\endcsname
  \else
    \hlend
  \fi
}%
\def\hl@asterisk@word{*}%
\def\hl@opts@word{opts}%
\newif\if@params@override
\def\hldest@groups{definexref,xrdef,li,eq,bib,foot,footback,idx}%
\def\hl@groups{ref,xref,eq,cite,foot,footback,idx,url,hrefint,hrefext}%
\def\hldesttype{%
  \def\hl@prefix{hldest}%
  \def\hl@param{type}%
  \let\hl@all@groups\hldest@groups
  \futurelet\hl@excl\hl@param@read@excl
}%
\def\hldestopts{%
  \def\hl@prefix{hldest}%
  \def\hl@param{opts}%
  \let\hl@all@groups\hldest@groups
  \futurelet\hl@excl\hl@param@read@excl
}%
\def\hltype{%
  \def\hl@prefix{hl}%
  \def\hl@param{type}%
  \let\hl@all@groups\hl@groups
  \futurelet\hl@excl\hl@param@read@excl
}%
\def\hlopts{%
  \def\hl@prefix{hl}%
  \def\hl@param{opts}%
  \let\hl@all@groups\hl@groups
  \futurelet\hl@excl\hl@param@read@excl
}%
\def\hl@param@read@excl{%
  \ifx!\hl@excl
    \let\next\hl@param@read@opt@arg
    \@params@overridetrue
  \else
    \def\next{\hl@param@read@opt@arg{!}}%
    \@params@overridefalse
  \fi
  \next
}%
\def\hl@param@read@opt@arg#1{
  \@getoptionalarg\hl@setparam
}%
\def\@hl@setparam#1{%
  \ifx\@optionalarg\empty
    \hl@setparam@default{#1}
  \else
    \let\hl@do@all@groups\gobble
    \For\hl@group:=\@optionalarg\do{%
      \ifx\hl@group\hl@asterisk@word
        \def\hl@do@all@groups{\let\@optionalarg\hl@all@groups \hl@setparam}%
      \else
        \hl@setparam@group{#1}%
      \fi
    }%
    \hl@do@all@groups{#1}%
  \fi
}%
\def\hl@setparam@group#1{%
  \ifx\hl@group\empty
    \hl@setparam@default{#1}%
  \else
    \expandafter\ifx\csname\hl@prefix @\hl@param @\hl@group\endcsname\relax
      \errmessage{Hyperlink group `\hl@prefix:\hl@param:\hl@group' is not defined}%
    \fi
    \ifx\hl@param\hl@opts@word
      \if@params@override
        \expandafter\let\csname\hl@prefix @\hl@param @\hl@group\endcsname\empty
      \fi
      \hl@update@opts@with@list{#1}
    \else
      \ece\def{\hl@prefix @\hl@param @\hl@group}{#1}%
    \fi
  \fi
}%
\def\hl@setparam@default#1{%
  \ifx\hl@param\hl@opts@word
    \For\hl@opt:=#1\do{%
      \ifx\hl@opt\empty \else
        \expandafter\hl@set@opt\hl@opt=,%
      \fi
    }%
  \else
    \expandafter\ifx\csname\hl@prefix @\hl@param\endcsname\relax
      \message{Default hyperlink parameter `\hl@prefix:\hl@param' is not defined}%
    \fi
    \ece\def{\hl@prefix @\hl@param}{#1}%
  \fi
}%
\def\hl@update@opts@with@list#1{%
  \global\expandafter\let\expandafter\hl@update@new@list
    \csname \hl@prefix @opts@\hl@group\endcsname
  \begingroup
    \For\hl@opt:=#1\do{%
      \hl@update@opts@with@opt
    }%
  \endgroup
  \ece\let{\hl@prefix @opts@\hl@group}\hl@update@new@list
}%
\def\hl@update@opts@with@opt{%
  \global\let\hl@update@old@list\hl@update@new@list
  \global\let\hl@update@new@list\empty
  \global\let\hl@update@new@opt\hl@opt
  \expandafter\hl@parse@opt@key\hl@opt=,%
  \let\hl@update@new@key\hl@update@key
  \global\let\hl@update@comma\empty
  \begingroup
    \for\hl@opt:=\hl@update@old@list\do{%
      \ifx\hl@opt\empty \else 
        \expandafter\hl@parse@opt@key\hl@opt=,%
        \toks@=\expandafter{\hl@update@new@list}%
        \ifx\hl@update@key\hl@update@new@key
          \ifx\hl@update@new@opt\empty \else 
            \toks@ii=\expandafter{\hl@update@new@opt}%
            \xdef\hl@update@new@list{\the\toks@\hl@update@comma\the\toks@ii}%
            \global\let\hl@update@new@opt\empty
            \global\def\hl@update@comma{,}%
          \fi
        \else
          \toks@ii=\expandafter{\hl@opt}%
          \xdef\hl@update@new@list{\the\toks@\hl@update@comma\the\toks@ii}%
          \global\def\hl@update@comma{,}%
        \fi
      \fi
    }%
  \endgroup
  \ifx\hl@update@new@opt\empty \else
    \toks@=\expandafter{\hl@update@new@list}%
    \toks@ii=\expandafter{\hl@update@new@opt}%
    \xdef\hl@update@new@list{\the\toks@\hl@update@comma\the\toks@ii}%
  \fi
}%
\def\hl@parse@opt@key#1=#2,{\def\hl@update@key{#1}}%
\def\hldest@opt@raise{\normalbaselineskip}%
\def\hl@opt@colormodel{cmyk}%
\def\hl@opt@color{0.28,1,1,0.35}%
\def\hldest@on@definexref{0}%
\def\hldest@on@xrdef{0}%
\def\hldest@on@li{0}%
\def\hldest@on@eq{0}
\def\hldest@on@bib{0}
\def\hldest@on@foot{0}
\def\hldest@on@footback{0}
\def\hldest@on@idx{0}
\let\hldest@type@definexref\empty
\let\hldest@type@xrdef\empty
\let\hldest@type@li\empty
\let\hldest@type@eq\empty 
\let\hldest@type@bib\empty 
\let\hldest@type@foot\empty 
\let\hldest@type@footback\empty 
\let\hldest@type@idx\empty 
\let\hldest@opts@definexref\empty
\let\hldest@opts@xrdef\empty
\let\hldest@opts@li\empty
\def\hldest@opts@eq{raise=1.7\normalbaselineskip}
\let\hldest@opts@bib\empty 
\let\hldest@opts@foot\empty 
\let\hldest@opts@footback\empty 
\let\hldest@opts@idx\empty 
\def\hl@on@ref{0}
\def\hl@on@xref{0}%
\def\hl@on@eq{0}
\def\hl@on@cite{0}
\def\hl@on@foot{0}
\def\hl@on@footback{0}
\def\hl@on@idx{0}%
\def\hl@on@url{0}
\def\hl@on@hrefint{0}
\def\hl@on@hrefext{0}
\let\hl@type@ref\empty 
\let\hl@type@xref\empty
\let\hl@type@eq\empty 
\let\hl@type@cite\empty 
\let\hl@type@foot\empty 
\let\hl@type@footback\empty 
\let\hl@type@idx\empty
\let\hl@type@url\empty 
\let\hl@type@hrefint\empty 
\let\hl@type@hrefext\empty 
\let\hl@opts@ref\empty 
\let\hl@opts@xref\empty
\let\hl@opts@eq\empty 
\let\hl@opts@cite\empty 
\let\hl@opts@foot\empty 
\let\hl@opts@footback\empty 
\let\hl@opts@idx\empty
\let\hl@opts@url\empty 
\let\hl@opts@hrefint\empty 
\let\hl@opts@hrefext\empty 
\def\@hlon{\@hlonoff@value@stub{hl}\@@hlon1 }%
\def\@hloff{\@hlonoff@value@stub{hl}\@@hloff0 }%
\def\@hldeston{\@hlonoff@value@stub{hldest}\@@hldeston1 }%
\def\@hldestoff{\@hlonoff@value@stub{hldest}\@@hldestoff0 }%
\def\@hlonoff@value@stub#1#2#3{%
  \def\hl@prefix{#1}%
  \let\hl@on@empty#2%
  \def\hl@value{#3}%
  \expandafter\let\expandafter\hl@all@groups
    \csname \hl@prefix @groups\endcsname
  \@getoptionalarg\@finhlswitch
}%
\def\@finhlswitch{%
  \ifx\@optionalarg\empty
    \hl@on@empty
  \fi
  \let\hl@do@all@groups\relax
  \For\hl@group:=\@optionalarg\do{%
    \ifx\hl@group\hl@asterisk@word
      \let\@optionalarg\hl@all@groups
      \let\hl@do@all@groups\@finhlswitch
    \else
      \ifx\hl@group\empty
        \hl@on@empty
      \else
        \expandafter\ifx\csname\hl@prefix @on@\hl@group\endcsname \relax
          \errmessage{Hyperlink group `\hl@prefix:on:\hl@group'
                      is not defined}%
        \fi
        \ece\edef{\hl@prefix @on@\hl@group}{\hl@value}%
      \fi
    \fi
  }%
  \hl@do@all@groups
}%
\def\@@hlon{%
  \let\hlstart\@hlstart
  \let\hlend\@hlend
}%
\def\@@hloff{%
  \def\hlstart##1##2##3{\leavevmode\ignorespaces}%
  \let\hlend\relax
}%
\def\@@hldeston{%
  \let\hldest\@hldest
}%
\def\@@hldestoff{%
  \def\hldest##1##2##3{\ignorespaces}%
}%
\def\hl@idxexact@word{idxexact}%
\def\hl@idxpage@word{idxpage}%
\def\hl@idxnone@word{idxnone}%
\def\hl@raw@word{raw}%
\def\enablehyperlinks{\@getoptionalarg\@finenablehyperlinks}%
\def\@finenablehyperlinks{%
  \let\hl@selecteddriver\empty
  \def\hldest@place@idx{0}%
  \for\hl@arg:=\@optionalarg\do{%
    \ifx\hl@arg\hl@idxexact@word
      \def\hldest@place@idx{1}%
    \else
      \ifx\hl@arg\hl@idxnone@word
        \def\hldest@place@idx{-1}%
      \else
        \ifx\hl@arg\hl@idxpage@word
          \def\hldest@place@idx{0}%
        \else
          \let\hl@selecteddriver\hl@arg
        \fi
      \fi
    \fi
  }%
  \ifx\hl@selecteddriver\empty
    \ifpdf
      \def\hl@selecteddriver{pdftex}%
      \message{^^JEplain: using `pdftex' hyperlink driver.}%
    \else
      \def\hl@selecteddriver{hypertex}%
      \message{^^JEplain: using `hypertex' hyperlink driver.}%
    \fi
  \else
    \expandafter\ifx\csname hldriver@\hl@selecteddriver\endcsname \relax
      \errmessage{No hyperlink driver `\hl@selecteddriver' available}%
    \fi
  \fi
  \let\hl@setparam\@hl@setparam
  \csname hldriver@\hl@selecteddriver\endcsname
  \def\@finenablehyperlinks{\errmessage{Hyperlink driver `\hl@selecteddriver'
                                        already selected}}%
  \let\hldriver@nolinks\undefined
  \let\hldriver@hypertex\undefined
  \let\hldriver@pdftex \undefined
  \let\hldriver@dvipdfm\undefined
  \let\hloff\@hloff
  \let\hlon\@hlon
  \let\hldestoff\@hldestoff
  \let\hldeston\@hldeston
  \hlon[*,]\hloff[foot,footback]%
  \hldeston[*,]\hldestoff[foot,footback]%
}%
\def\hldriver@nolinks{%
  \def\@hldest##1##2##3{%
    \edef\temp{\write-1{hldest: ##3}}%
    \ifvmode
      \temp
    \else
      \allowhyphens
      \expandafter\smash\expandafter{\temp}%
      \allowhyphens
    \fi
    \ignorespaces
  }%
  \def\@hlstart##1##2##3{%
    \leavevmode
    \begingroup 
    \edef\temp{\write-1{hlstart: ##3}}%
    \temp
    \ignorespaces
  }%
  \def\@hlend{%
    \edef\temp{\write-1{hlend}}%
    \temp
    \endgroup 
  }%
  \let\hl@setparam\gobble
}%
{\catcode`\#=\other
\gdef\hlhash{#}}%
\def\hldriver@hypertex{%
  \def\hldest@type{xyz}%
  \let\hldest@opt@cmd \empty
  \def\hldest@driver{%
    \ifx\@hltype\hl@raw@word
      \csname \hldest@opt@cmd \endcsname
    \else
    \fi
  }%
  \let\hldest@typeh@raw \empty
  \let\hldest@typeh@xyz \empty
  \def\hl@type{name}%
  \ifx\hl@type@url\empty
    \def\hl@type@url{url}%
  \fi
  \ifx\hl@type@hrefext\empty
    \def\hl@type@hrefext{url}%
  \fi
  \let\hl@opt@cmd  \empty
  \let\hl@opt@ext  \empty
  \let\hl@opt@file \empty
  \def\hl@driver{%
    \ifx\@hltype\hl@raw@word
      \csname \hl@opt@cmd \endcsname
    \else
      \def\hlstart@preamble{html:<a href="}%
      \csname hl@typeh@\@hltype\endcsname
    \fi
  }%
  \let\hl@typeh@raw \empty
  \def\hl@typeh@name{\special{\hlstart@preamble \hlhash\@hllabel">}}%
  \def\hl@typeh@filename{%
    \special{%
      \hlstart@preamble
        file:\hl@opt@file\hl@opt@ext
        \ifempty\@hllabel \else \hlhash\@hllabel\fi
      ">%
    }%
  }%
  \def\hl@typeh@url{%
    \special{%
      \hlstart@preamble
        \@hllabel
      ">%
    }%
  }%
  \def\@hlend{\endgroup}
}%
\def\hldriver@pdftex{%
\ifpdf 
  \def\hldest@type{xyz}%
  \let\hldest@opt@width  \empty
  \let\hldest@opt@height \empty
  \let\hldest@opt@depth  \empty
  \let\hldest@opt@zoom   \empty
  \let\hldest@opt@cmd    \empty
  \def\hldest@driver{%
    \ifx\@hltype\hl@raw@word
      \csname \hldest@opt@cmd \endcsname
    \else
      \pdfdest name{\@hllabel}\@hltype
        \csname hldest@typeh@\@hltype\endcsname
    \fi
  }%
  \let\hldest@typeh@raw   \empty
  \let\hldest@typeh@fit   \empty
  \let\hldest@typeh@fith  \empty
  \let\hldest@typeh@fitv  \empty
  \let\hldest@typeh@fitb  \empty
  \let\hldest@typeh@fitbh \empty
  \let\hldest@typeh@fitbv \empty
  \def\hldest@typeh@fitr{%
    \ifx\hldest@opt@width  \empty \else width  \hldest@opt@width  \fi
    \ifx\hldest@opt@height \empty \else height \hldest@opt@height \fi
    \ifx\hldest@opt@depth  \empty \else depth  \hldest@opt@depth  \fi
  }%
  \def\hldest@typeh@xyz{%
    \ifx\hldest@opt@zoom\empty \else zoom \hldest@opt@zoom \fi
  }%
  \def\hl@type{name}%
  \ifx\hl@type@url\empty
    \def\hl@type@url{url}%
  \fi
  \ifx\hl@type@hrefext\empty
    \def\hl@type@hrefext{url}%
  \fi
  \let\hl@opt@width   \empty
  \let\hl@opt@height  \empty
  \let\hl@opt@depth   \empty
  \def\hl@opt@bstyle  {S}%
  \def\hl@opt@bwidth  {1}%
  \let\hl@opt@bcolor  \empty
  \let\hl@opt@hlight  \empty
  \let\hl@opt@bdash   \empty
  \let\hl@opt@pagefit \empty
  \let\hl@opt@cmd     \empty
  \let\hl@opt@file    \empty
  \let\hl@opt@newwin  \empty
  \def\hl@driver{%
    \ifx\@hltype\hl@raw@word
      \csname \hl@opt@cmd \endcsname
    \else
      \let\hl@BSspec\relax 
      \ifx\hl@opt@bstyle \empty
        \ifx\hl@opt@bwidth \empty
          \ifx\hl@opt@bdash \empty
            \let\hl@BSspec\empty 
          \fi
        \fi
      \fi
      \def\hlstart@preamble{%
        \pdfstartlink
          \ifx\hl@opt@width  \empty \else width  \hl@opt@width  \fi
          \ifx\hl@opt@height \empty \else height \hl@opt@height \fi
          \ifx\hl@opt@depth  \empty \else depth  \hl@opt@depth \fi
          attr{%
            \ifx\hl@opt@bcolor\empty\else /C[\hl@opt@bcolor]\fi
            \ifx\hl@opt@hlight\empty\else /H/\hl@opt@hlight\fi
            \ifx\hl@BSspec\relax
              /BS<<%
                /Type/Border%
                \ifx\hl@opt@bstyle\empty\else /S/\hl@opt@bstyle\fi
                \ifx\hl@opt@bwidth\empty\else /W \hl@opt@bwidth\fi
                \ifx\hl@opt@bdash\empty \else /D[\hl@opt@bdash]\fi
              >>%
            \fi
          }%
      }%
      \csname hl@typeh@\@hltype\endcsname
    \fi
  }%
  \let\hl@typeh@raw\empty
  \def\hl@typeh@name{\hlstart@preamble goto name{\@hllabel}}%
  \def\hl@typeh@num{\hlstart@preamble  goto num \@hllabel}%
  \def\hl@typeh@page{%
    \count@=\@hllabel
    \advance\count@ by-1
    \hlstart@preamble
    user{%
      /Subtype/Link%
      /Dest%
        [\the\count@
          \ifx\hl@opt@pagefit\empty/Fit\else\hl@opt@pagefit\fi]%
    }%
  }%
  \def\hl@typeh@filename{\hl@file{(\@hllabel)}}%
  \def\hl@typeh@filepage{%
    \count@=\@hllabel
    \advance\count@ by-1
    \hl@file{%
      [\the\count@ \ifx\hl@opt@pagefit\empty/Fit\else\hl@opt@pagefit\fi]%
    }%
  }%
  \def\hl@file##1{%
    \hlstart@preamble
    user{%
      /Subtype/Link%
      /A<<%
        /Type/Action%
        /S/GoToR%
        /D##1%
        /F(\hl@opt@file)%
        \ifx\hl@opt@newwin\empty \else
          /NewWindow \ifcase\hl@opt@newwin false\else true\fi
        \fi
      >>%
    }%
  }%
  \def\hl@typeh@url{%
    \hlstart@preamble
    user{%
      /Subtype/Link%
      /A<<%
        /Type/Action%
        /S/URI%
        /URI(\@hllabel)%
      >>%
    }%
  }%
  \def\@hlend{\pdfendlink\endgroup}
\else 
  \message{Eplain warning: `pdftex' hyperlink driver: PDF output is^^J
           \space not enabled, falling back on `nolinks' driver.}%
  \hldriver@nolinks
\fi
}%
\def\hldriver@dvipdfm{%
  \def\hldest@type{xyz}%
  \let\hldest@opt@left   \empty
  \let\hldest@opt@top    \empty
  \let\hldest@opt@right  \empty
  \let\hldest@opt@bottom \empty
  \let\hldest@opt@zoom   \empty
  \let\hldest@opt@cmd    \empty
  \def\hldest@driver{%
    \ifx\@hltype\hl@raw@word
      \csname \hldest@opt@cmd \endcsname
    \else
      \def\hldest@preamble{%
        pdf: dest (\@hllabel) [@thispage
      }%
      \csname hldest@typeh@\@hltype\endcsname
    \fi
  }%
  \let\hldest@typeh@raw\empty
  \def\hldest@typeh@fit{%
    \special{\hldest@preamble /Fit]}%
  }%
  \def\hldest@typeh@fith{%
    \special{\hldest@preamble /FitH
      \ifx\hldest@opt@top\empty @ypos \else \hldest@opt@top \fi]}%
  }%
  \def\hldest@typeh@fitv{%
    \special{\hldest@preamble /FitV
      \ifx\hldest@opt@left\empty @xpos \else \hldest@opt@left \fi]}%
  }%
  \def\hldest@typeh@fitb{%
    \special{\hldest@preamble /FitB]}%
  }%
  \def\hldest@typeh@fitbh{%
    \special{\hldest@preamble /FitBH
      \ifx\hldest@opt@top\empty @ypos \else \hldest@opt@top \fi]}%
  }%
  \def\hldest@typeh@fitbv{%
    \special{\hldest@preamble /FitBV
      \ifx\hldest@opt@left\empty @xpos \else \hldest@opt@left \fi]}%
  }%
  \def\hldest@typeh@fitr{%
    \special{\hldest@preamble /FitR
      \ifx\hldest@opt@left\empty @xpos\else\hldest@opt@left\fi\space
      \ifx\hldest@opt@bottom\empty @ypos\else\hldest@opt@bottom\fi\space
      \ifx\hldest@opt@right\empty @xpos\else\hldest@opt@right\fi\space
      \ifx\hldest@opt@top\empty @ypos\else\hldest@opt@top \fi]}%
  }%
  \def\hldest@typeh@xyz{%
    \begingroup
      \ifx\hldest@opt@zoom\empty
        \count1=\z@ \count2=\z@
      \else
        \count2=\hldest@opt@zoom
        \count1=\count2 \divide\count1 by 1000
        \count3=\count1 \multiply\count3 by 1000
        \advance\count2 by -\count3
      \fi
      \special{\hldest@preamble /XYZ
        \ifx\hldest@opt@left\empty @xpos\else\hldest@opt@left\fi\space
        \ifx\hldest@opt@top\empty @ypos\else\hldest@opt@top\fi\space
        \the\count1.\the\count2]}%
    \endgroup
  }%
  \def\hl@type{name}%
  \ifx\hl@type@url\empty
    \def\hl@type@url{url}%
  \fi
  \ifx\hl@type@hrefext\empty
    \def\hl@type@hrefext{url}%
  \fi
  \def\hl@opt@bstyle  {S}%
  \def\hl@opt@bwidth  {1}%
  \let\hl@opt@bcolor  \empty
  \let\hl@opt@hlight  \empty
  \let\hl@opt@bdash   \empty
  \let\hl@opt@pagefit \empty
  \let\hl@opt@cmd     \empty
  \let\hl@opt@file    \empty
  \let\hl@opt@newwin  \empty
  \def\hl@driver{%
    \ifx\@hltype\hl@raw@word
      \csname \hl@opt@cmd \endcsname
    \else
      \let\hl@BSspec\relax 
      \ifx\hl@opt@bstyle \empty
        \ifx\hl@opt@bwidth \empty
          \ifx\hl@opt@bdash \empty
            \let\hl@BSspec\empty 
          \fi
        \fi
      \fi
      \def\hlstart@preamble{%
        pdf: beginann
          <<%
            /Type/Annot%
            /Subtype/Link%
            \ifx\hl@opt@bcolor\empty\else /C[\hl@opt@bcolor]\fi
            \ifx\hl@opt@hlight\empty\else /H/\hl@opt@hlight\fi
            \ifx\hl@BSspec\relax
              /BS<<%
                /Type/Border%
                \ifx\hl@opt@bstyle\empty\else /S/\hl@opt@bstyle\fi
                \ifx\hl@opt@bwidth\empty\else /W \hl@opt@bwidth\fi
                \ifx\hl@opt@bdash\empty \else /D[\hl@opt@bdash]\fi
              >>%
            \fi
      }%
      \csname hl@typeh@\@hltype\endcsname
    \fi
  }%
  \let\hl@typeh@raw\empty
  \def\hl@typeh@name{\special{\hlstart@preamble /Dest(\@hllabel)>>}}%
  \def\hl@typeh@page{%
    \count@=\@hllabel
    \advance\count@ by-1
    \special{%
      \hlstart@preamble
      /Dest[\the\count@
            \ifx\hl@opt@pagefit\empty/Fit\else\hl@opt@pagefit\fi]%
     >>%
    }%
  }%
  \def\hl@typeh@filename{\hl@file{(\@hllabel)}}%
  \def\hl@typeh@filepage{%
    \count@=\@hllabel
    \advance\count@ by-1
    \hl@file{%
      [\the\count@ \ifx\hl@opt@pagefit\empty/Fit\else\hl@opt@pagefit\fi]%
    }%
  }%
  \def\hl@file##1{%
    \special{%
      \hlstart@preamble
      /A<<%
        /Type/Action%
        /S/GoToR%
        /D##1%
        /F(\hl@opt@file)%
        \ifx\hl@opt@newwin\empty \else
          /NewWindow \ifcase\hl@opt@newwin false\else true\fi
        \fi
      >>%
     >>%
    }%
  }%
  \def\hl@typeh@url{%
    \special{%
      \hlstart@preamble
      /A<<%
        /Type/Action%
        /S/URI%
        /URI(\@hllabel)%
      >>%
     >>%
    }%
  }%
  \def\@hlend{\endgroup}
}%
\def\href{%
  \bgroup
    \uncatcodespecials
    \catcode`\{=1 \catcode`\}=2
    \@href
}%
\def\@href#1{
  \egroup
  \edef\@hreftmp{\ifempty{#1}{}\fi}
  \expandafter\@@href\@hreftmp#1\@@
}%
\def\href@end@int{\hlend@impl{hrefint}}%
\def\href@end@ext{\hlend@impl{hrefext}}%
\def\@@href#1#2\@@{%
  \def\@hreftmp{#1}%
  \ifx\@hreftmp\hlhash
    \let\href@end\href@end@int
    \hlstart@impl{hrefint}{#2}%
  \else
    \let\href@end\href@end@ext
    \hlstart@impl{hrefext}{#1#2}%
  \fi
  \@@@href
}%
\def\@@@href{%
  \futurelet\@hreftmp\href@
}%
\def\href@{%
  \ifcat\bgroup\noexpand\@hreftmp
    \let\@hreftmp\href@@
  \else
    \let\@hreftmp\href@@@
  \fi
  \@hreftmp
}%
\def\href@@{\bgroup\aftergroup\href@end \let\@hreftmp}%
\def\href@@@#1{#1\href@end}%
\def\hldeston{\errmessage{Please enable hyperlinks with
  \string\enablehyperlinks\space before using hyperlink commands
  (consider selecting the `nolinks' driver to ignore all hyperlink
  commands in your document)}}%
\let\hldestoff\hldeston \let\hlon\hldeston \let\hloff\hldeston
\let\hlstart\hldeston \let\hlend\hldeston \let\hldest\hldeston
\let\hl@setparam\hldeston
\@hloff[*]\@hldestoff[*]%
\newif\ifusepkg@miniltx@loaded
\newcount\usepkg@recursion@level
\def\usepkg@rcrs{\the\usepkg@recursion@level}%
\let\usepkg@at@begin@document\empty
\let\usepkg@at@end@of@package\empty
\def\usepkg@word@autopict{autopict}%
\def\usepkg@word@psfrag{psfrag}%
\long\def\beginpackages#1\endpackages{%
  \let\usepackage\real@usepackage
  \let\DoNotLoadEpstopdf=t
  \let\eplaininput=\input
  #1%
  \usepkg@at@begin@document
  \global\let\usepkg@at@begin@document\empty
  \global\let\usepackage\fake@usepackage
  \let\packageinput=\input
  \let\input=\eplaininput
}%
\def\fake@usepackage{\errmessage{You should not use \string\usepackage\space outside of^^J
  \@spaces\string\beginpackages...\string\endpackages\space environment}%
}%
\def\eplain@RequirePackage{%
  \global\ece\let{usepkg@save@pkg\usepkg@rcrs}\usepkg@pkg
  \global\ece\let{usepkg@save@options\usepkg@rcrs}\usepkg@options
  \global\ece\let{usepkg@save@date\usepkg@rcrs}\usepkg@date
  \global\ece\let{usepkg@at@end@of@package\usepkg@rcrs}\usepkg@at@end@of@package
  \global\advance\usepkg@recursion@level by\@ne
  \real@usepackage
}%
\let\usepackage\fake@usepackage
\def\real@usepackage{\@getoptionalarg\@finusepackage}%
\def\@finusepackage#1{%
  \let\usepkg@options\@optionalarg
  \ifempty{#1}%
    \errmessage{No packages specified}%
  \fi
  \ifusepkg@miniltx@loaded \else
    \testfileexistence[miniltx]{tex}%
    \if@fileexists
      \input miniltx.tex
      \global\usepkg@miniltx@loadedtrue
      \global\let\RequirePackage\eplain@RequirePackage
      \global\let\DeclareOption\eplain@DeclareOption
      \global\let\PassOptionsToPackage\eplain@PassOptionsToPackage
      \global\let\ExecuteOptions\eplain@ExecuteOptions
      \gdef\ProcessOptions{\@ifstar\eplain@ProcessOptions
                                   \eplain@ProcessOptions}%
      \global\let\AtBeginDocument\eplain@AtBeginDocument
      \global\let\AtEndOfPackage\eplain@AtEndOfPackage
      \global\let\ProvidesFile\eplain@ProvidesFile
      \global\let\ProvidesPackage\eplain@ProvidesPackage
    \else
      \errmessage{miniltx.tex not found, cannot load LaTeX packages}%
    \fi
  \fi
  \@ifnextchar[
    {\@finfinusepackage{#1}}%
    {\@finfinusepackage{#1}[]}%
}%
\def\@finfinusepackage#1[#2]{%
  \edef\usepkg@date{#2}%
  \let\usepkg@load@list\empty
  \for\usepkg@pkg:=#1\do{%
    \toks@=\expandafter{\usepkg@load@list}%
    \edef\usepkg@load@list{%
      \the\toks@
      \noexpand\usepkg@load@pkg{\usepkg@pkg}%
    }%
  }%
  \usepkg@load@list
  \ifnum\usepkg@recursion@level>0
    \global\advance\usepkg@recursion@level by\m@ne
    \expandafter\let\expandafter\usepkg@pkg\csname usepkg@save@pkg\usepkg@rcrs\endcsname
    \expandafter\let\expandafter\usepkg@options\csname usepkg@save@options\usepkg@rcrs\endcsname
    \expandafter\let\expandafter\usepkg@date\csname usepkg@save@date\usepkg@rcrs\endcsname
    \expandafter\let\expandafter\usepkg@at@end@of@package\csname usepkg@at@end@of@package\usepkg@rcrs\endcsname
    \global\ece\let{usepkg@save@pkg\usepkg@rcrs}\undefined
    \global\ece\let{usepkg@save@options\usepkg@rcrs}\undefined
    \global\ece\let{usepkg@save@date\usepkg@rcrs}\undefined
    \global\ece\let{usepkg@at@end@of@package\usepkg@rcrs}\undefined
  \fi
}%
\def\usepkg@load@pkg#1{%
  \def\usepkg@pkg{#1}%
  \ifx\usepkg@pkg\usepkg@word@autopict
    \testfileexistence[picture]{tex}%
    \if@fileexists \else
      \errmessage{Loader `picture.tex' for package `\usepkg@pkg' not found}%
    \fi
  \else
    \ifx\usepkg@pkg\usepkg@word@psfrag
      \testfileexistence[psfrag]{tex}%
      \if@fileexists \else
        \errmessage{Loader `psfrag.tex' for package `\usepkg@pkg' not found}%
      \fi
    \fi
  \fi
  \ifundefined{ver@\usepkg@pkg.sty}%
    \expandafter\@finusepkg@load@pkg
  \else
    \immediate\write-1{^^J\linenumber Eplain: package `\usepkg@pkg' already
             loaded, skipping reloading}%
  \fi
}%
\def\@finusepkg@load@pkg{%
  \testfileexistence[\usepkg@pkg]{sty}%
  \if@fileexists \else
    \errmessage{Package `\usepkg@pkg' not found}%
  \fi
  \expandafter\let\expandafter\temp\csname usepkg@options@\usepkg@pkg\endcsname
  \ifx\temp\relax
    \let\temp\empty
  \fi
  \ifx\temp\empty
    \let\temp\usepkg@options
  \else
    \ifx\usepkg@options\empty \else
      \edef\temp{\temp,\usepkg@options}%
    \fi
  \fi
  \global\ece\let{usepkg@options@\usepkg@pkg}\temp
  \let\usepackage\eplain@RequirePackage
  \global\let\usepkg@at@end@of@package\empty
  \ifx\usepkg@pkg\usepkg@word@autopict
    \input picture.tex
  \else
    \ifx\usepkg@pkg\usepkg@word@psfrag
      \input \usepkg@pkg.tex
    \else
      \input \usepkg@pkg.sty
    \fi
  \fi
  \usepkg@at@end@of@package
  \global\let\usepkg@at@end@of@package\empty
  \let\usepackage\real@usepackage
  \global\ece\let{usepkg@options@\usepkg@pkg}\undefined
  \def\Url@HyperHook##1{\hlstart@impl{url}{\Url@String}##1\hlend@impl{url}}%
}%
\def\eplain@DeclareOption#1#2{%
  \toks@{#2}%
  \expandafter\xdef\csname usepkg@option@\usepkg@pkg @#1\endcsname{\the\toks@}%
}%
\def\eplain@PassOptionsToPackage#1#2{%
  \ifempty{#1}\else
    \for\usepkg@temp:=#2\do{%
      \expandafter\let\expandafter\temp\csname usepkg@options@\usepkg@temp\endcsname
      \ifx\temp\relax
        \let\temp\empty
      \fi
      \ifx\temp\empty
        \edef\temp{#1}%
      \else
        \edef\temp{\temp,#1}%
      \fi
      \global\ece\let{usepkg@options@\usepkg@temp}\temp
    }%
  \fi
}%
\def\usepkg@exec@options#1{%
  \for\CurrentOption:=#1\do{%
    \expandafter\let\expandafter\usepkg@temp
      \csname usepkg@option@\usepkg@pkg @\CurrentOption\endcsname
    \ifx\usepkg@temp\relax
      \expandafter\let\expandafter\temp\csname usepkg@option@\usepkg@pkg @*\endcsname
      \ifx\temp\relax
        \errmessage{Unknown option `\CurrentOption' to package `\usepkg@pkg'}%
      \else
        \temp
      \fi
    \else
      \usepkg@temp
    \fi
  }%
}%
\let\eplain@ExecuteOptions\usepkg@exec@options
\def\eplain@ProcessOptions{%
  \expandafter\usepkg@exec@options\csname usepkg@options@\usepkg@pkg\endcsname
}%
\def\usepkg@accumulate@cmds#1#2{%
  \toks@=\expandafter{#1}%
  \toks@ii={#2}%
  \xdef#1{\the\toks@\the\toks@ii}%
}%
\def\eplain@AtBeginDocument{\usepkg@accumulate@cmds\usepkg@at@begin@document}%
\def\eplain@AtEndOfPackage{\usepkg@accumulate@cmds\usepkg@at@end@of@package}%
\def\eplain@ProvidesPackage#1{%
  \@ifnextchar[
    {\eplain@pr@videpackage{#1.sty}}{\eplain@pr@videpackage#1[]}%
}%
\def\eplain@pr@videpackage#1[#2]{%
  \wlog{#1: #2}%
  \expandafter\xdef\csname ver@#1\endcsname{#2}%
  \@ifl@t@r{#2}\usepkg@date{}%
    {\message{Warning: you have requested package `\usepkg@pkg', version \usepkg@date,^^J
       \@spaces but only version `\csname ver@#1\endcsname' is available.}}%
}%
\def\eplain@ProvidesFile#1{%
  \@ifnextchar[
    {\eplain@pr@videfile{#1}}{\eplain@pr@videfile#1[]}%
}%
\def\eplain@pr@videfile#1[#2]{%
  \wlog{#1: #2}%
  \expandafter\xdef\csname ver@#1\endcsname{#2}%
}%
\def\@ifl@ter#1#2{%
  \expandafter\@ifl@t@r
    \csname ver@#2.#1\endcsname
}%
\def\@ifl@t@r#1#2{%
  \ifnum\expandafter\@parse@version#1//00\@nil<%
        \expandafter\@parse@version#2//00\@nil
    \expandafter\@secondoftwo
  \else
    \expandafter\@firstoftwo
  \fi
}%
\def\@parse@version#1/#2/#3#4#5\@nil{#1#2#3#4 }%

\def\strip@prefix#1>{}%
\def\@ifpackageloaded#1{%
  \expandafter\ifx\csname ver@#1.sty\endcsname\relax
    \expandafter\@secondoftwo
  \else
    \expandafter\@firstoftwo
  \fi
}%
\long\def\g@addto@macro#1#2{%
  \begingroup
    \toks@\expandafter{#1#2}%
    \xdef#1{\the\toks@}%
  \endgroup
}%
\def\PackageWarning#1#2{%
  \begingroup
    \newlinechar=10 %
    \def\MessageBreak{%
      ^^J(#1)\@spaces\@spaces\@spaces\@spaces
    }%
    \immediate\write16{^^JPackage #1 Warning: #2\on@line.^^J}%
  \endgroup
}%
\def\PackageWarningNoLine#1#2{%
  \PackageWarning{#1}{#2\@gobble}%
}%
\def\on@line{ on input line \the\inputlineno}%
\def\@spaces{\space\space\space\space}%
\def\@inmatherr#1{%
   \relax
   \ifmmode
     \errmessage{The command is invalid in math mode}%
   \fi
}%
\let\protected@edef\edef
\let\wlog = \@plainwlog
\catcode`@ = \@eplainoldatcode
\def\eplain{t}%
{\edef\plainversion{\fmtversion}%
 \xdef\fmtversion{3.4: 21 February 2010 (and plain \plainversion)}%
}%

\lefteqnumbers
   \def\testd{oui}
   \def\choixlat{\ifx\numadroite\testd\righteqnumbers
            \else  \lefteqnumbers\fi}
    \choixlat

\catcode`@=\letter
\def\@eplainfootnote#1{\let\@sf\empty 
  \ifhmode\edef\@sf{\spacefactor\the\spacefactor}\/\fi
  \global\advance\hlfootlabelnumber by 1
  \hlstart@impl{foot}{\hlfootlabel}%
  \hldest@impl{footback}{\hlfootbacklabel}%
  \hbox{$^{(#1)}$}%
  \hlend@impl{foot}%
  \@sf\vfootnote{#1.}%
}%
\catcode`@=\other

  \interfootnoteskip=0pt
  
  \everyfootnote={\eightpoint\leftskip=5truemm\rightskip5truemm}
  
  \hsize150truemm\vsize 240truemm\hoffset=5truemm

  \pretolerance=500\tolerance=1000\brokenpenalty=5000
  \parindent3mm
  
  \countdef\temps=170
  \temps=\time
  \countdef\nminutes=171{\nminutes = \time}
  \countdef\nheure=172
  \def\heure{\begingroup                     
     \temps = \time \divide\temps by 60
     \nheure = \temps                        
     \nminutes = \time
     \multiply\temps by 60
     \advance\nminutes by -\temps            
     \ifnum\nminutes<10 \toks1 = {0}%
     \else\toks1 = {}%
     \fi
     \number\nheure h\the\toks1 \number\nminutes  
  \endgroup}%

  \newcount\chstart
  \chstart=\pageno
 \headline={\ifnum\pageno=\chstart {\hfill} \else {\hss \tenrm --\ \folio\ --\hss}\fi}
  \footline={\hfill}
  \normalbaselines
  \frenchspacing
    \def\dater{\vglue-10mm\rightline{(\the\day/\the\month/\the\year)}}
  \def\dateheure{\vglue-10mm\rightline{(\the\day/\the\month/\the\year,\ \heure)}}

  \newif\ifpagetitre \pagetitretrue
\newtoks\hautpagetitre \hautpagetitre={\hfill}
\newtoks\baspagetitre \baspagetitre={\hfill}
\newtoks\auteurcourant \auteurcourant={\hfill}
\newtoks\titrecourant \titrecourant={\hfill}
\newtoks\hautpagegauche
\newtoks\hautpagedroite
\newtoks\hautpagemilieu
\hautpagemilieu={\tenrm\hfil -- \folio\ -- \hfil}
\hautpagegauche={\ifx\midfolio\oui\the\hautpagemilieu\else\tenrm\folio\hfill\the\auteurcourant\hfill\fi}
\hautpagedroite={\ifx\midfolio\oui\the\hautpagemilieu\else\hfill\the\titrecourant\hfill\tenrm\folio\fi}
\newtoks\baspagegauche \baspagegauche={\hfil}
\newtoks\baspagedroite \baspagedroite={\hfil}
\headline={\ifpagetitre\the\hautpagetitre
\else\ifodd\pageno\the\hautpagedroite\else\the\hautpagegauche\fi\fi }
\footline={\ifpagetitre\the\baspagetitre
\else\ifodd\pageno\the\baspagedroite
\else\the\baspagegauche\fi\fi \global\pagetitrefalse}

\def\pageblanche{\vfill\eject\pagetitretrue
\null\vfill\eject
\pagetitretrue
}
\def\chgtpage{\ifodd\pageno \else
\pageblanche \fi \pagetitretrue\titreun=0\footnotenumber=0}

\def\majnombres{\ifodd\pageno \else
\pageblanche \fi \pagetitretrue\hautpoly\titreun=0\footnotenumber=0}

\ifnum\chstart=\pageno \pagetitretrue\fi
  


  \def\leftnote#1{\vadjust{\setbox1=\vtop{\hsize 20mm\parindent=0pt\eightpoint
  \baselineskip=9pt\rightskip=4mm plus 4mm\vskip-4.7mm#1}\hbox{\kern-2cm\smash{\box1}}}}

  
  \def\raggedcenter{\leftskip=20pt plus 10em  
       \rightskip=\leftskip 
        \parfillskip=0pt 
         \spaceskip=.3333em \xspaceskip=.5em 
          \pretolerance=9999 \tolerance=9999
           \hyphenpenalty=9999 \exhyphenpenalty=9999 }


\def\boxit#1#2{\vbox{\hrule\hbox{\vrule\vbox spread#1{\vfil\hbox spread#1{\hfil#2\hfil}\vfil}%
\vrule}\hrule}}


  \def\oui{oui}
  
   \def\fontetitreun{\ifx\paradouze\oui\douzepts\gpdouze\twelvebf\textfont1=\twelveib\else
\quatorzepts\gpquatorze\fourteenbf\fi}

\def\fontetitreunl{\douzepts\textfont1=\twelveib\scriptfont1=\tenib\fourteenti}
 
 \def\fontetitredeux{\textfont1=\eleveni\ifx\paradouze\oui\onzepts\scriptfont1=\ninei\elevenit\else
                        \douzepts\twelveit\fi}
 
   \def\fontetitredeuxb{\ifx\paradouze\oui\onzepts\eleventi\gponze\textfont1=\elevenib\scriptfont1=\nineib
                         \else\douzepts\twelveti\scriptfont1=\twelveib\scriptfont1=\tenib\gpdouze\fi}
                         
\def\fontetitredeuxl{\onzepts\textfont1=\elevenbf\scriptfont1=\ninebf\twelvebf}
  
\def\fontetitretrois{\textfont0=\elevenrm\scriptfont0=\eightrm\textfont1=\eleveni
                      \scriptfont1=\eighti\scriptscriptfont1=\sixi\elevenit}
                      
\def\fontetitrequatre{\textfont0=\elevenrm\scriptfont0=\eightrm\textfont1=\eleveni
                      \scriptfont1=\eighti\scriptscriptfont1=\sixi\elevenrm}
  
  \newcount\titreun\titreun=0
  \newcount\titredeux\titredeux=0
  \newcount\titretrois\titretrois=0
  \newcount\titrequatre\titrequatre=0
  \newcount\enonce\enonce=0
  
  \def\incr#1{\global\advance#1 by 1 {\the #1}}
  \def\avance#1{\global\advance#1 by 1}
  \def\init#1{\global#1=0}
  
  \long\def\Indentation#1#2{\setbox10=\hbox{\fontetitreun#1}
  	                    \ifdim\wd10 < 4mm
                         \setbox10=\hbox to 4mm{\box10\hfill}
                       \else\ifdim\wd10 < 6mm
                         \setbox10=\hbox to 6mm{\box10\hfill}
  	                    \else\ifdim\wd10 < 8mm
                         \setbox10=\hbox to 8mm{\box10\hfill}
                       \else\ifdim\wd10 < 12mm
                         \setbox10=\hbox to 12mm{\box10\hfill}
                       \fi\fi\fi\fi
                       \dimen10=\hsize
                       \advance \dimen10 by -\wd10
                       \noindent \box10 %
                       \ignorespaces
                       \hbox{\vtop{\hsize=\dimen10\raggedright\noindent\fontetitreun#2}}}

  \long\def\paraun#1{\removelastskip\par\medskip\goodbreak\vskip0pt plus.01\vsize\penalty-100
                \vskip0pt plus-.01\vsize
  	              \init{\titredeux}\ifnum\optionparag=1{\init\eqnumber\init\enonce}\else{}\fi
                  \goodbreak{\fontetitreun
  	                \Indentation{\incr{\titreun}.\ }{\fontetitreun #1\par}}\nobreak\medskip}

 %
 %
 \long\def\paraunc#1{\removelastskip\par\bigskip\goodbreak\vskip0pt plus.01\vsize\penalty-100
                \vskip0pt plus-.01\vsize
  	              \init{\titredeux}
                 \ifnum\optionparag=1{\init{\eqnumber}\init\enonce}\else{}\fi
                  \goodbreak
  	                {\parindent0mm\raggedcenter\fontetitreun\incr{\titreun}.\ 
                     \fontetitreun #1\par}\nobreak\medskip}
                     
\newtoks\titreunl
\titreunl={\ifnum\titreun=1{I}\fi%
\ifnum\titreun=2{II}\fi%
\ifnum\titreun=3{III}\fi%
\ifnum\titreun=4{IV}\fi%
\ifnum\titreun=5{V}\fi%
\ifnum\titreun=6{VI}\fi%
\ifnum\titreun=7{VII}\fi%
\ifnum\titreun=8{VIII}\fi%
\ifnum\titreun=9{IX}\fi%
\ifnum\titreun=10{X}\fi%
\ifnum\titreun=11{XI}\fi%
\ifnum\titreun=12{XII}\fi%
\ifnum\titreun=13{XIII}\fi%
}
\long\def\paraunl#1{\removelastskip\par\bigskip\bigskip\goodbreak\vskip0pt plus.01\vsize\penalty-100
                \vskip0pt plus-.01\vsize
  	              \init{\titredeux}\ifnum\optionparag=1{\init\eqnumber\init\enonce}\else{}\fi
                  \goodbreak{\fontetitreunl
  	                \Indentation{\global\advance\titreun by 1{\the\titreunl}.\ }{\fontetitreunl #1\par}}\nobreak\smallskip}

  
  \long\def\paradeux#1{\init{\titretrois}\vskip0pt plus.01\vsize\penalty-10
                \vskip0pt plus-.01\vsize\ifx \elie\oui\medskip\ifnum\titredeux>0\medskip\fi\fi
                 \Indentation{\fontetitredeux\the\titreun${\cdot}$\incr{\titredeux}.}
                              {\fontetitredeux\textfont1=\eleveni#1}\nobreak\par }
  
  \long\def\paradeuxb#1{\init{\titretrois}\vskip0pt plus.001\vsize\penalty-10
                \vskip0pt plus-.01\vsize{\ifx \elie\oui\medskip\ifnum\titredeux>0\medskip\fi\fi
                  \Indentation
  {\fontetitredeuxb\the\titreun${\cdot}$\incr{\titredeux}.}{ \fontetitredeuxb#1}}\nobreak\smallskip
  \par}

\newtoks\titredeuxl
\titredeuxl={\ifnum\titredeux=1{A}\fi%
\ifnum\titredeux=2{B}\fi%
\ifnum\titredeux=3{C}\fi%
\ifnum\titredeux=4{D}\fi%
\ifnum\titredeux=5{E}\fi%
\ifnum\titredeux=6{F}\fi%
\ifnum\titredeux=7{G}\fi%
\ifnum\titredeux=8{H}\fi%
\ifnum\titredeux=9{I}\fi%
\ifnum\titredeux=10{J}\fi%
\ifnum\titredeux=11{K}\fi%
\ifnum\titredeux=12{L}\fi%
\ifnum\titredeux=13{M}\fi%
}
 \long\def\paradeuxl#1{\init{\titretrois}\vskip0pt plus.001\vsize\penalty-10
                \vskip0pt plus-.01
                \vsize \bigskip%
                  \Indentation
     {\fontetitredeuxl\global\advance\titredeux by 1
  \quad \the\titreunl${\cdot}$\the\titredeuxl.}{ \fontetitredeuxl#1}
  \removelastskip\medskip\nobreak\par}
  

  \long\def\paratrois#1{\init{\titrequatre}\ifdim\lastskip<\smallskipamount
                \removelastskip\smallskip\fi
                 \vskip0pt plus.01\vsize\penalty-10
                  \vskip0pt
plus-.01\vsize{\ifx \elie\oui\ifnum\titretrois>0\medskip\fi\fi
\Indentation{\fontetitretrois\the\titreun${\cdot}$\the\titredeux${\cdot}$\incr{\titretrois}.\ }
  {\hskip0mm\baselineskip=14pt\fontetitretrois#1}\smallskip\nobreak\par}}
  
  
  \long\def\paratroisl#1{\init{\titrequatre}\ifdim\lastskip<\smallskipamount
                \removelastskip\fi
                 \vskip0pt plus.01\vsize\penalty-10
                  \vskip0pt
plus-.01\vsize\ifx \elie\oui\bigskip
\fi
\Indentation{\fontetitretrois\quad \quad \the\titreunl{${\cdot}$}\the\titredeuxl${\cdot}$\incr{\titretrois}.\ }
  {\hskip0mm\fontetitretrois#1}\smallskip\nobreak\par}


  \long\def\paraquatre#1{\ifdim\lastskip<\smallskipamount
                \removelastskip\smallskip\fi
                 \vskip0pt plus.01\vsize\penalty-10
                  \vskip0pt
                  plus-.01\vsize\par
 
\Indentation{\fontetitrequatre \the\titreun{${\cdot}$}\the\titredeux${\cdot}$\the\titretrois${\cdot}$\incr{\titrequatre}.\ }
{\hskip0mm\fontetitrequatre#1}\smallskip\nobreak\par}


\newtoks\titrequatrel
\titrequatrel={\ifnum\titrequatre=1{a}\fi%
\ifnum\titrequatre=2{b}\fi%
\ifnum\titrequatre=3{c}\fi%
\ifnum\titrequatre=4{d}\fi%
\ifnum\titrequatre=5{e}\fi%
\ifnum\titrequatre=6{f}\fi%
\ifnum\titrequatre=7{g}\fi%
\ifnum\titrequatre=8{h}\fi%
\ifnum\titrequatre=9{i}\fi%
\ifnum\titrequatre=10{j}\fi%
\ifnum\titrequatre=11{k}\fi%
\ifnum\titrequatre=12{l}\fi%
\ifnum\titrequatre=13{m}\fi%
}
\long\def\paraquatrel#1{\ifdim\lastskip<\smallskipamount
                \removelastskip\smallskip\fi
                 \vskip0pt plus.01\vsize\penalty-10
                  \vskip0pt
                  plus-.01\vsize{\bigskip
\Indentation{\global\advance\titrequatre by 1
\fontetitrequatre\quad \quad \quad \the\titreunl${\cdot}$\the\titredeuxl${\cdot}$\the\titretrois${\cdot}$\the\titrequatrel.\ }
{\hskip0mm\fontetitrequatre#1}\smallskip\nobreak\par}}

\ifx\optionkeys\oui
\def\drefun#1{\definexref{¤#1}{{\the\titreun}}{}} 
\def\drefdeux#1{\definexref{¤#1}{{\the\titreun}.{\the\titredeux}}{}}
\def\dreftrois#1{\definexref{¤#1}{{\the\titreun}.{\the\titredeux}.{\the\titretrois}}{}}
\else
\def\drefun#1{\definexref{prg#1}{{\the\titreun}}{}} 
\def\drefdeux#1{\definexref{prg#1}{{\the\titreun}.{\the\titredeux}}{}}
\def\dreftrois#1{\definexref{prg#1}{{\the\titreun}.{\the\titredeux}.{\the\titretrois}}{}}
\fi

%


  \long\def\propdeux#1#2#3#4{%
       \avance{\enonce}
       \leavevmode\edef\temp{#2}%
         \ifx\temp\empty 
          \else
           \definexref{#2}{#1~{\the\titreun.\the\enonce}}{enonces}
            \definexref{s#2}{{\the\titreun.\the\enonce}}{enonces}
             \fi
\medbreak
      \noindent{\bf#1\ {\bf\the\titreun.\the\enonce{#3}.}\enspace}{\sl#4\par}%
      \ifdim\lastskip<\medskipamount \removelastskip\penalty55\medskip\fi
   }

  \long\def\propun#1#2#3#4{%
      \avance{\enonce}
       \leavevmode\edef\temp{#2}%
        \ifx\temp\empty 
          \else
           \definexref{#2}{#1~{\the\enonce}}{enonces}
            \definexref{{s#2}}{{\the\enonce}}{enonces}
             \fi
   \medbreak
     \noindent{\bf#1\ {\bf\the\enonce{#3}.}\enspace}{\sl#4\par}%
     \ifdim\lastskip<\medskipamount \removelastskip\penalty55\medskip\fi
  }
  
  \long\def\prop#1#2#3#4{\ifnum\optionparag=1
                          \propdeux{#1}{#2}{\textfont1=\elevenib#3}{#4} \else\propun{#1}{#2}{\textfont1=\elevenib#3}{#4}\fi}

  \long\def\propt#1#2#3{\ifx\tpf\oui \prop{Th\'eo\-r\`eme}{#1}{#2}{#3}\par
                       \else\prop{Theorem}{#1}{#2}{#3}\par\fi}
  
  \long\def\propl#1#2#3{\ifx\tpf\oui\prop{Lem\-me}{#1}{#2}{#3}\par
                         \else\prop{Lemma}{#1}{#2}{#3}\par\fi}
  
  \long\def\propc#1#2#3{\ifx\tpf\oui\prop{Corol\-laire}{#1}{#2}{#3}\par
                         \else\prop{Corollary}{#1}{#2}{#3}\par\fi}

  \long\def\propd#1#2#3{\ifx\tpf\oui\prop{D\'efi\-nition}{#1}{#2}{#3}\par
                       \else\prop{Definition}{#1}{#2}{#3}\par\fi} 
  
  \long\def\proptd#1#2#3{\ifx\tpf\oui\prop{Th\'eor\`eme et d\'efi\-nition}{#1}{#2}{#3}\par
                       \else\prop{Theorem and definition}{#1}{#2}{#3}\par\fi}


  
  \newcount\optionparag\optionparag=1
  
  \long\def\section#1#2{\ifnum\optionparag=1 \paraun{#2} 
                        \else\goodbreak{\fontetitreun
  	                \Indentation{#1.\ }{#2}}\nobreak\smallskip\fi}

  \def\eqconstruct#1{\ifnum\optionparag=1{\the\titreun\hbox{$\cdot$}#1}\else{#1}\fi}
  \def\eqprint#1{{\rm (#1)}}

  
  
  \def\numref{oui}  
  
  \newcount\mesref\mesref=0 
  \def\defbib#1{\ifx\numref\oui\global\advance\mesref by 1 \definexref{#1}{{\the
                 \mesref}}{}\else\definexref{#1}{#1}{}\fi}

  
  \font\fourteenmsa=msam10 at 14pt
  \font\twelvemsa=msam10 at 12pt
  \font\tenmsa=msam10                 
  \font\ninemsa=msam10 at 9pt 
  \font\eightmsa=msam10 at 8pt 
  \font\sevenmsa=msam7 
  \font\sixmsa=msam10 at 6pt
  \font\fivemsa=msam5
  \newfam\msafam\textfont\msafam=\tenmsa\scriptfont\msafam=\sevenmsa\scriptscriptfont\msafam=\fivemsa
  
  \font\fourteenbb=msbm10 at 14pt
  \font\twelvebb=msbm10 at 12pt
  \font\tenbb=msbm10                   
  \font\ninebb=msbm10 at 9pt 
  \font\eightbb=msbm10 at 8pt 
  \font\sevenbb=msbm7 
  \font\sixbb=msbm10 at 6pt
  \font\fivebb=msbm5 
  \newfam\bbfam\textfont\bbfam=\tenbb\scriptfont\bbfam=\sevenbb\scriptscriptfont\bbfam=\fivebb
  \def\bb{\fam\bbfam\tenbb}%

  \font\twelvescaln=eusm10 at 12pt
  \font\tenscaln=eusm10                
  \font\ninescaln=eusm10 scaled 900
  \font\eightscaln=eusm10 scaled 800
  \font\sevenscaln=eusm10 scaled 700
  \font\sixscaln=eusm10 scaled 600
   
  \newfam\scalnfam\textfont\scalnfam=\tenscaln\scriptfont\scalnfam=\sevenscaln\scriptscriptfont\scalnfam=\sixscaln
  \def\scaln{\fam\scalnfam\tenscaln}%
  \def\scal{\scaln}
  
  \font\tenscalb=eusb10                

  \font\sevenscalb=eusb10 scaled 700

  \newfam\scalbfam\textfont\scalbfam=\tenscalb\scriptfont\scalbfam=\sevenscalb
  %
  
  %
  %
  \font\fourteenrm=cmr12 scaled 1200
  \font\elevenrm=cmr10 at 11pt
  \font\twelverm=cmr12
  \font\ninerm=cmr9
  \font\eightrm=cmr8      
  \font\sevenrm=cmr7
  \font\sixrm=cmr6

  \font\tenpcap=cmcsc10                        
  
  \font\eightpcap=cmcsc8
  \font\sevenpcap=cmcsc10 scaled 700
  
  \newfam\pcapfam\textfont\pcapfam=\tenpcap\scriptfont\pcapfam=\sevenpcap
  \def\pcap{\fam\pcapfam\tenpcap}
  
  \font\gras=cmbx10 scaled 1100             
  \font\TIT=cmbx12 scaled 1600              

  \font\fourteenbf=cmbx10 scaled 1400
  
  \font\twelvebf=cmbx12
  \font\elevenbf=cmbx10 at 11pt
  \font\ninebf=cmbx9  
  \font\eightbf=cmbx8
  \font\sixbf=cmbx6
  
  \font\tengot=eufm10                           
   
   at 8truept 
  \font\sevengot=eufm7 
   at 6 truept 
   
  \newfam\gotfam
  \textfont\gotfam=\tengot\scriptfont\gotfam=\sevengot
  %



  \font\fourteenti=cmbxti10 at 14pt
  
  \font\twelveti=cmbxti10 scaled 1200
  \font\eleventi=cmbxti10 at 11pt

  %
  %
  \font\twelveit=cmti12	
  \font\elevenit=cmti10 scaled 1100
  \font\nineit=cmti9
  \font\eightit=cmti8
  \font\sevenit=cmti7

  %
  %
  
  \font\fourteenib=cmmib10 scaled 1400
  \font\twelveib=cmmib10 scaled 1200
  \font\elevenib=cmmib10 scaled 1100
  \font\tenib=cmmib10
  \font\nineib=cmmib10 scaled 900

  %
  %
  
  \font\eleveni=cmmi10 scaled 1100
  \font\ninei=cmmi9
  \font\eighti=cmmi8 
  \font\seveni=cmmi7 			                
  \font\sixi=cmmi6
  
  \font\ninesl=cmsl9                    
  \font\eightsl=cmsl8 
  \font\sevensl=cmsl10 at 7pt

  \font\ninett=cmtt9                    
  \font\eighttt=cmtt8
  \font\seventt=cmtt10 scaled 700

  \font\fourteensy=cmsy10 scaled 1400
  \font\twelvesy=cmsy10 scaled 1176
  
  \font\ninesy=cmsy9                      
  \font\eightsy=cmsy8
  \font\sixsy=cmsy6
  
  \font\fourteenex=cmex10 at 14pt
  \font\twelveex=cmex10 at 12pt
  \font\nineex=cmex10 at 9pt
  \font\eightex=cmex10 at 8pt
  \font\sevenex=cmex10 at 7pt
  \font\sixex=cmex10 at 6pt
  \font\fiveex=cmex10 at 5pt
  
   
  \font\fourteengp=cmmi10 at 14pt
  
  \font\twelvegp=cmmib10 at 12pt
  \font\elevengp=cmmib10 at 11pt
  \font\tengp=cmmib10                          
  \font\ninegp=cmmib10 at 9pt 
  \font\eightgp=cmmib8 
   
  \font\sixgp=cmmib6

  \def\gponze{\textfont0=\elevenbf\scriptfont0=\eightbf\scriptscriptfont0=\sixbf
           \textfont1=\elevengp\scriptfont1=\eightgp\scriptscriptfont1=\sixgp}
  \def\gpdouze{\textfont0=\twelvebf\scriptfont0=\tenbf\scriptscriptfont0=\ninebf
           \textfont1=\twelvegp\scriptfont1=\tengp\scriptscriptfont1=\ninegp}        
  
 \def\gpquatorze{\textfont0=\fourteenbf\scriptfont0=\twelvebf\scriptscriptfont0=\elevenbf
           \textfont1=\fourteengp\scriptfont1=\twelvegp\scriptscriptfont1=\elevengp}

  
  \expandafter\chardef\csname pre amssym.def at\endcsname=\the\catcode`\@
  \catcode`\@=11
  \def\undefine#1{\let#1\undefined}
  \def\newsymbol#1#2#3#4#5{\let\next@\relax
   \ifnum#2=\@ne\let\next@\msafam@\else
   \ifnum#2=\tw@\let\next@\bbfam@\fi\fi
   \mathchardef#1="#3\next@#4#5}
  \def\mathhexbox@#1#2#3{\relax
   \ifmmode\mathpalette{}{\m@th\mathchar"#1#2#3}%
   \else\leavevmode\hbox{$\m@th\mathchar"#1#2#3$}\fi}
  \def\hexnumber@#1{\ifcase#1 0\or 1\or 2\or 3\or 4\or 5\or 6\or 7\or 8\or
   9\or A\or B\or C\or D\or E\or F\fi}
  
  \def\setboxz@h{\setbox\z@\hbox}
  \def\wdz@{\wd\z@}
  \def\boxz@{\box\z@}
  
  \edef\msafam@{\hexnumber@\msafam}
  \mathchardef\dabar@"0\msafam@39
  
  \edef\bbfam@{\hexnumber@\bbfam}
  \def\widehat#1{\setboxz@h{$\m@th#1$}%
   \ifdim\wdz@>\tw@ em\mathaccent"0\bbfam@5B{#1}%
   \else\mathaccent"0362{#1}\fi}
  \def\widetilde#1{\setboxz@h{$\m@th#1$}%
   \ifdim\wdz@>\tw@ em\mathaccent"0\bbfam@5D{#1}%
   \else\mathaccent"0365{#1}\fi}
  \newsymbol\leqq 1335          
  \newsymbol\leqslant 1336
  \newsymbol\lessgtr 1337       
  \newsymbol\backprime 1038     
  \newsymbol\risingdotseq 133A  
  \newsymbol\fallingdotseq 133B 
  \newsymbol\succcurlyeq 133C   
  \newsymbol\geqq 133D          
  \newsymbol\geqslant 133E
  \newsymbol\nmid 232D
  \newsymbol\nexists 2040
  \newsymbol\smallsetminus 2272
  \newsymbol\varnothing 203F
  
  \catcode`\@=\active

  
  \catcode`\@=11
  
  \newcount\typofr\typofr=1
  
  \catcode`\;=\active
  \def;{\ifnum\typofr=1\relax\ifhmode\ifdim\lastskip>\z@\unskip\fi
     \kern.2em\fi\string;\else\string;\fi}
  
  \catcode`\:=\active
  \def:{\ifnum\typofr=1\relax\ifhmode\ifdim\lastskip>\z@\unskip\fi
  \penalty\@M\ \fi\string:\else\string:\fi}
  
  \catcode`\!=\active
  \def!{\ifnum\typofr=1\relax\ifhmode\ifdim\lastskip>\z@\unskip\fi
     \kern.2em\fi\string!\else\string!\fi}
  
  \catcode`\?=\active
  \def?{\ifnum\typofr=1\relax\ifhmode\ifdim\lastskip>\z@\unskip\fi
     \kern.2em\fi\string?\else\string?\fi}

  \def\francais{\typofr=1\def\tpf{oui}}
  
  \def\oui{oui}
  \francais
  
  \catcode`\@=12


  
  \def\og{\leavevmode\raise.24ex\hbox{$\scriptscriptstyle\langle\!\langle\>$}}    
  \def\fg{\leavevmode\raise.24ex\hbox{$\scriptscriptstyle\>\rangle\!\rangle$}}    

  \def\d{\,{\rm d}}

  \def\z{{\bb Z}}
  \def\r{{\bb R}}
  \def\CC{{\bb C}}
  \def\N{{\bb N}}

  \def\n{{\scal N}}

  \def\frac#1#2{{#1\over #2}}

  \def\qedbox{$\rlap{$\sqcap$}\sqcup$}           
  \def\qed{\nobreak\hfill\penalty250 \hbox{}\nobreak\hfill\qedbox\par\medskip}

  \def\¤{\S\thinspace}

  \def\¥{$\bullet$ }
  
  
  \def\e{{\rm e}}

  \def\ep{\varepsilon}
  \def\epsilon{\varepsilon}

  \def\phi{\varphi}
  \def\theta{\vartheta}
  \def\rho{\varrho}

  \def\pf{\noi{\it Proof. }}
  \def\nid{\ifnum\typofr=1\par\noindent{\it D\'emonstration. }\else\pf\fi}
  \def\noi{\noindent}
  \def\rem{\ifnum\typofr=1\noi{\it Remarque.}\ \else\noi{\it Remark.}\ \fi}
  \def\rems{\ifnum\typofr=1\noi{\it Remarques.}\ \else\noi{\it Remarks.}\ \fi}

  \def\1{{\bf 1}}
  \def\|{\Vert}

  \def\le{\leqslant}
  \def\ge{\geqslant}


  \def\li{\mathop{\rm li}\nolimits}

  \def\log{\mathop{\rm log}\nolimits}




  \def\pmb#1{\setbox0=\hbox{#1}%
  \kern-.025em\copy0\kern-\wd0\kern.05em\copy0\kern-\wd0\kern-.025em\raise .0433em\box0 }

  
  \skewchar\eighti='177 \skewchar\sixi='177
  \skewchar\eightsy='60 \skewchar\sixsy='60
  
  \def\eightpoint{%
  \textfont0=\eightrm\scriptfont0=\sixrm\scriptscriptfont0=\fiverm
  \def\rm{\fam0\eightrm}%
  \textfont1=\eighti\scriptfont1=\sixi
  \scriptscriptfont1=\fivei\def\oldstyle{\fam1\seveni}%
  \textfont2=\eightsy\scriptfont2=\sixsy\scriptscriptfont2=\fivesy
  \textfont3=\eightex\scriptfont3=\sixex
  \textfont\itfam=\eightit
  \def\it{\fam\itfam\eightit}%
  \textfont\slfam=\eightsl
  \def\sl{\fam\slfam\eightsl}%
  \textfont\bbfam=\eightbb \scriptfont\bbfam=\sixbb\scriptscriptfont\bbfam=\fivebb
  \def\bb{\fam\bbfam\eightbb}%
  \textfont\msafam=\eightmsa\scriptfont\msafam=\sixmsa
  \textfont\scalnfam=\eightscaln
  \def\scaln{\fam\scalnfam\eightscaln}
  \textfont\ttfam=\eighttt
  \def\tt{\fam\ttfam\eighttt}%
  \textfont\bffam=\eightbf\scriptfont\bffam=\sixbf\scriptscriptfont\bffam=\fivebf
  \def\bf{\fam\bffam\eightbf}%
  \textfont\pcapfam=\eightpcap
  \def\pcap{\fam\pcapfam\eightpcap}
  \abovedisplayskip=2pt plus2pt minus 2pt
  \belowdisplayskip=2pt plus1pt minus 2pt
  \abovedisplayshortskip= 1pt plus 2pt minus 1pt
  \belowdisplayshortskip= 1pt plus 2pt minus 1pt
  \smallskipamount=2pt plus 1pt minus 2pt
  \medskipamount=3pt plus 2pt minus 2pt
  \bigskipamount=7pt plus 3pt minus 3pt
  \setbox\strutbox=\hbox{\vrule height 5pt depth 2pt width 0pt}%
  \normalbaselineskip=9pt\normalbaselines\rm}

  \def\({\left(}
  \def\){\right)}
  
  \def\footnoterule{\kern -2pt\hrule width 7truecm\kern 2.4pt}
  
  \def\xnotedef#1{\definexref{#1}{\noexpand\number\footnotenumber}{Note}}%

  
  
  \def\ninepoint{%
  \textfont0=\ninerm\scriptfont0=\sixrm\scriptscriptfont0=\fiverm
  \def\rm{\fam0\ninerm}%
  \textfont1=\ninei\scriptfont1=\sixi
  \scriptscriptfont1=\fivei\def\oldstyle{\fam1\ninei}%
  \textfont2=\ninesy\scriptfont2=\sixsy\scriptscriptfont2=\fivesy
  \textfont3=\nineex\scriptfont3=\sixex
  \textfont\itfam=\nineit
  \def\it{\fam\itfam\nineit}%
  \textfont\slfam=\ninesl
  \def\sl{\fam\slfam\ninesl}%
  \textfont\bbfam=\ninebb\scriptfont\bbfam=\sixbb\scriptscriptfont\bbfam=\fivebb
  \def\bb{\fam\bbfam\ninebb}%
  \textfont\msafam=\ninemsa\scriptfont\msafam=\sixmsa\scriptscriptfont\msafam=\fivemsa
  \textfont\scalnfam=\ninescaln
  \def\scaln{\fam\scalnfam\ninescaln}
  \textfont\ttfam=\ninett
  \def\tt{\fam\ttfam\ninett}%
  \textfont\bffam=\ninebf\scriptfont\bffam=\sixbf\scriptscriptfont\bffam=\fivebf
  \def\bf{\fam\bffam\ninebf}%
  \abovedisplayskip=3pt plus2pt minus 2pt
  \belowdisplayskip=3pt plus1pt minus 2pt
  \abovedisplayshortskip= 2pt plus 2pt minus 1pt
  \belowdisplayshortskip= 2pt plus 2pt minus 1pt
  \smallskipamount=2pt plus 1pt minus 2pt
  \medskipamount=3pt plus 2pt minus 2pt
  \bigskipamount=7pt plus 3pt minus 3pt
  \setbox\strutbox=\hbox{\vrule height 5pt depth 2pt width 0pt}%
  \normalbaselineskip=10.5pt plus.3pt minus.3pt\normalbaselines\rm}

  \def\sevenpoint{%
  \textfont0=\sevenrm\scriptfont0=\sixrm\scriptscriptfont0=\fiverm
  \def\rm{\fam0\sevenrm}%
  \textfont1=\seveni\scriptfont1=\sixi
  \scriptscriptfont1=\fivei\def\oldstyle{\fam1\seveni}%
  \textfont2=\sevensy\scriptfont2=\sixsy\scriptscriptfont2=\fivesy
  \textfont3=\sevenex\scriptfont3=\fiveex
  \textfont\itfam=\sevenit
  \def\it{\fam\itfam\sevenit}%
  \textfont\slfam=\sevensl
  \def\sl{\fam\slfam\sevensl}%
  \textfont\bbfam=\sevenbb \scriptfont\bbfam=\sixbb\scriptscriptfont\bbfam=\fivebb
  \def\bb{\fam\bbfam\sevenbb}%
  \textfont\msafam=\sevenmsa\scriptfont\msafam=\sixmsa
  \textfont\scalnfam=\sevenscaln
  \def\scaln{\fam\scalnfam\sevenscaln}
  \textfont\bffam=\sevenbf\scriptfont\bffam=\sixbf\scriptscriptfont\bffam=\fivebf
  \def\bf{\fam\bffam\sevenbf}%
  \textfont\ttfam=\seventt
  \abovedisplayskip=2pt plus2pt minus 2pt
  \belowdisplayskip=2pt plus1pt minus 2pt
  \abovedisplayshortskip= 1pt plus 2pt minus 1pt
  \belowdisplayshortskip= 1pt plus 2pt minus 1pt
  \smallskipamount=2pt plus 1pt minus 2pt
  \medskipamount=3pt plus 2pt minus 2pt
  \bigskipamount=7pt plus 3pt minus 3pt
  \setbox\strutbox=\hbox{\vrule height 5pt depth 2pt width 0pt}%
  \normalbaselineskip=9pt\normalbaselines\rm}

 \def\onzepts{%
 \textfont0=\elevenrm\scriptfont0=\ninerm
 \textfont1=\elevenib\scriptfont1=\ninei}

\def\douzepts{%
  \textfont0=\twelverm\scriptfont0=\tenrm\def\rm{\fam0\twelverm}%
  \textfont1=\twelveib\scriptfont1=\teni%
  \textfont2=\twelvesy\scriptfont2=\tensy\scriptscriptfont2=\eightsy
  \textfont3=\twelveex
  \textfont\itfam=\twelveti
  \def\it{\fam\itfam\twelveti}%
  \textfont\bffam=\twelvebf\scriptfont\bffam=\tenbf\scriptscriptfont\bffam=\eightbf
  \def\bf{\fam\bffam\twelvebf}%
  \textfont\bbfam=\twelvebb \scriptfont\bbfam=\tenbb
  \def\bb{\fam\bbfam\twelvebb}%
  \textfont\msafam=\twelvemsa\scriptfont\msafam=\tenmsa
  \textfont\scalnfam=\twelvescaln
  \normalbaselineskip=15pt\normalbaselines\rm}

\def\quatorzepts{%
  \textfont0=\fourteenrm\scriptfont0=\twelverm\def\rm{\fam0\fourteenrm}%
  \textfont1=\fourteenib\scriptfont1=\twelveib%
  \textfont2=\fourteensy\scriptfont2=\twelvesy\scriptscriptfont2=\tensy
  \textfont3=\fourteenex
  \textfont\itfam=\fourteenti
  \def\it{\fam\itfam\fourteenti}%
  \textfont\bffam=\fourteenbf\scriptfont\bffam=\twelvebf\scriptscriptfont\bffam=\tenbf
  \def\bf{\fam\bffam\fourteenbf}%
  \textfont\bbfam=\fourteenbb \scriptfont\bbfam=\twelvebb
  \def\bb{\fam\bbfam\fourteenbb}%
  \textfont\msafam=\fourteenmsa\scriptfont\msafam=\twelvemsa
  \textfont\scalnfam=\twelvescaln
  \normalbaselineskip=18pt\normalbaselines\rm}


\def\picture #1 by #2 (#3){\leavevmode\vbox to #2{
     \hrule width #1 height 0pt depth 0pt
      \vfill
       \special{picture #3}}}

\def\illustration #1 by #2 (#3) scaled #4{\dimen1=#2
  \divide\dimen1 by 1000
  \multiply\dimen1 by #4
  \vtop to \dimen1{\dimen1=#1
  \divide\dimen1 by 1000
  \multiply\dimen1 by #4
  \hsize=\dimen1\vss
  \noindent\special{illustration #3 scaled #4}}}

\def\frac#1#2{{#1\over#2}}
\centerline {\TIT{A fourth derivative test}}
\vskip 5 mm
\centerline {\TIT {for exponential sums}}
\vskip 1 cm
\centerline {\fourteenrm {  by O. Robert and P. Sargos}}
\bigskip
\bigskip
\bigskip
\centerline {\pcap {Abstract}}
\medskip
\ninerm
We give an upper bound for the exponential sum $\sum_{m=1}^M\e^{2i\pi f(m)}$
in terms of $M$ and $\lambda$, where $\lambda$ is a small positive
number which denotes the size of the fourth derivative of the real
valued function $f$. The classical van der Corput's exponent 1/14 is
improved into 1/13 by reducing the problem to a mean square value
theorem for triple exponential sums.
\bigskip
\bigskip

\tenrm
\bigskip
\anote{}{1991 {\it Mathematics Subject Classification} : 11L07}
\anote{}{{\bf Key words} : Exponential sums, Diophantine systems.}
\paraun{Introduction and statement of the result }
\medskip
The aim of this paper is to find upper bounds for the exponential sum
$$S_M=\sum_{m=1}^M\e(f(m)),\eqdef{sum}$$
where we have set $\e(x)$ for $\e^{2i\pi x}$ and where $M$ is a large
integer and $f : [1,M] \to \r $ is a four times continuously
differentiable function which satisfies  van der Corput's condition
:
$$\lambda\le f^{(4)}(x) \ll \lambda , \ {\rm for } \ 1\le x\le M \eqdef{VdCcond}$$
where $\lambda $ is a small positive number, and where the Vinogradov's  symbol $u\ll v$ means that there exists an absolute positive
constant $C$ such that $|u|\le C v$. Under the condition \eqref{VdCcond}, van der Corput has
 obtained the classical following bound (cf [3], Theorem 2.8)
$$S_M\ll M\lambda^{1/14}, \ {\rm provided \ that \ }
M\gg\lambda^{-4/7}. \eqdef{quatorze}$$
\medskip
The proof consists in applying twice Weyl and van der Corput's
A-process (cf [3], Lemma 2.5), and then van der Corput's inequality (cf
[3], Theorem 2.2).
Slight improvements on \eqref{quatorze} have been obtained later, but
only under stronger hypothesis (see e.g. [3] or [4]). It is interesting
to notice that the bound \eqref{quatorze} can be improved without any
new hypothesis, as a consequence of a strong result of Bombieri and
Iwaniec [2] on the mean value of eighth powers of simple cubic
exponential sums. The deduction has been made in [7] with the bound : 
$$S_M\ll_{\ep} M^{1+\ep}\lambda^{3/40}, {\rm provided \ that \ }
M\gg\lambda^{-3/5}. \eqdef{troisquarante}$$
Here and in the sequel, the symbol $\ll_{\ep}$ means that the inequality holds for each
$\ep >0$ and that the implied constant depends at most on $\ep $ and
on the previous implied constants.
Our result can be stated as follows :

\proclaim Theorem 1. If the condition \eqref{VdCcond} is satisfied,
then the two equivalent properties :
$$ S_M\ll_{\ep}M^{1+\ep}\lambda^{1/13}, 
{\rm \ provided \ that \ }  M\gg \lambda^{-8/13}.\eqdef{th1}$$
and
$$ S_M\ll_{\ep}M^{\ep}(M\lambda^{1/13}+\lambda^{-7/13})\eqdef{th1'}$$
hold true.
\par
\bigskip
We conclude this section with some remarks and comments while sections
2,3,4 are entirely devoted to the proof of Theorem 1.
\medskip
{\bf Exponent pairs} Most of problems in analytic number theory where
exponential sums occur, involve phase functions which satisfy much
more than \eqref{VdCcond}. Namely, these functions $f:[M,2M]\to \r$
satisfy conditions (3.3.3) of [3](we shall call them \og semi-monomial
functions \fg ).
Bounds for exponential sums $S_M=\sum_{m=M+1}^{2M}\e(f(m))$ are then
obtained in terms of exponent pairs (see \S 3.3 of [3]). For
semi-monomial functions, the bound :
$$S_M\ll_{\ep}M^{1+\ep}\lambda^{\theta}\eqdef{exptheta}$$
corresponds to the property :
$$ (\theta+\ep,1-3\theta+\ep) {\rm \ is \ an \ exponent \ pair \ for
\ each \ } \ep>0. \eqdef{exppair}$$ 
Thus, our Theorem 1 implies that \eqref{exppair} holds for
$\theta=1/13$ (to see this, we only have to complete the proof in the
case $M\ll \lambda^{-8/13}$ by means of the classical exponent pair
$(2/18,13/18)=ABA^2B(0,1)$).
But this value of $\theta$ is not the best known.
Indeed, the refinement by Huxley and Kolesnik [5] (see also [4], \S
19.3) of Huxley's deep method for exponential sums with a large second
derivative ([4], \S 17.4), yields a better value of $\theta $. Namely,
one can take out from table 19.2 of [4] the following result :
\proclaim Theorem {\rm(Huxley and Kolesnik)}. The property \eqref{exppair}
holds for $\theta=\frac{516247}{6629696}=\frac{1}{12.84...}.$
\par
\medskip
The interest of our Theorem 1 consists, on the one hand, in the
simplicity of its proof and, on the other hand, in the wider range of
its applications, particularly to short exponential sums.  
\medskip
{\bf Van der Corput's exponent} The exponent $1/14$ in
\eqref{quatorze} can be sharpened into $1/13$, at least with some
restrictions on the relative size of $M$ and $\lambda $. The question
of knowing how much van der Corput's exponent $1/14$ can be bettered
and under which conditions, arises naturally.
\smallskip
We have heuristic proofs of the two following assertions that we state
as conjectures :
\medskip
\proclaim Conjecture 1. Under the hypothesis \eqref{VdCcond}, we have
$$S_M\ll_{\ep}M^{1+\ep}\lambda^{3/38},{\rm \ provided \ that \
}M\gg\lambda^{-13/19}.\eqdef{conjecture1}$$
\par
\medskip
\proclaim Conjecture 2. Under the hypothesis \eqref{VdCcond}, we have
$$S_M\ll_{\ep}M^{1+\ep}\lambda^{1/12},{\rm \ provided \ that \
}M\gg\lambda^{-1}.\eqdef{conjecture2}$$
\par
\medskip
This last conjecture, if true, is far from implying that the pair
$(1/12+\ep,9/12+\ep)$ is an exponent pair for each $\ep>0$. The
restriction $M\gg\lambda^{-1}$ in \eqref{conjecture2} is quite
constraining and we think that, perhaps, it cannot be
weakened. Furthermore, if we restrict conjecture 2 to semi-monomial
phase functions, then Huxley's results already imply
\eqref{conjecture2} (cf [4], \S 17.4).
\medskip
{\bf Very short exponential sums} In the opposite direction, we have
the following improvement of \eqref{quatorze} (cf [7], Lemma 2.6) :
$$S_M\ll M\lambda^{1/14},{\rm \ provided \ that \
}M\gg\lambda^{-3/7},\eqdef{3quatorze}$$
which concerns shorter exponential sums. It would be of interest,
both in itself and for the applications, to find the infimum of positive
real $\beta $ such that the bound
$$S_M\ll M\lambda^{1/14},{\rm \ provided \ that \
}M\gg\lambda^{-\beta},\eqdef{beta}$$
holds under the hypothesis \eqref{VdCcond}. The example
$$f(m)=\frac{m^4}{100M^4},$$ in which we have $|S_M|\gg
\lambda^{-1/4},$ shows that $\beta\ge 9/28,$ so we have
$$ 9/28\le \beta\le 3/7\eqdef{betaenca}$$
\medskip
{\bf Outline of proof } At first, we apply van der Corput's A-process to the
initial sum \eqref{sum} and get a double sum in the variables $h$ and
$m$. Then we apply $A\times A-$process to the new double sum and get a
quadruple sum in the variables $r,q,h,m$. At last, we shift the main
variable $m$ to produce a new variable $n$. This can be sketched in
the following diagram :
$$\left\{\matrix{\displaystyle
{\ \ S_M=\sum_m\e(f(m))\mathop{\longrightarrow}^{A}\sum_h\sum_m\e(\Delta_hf(m))}\hfill\cr
\displaystyle{\ \ \mathop{\longrightarrow}^{A\times A}
\sum_r\sum_q\sum_h\sum_m\e(\Delta_{h+r}f(m)-\Delta_hf(m+q))}\hfill\cr
\displaystyle{\ \  \mathop{\longrightarrow}^{{\rm shift}}\sum_r\sum_m\left|\sum_h\sum_n\sum_q\e(\Delta_{h+r}f(m+n)-\Delta_hf(m+n+q))\right|}\hfill\cr}
\right. \eqdef{outline}$$
where we have set $\Delta_hf(m)$ for $f(m+h)-f(m-h)$.
By expanding the phase in the last exponential sum by means of
Taylor's formula, we are in a position to apply Bombieri and Iwaniec's
double large sieve [1]. Thus we have reduced the initial problem into
that of counting the number of solutions of a (very particular)
diophantine system, which is the purpose of our Theorem 2. 
The whole proof is self contained and elementary. 

\medskip
\paraun{Preliminary lemmas}
\medskip
We recall some basic lemmas.
\medskip
{\twelvebf 2.1. Weyl and van der Corput A x A lemma}
\medskip
\proclaim Lemma 1.
Let $M$ and $H$ be positive integers and let $(a(m,h))_{(m,h)\in \z^2}$ be complex numbers which are zero whenever
$(m,h)$ is outside the compact $[1,M]\times [1,H]$. We set
$$S=\sum_{(m,h)\in \z^2}a(m,h)$$
and we choose two integers $Q$ and $R$ such that $1\le Q\le M$ and
$1\le R\le H$. We then have :

$$S^2\ll\frac{MH}{QR}\sum_{|q|<Q}\sum_{|r|<R}\left(1-\frac{|q|}{Q}\right)\left(1-\frac{|r|}{R}\right)\sum_{(m,h)\in
\z^2}a(m+q,h)\overline{a(m,h+r)}\eqdef{lemAA}$$
\par
\medskip
For the proof, see [3], Lemma 6.1.\qed
\medskip
{\twelvebf 2.2. Partial summation for multiple sums}
\medskip
We give a general statement of partial summation for $k$-dimensional
sums, where $k$ is a positive integer. We need some notations.
\smallskip
Let $M_1,...,M_k$ be positive integers and set :
$${\cal P}=[1,M_1]\times ... \times [1,M_k]\subset \r^k. \eqdef{pave}$$
Let $I$ be any finite set and, for each fixed $i\in I$, let
$\phi_i:{\cal P}\to \CC $ be a function which satisfies the following
regularity condition.
\smallskip 
For each integer $r\ (0\le r\le k)$, for each $(j_1,...,j_r)$ such
that $1\le j_s\le k \ (1\le s \le r)$ and $j_s\not =j_t$ for $s\not
=t$, the $r$th order derivative $\frac{\partial^r\phi_i}{\partial
x_{j_1}...\partial x_{j_r}}$ exists and is continuous on ${\cal P}$
and satisfies the bound :
$$\left|\frac{\partial^r\phi_i}{\partial
x_{j_1}...\partial x_{j_r}}({\bf x})\right|\le
\frac{D}{M_{j_1}...M_{j_r}} \ {\rm whenever }\ i\in I, \ {\bf x}\in
{\cal P}, \ r\in\{0,...,k\},\ 1\le j_1<...<j_r\le k\eqdef{abelcond}$$
for some $D>0$. We recall that the bound \eqref{abelcond} in case
$r=0$ means that $|\phi_i({\bf x})|\le D$ for each $i\in I$ and ${\bf
x}\in{\cal P}$.
\smallskip
Let us now consider the $k$-dimensional sum :
$$S_0=\sum_{i\in I}\left|\sum_{{\bf m}\in {\cal P}\cap \N^k}a_i({\bf
m})\phi_i({\bf m})\right|\eqdef{abelsum}$$
where $(a_i({\bf m}))_{i\in I,{\bf m}\in {\cal P}\cap \N^k}$ is any
given family of complex numbers. We can now state our lemma for
$k$-dimensional partial summation.
\medskip
\proclaim Lemma 2. Let the above notations and hypothesis hold. We
then have :
$$S_0\le 2^k D \max_{{\cal P}'}\sum_{i\in I}\left|\sum_{{\bf m}\in {\cal P}'\cap \N^k}a_i({\bf
m})\right|\eqdef{abel}$$
where the maximum has to be taken over all possible sets of the form 
${\cal P}'=[1,M'_1]\times ...\times [1,M'_k]\subset {\cal P}$.
\par
\medskip
{\it Proof.} The proof goes by recurrence on $k$. When $k=1$, the result is nothing but the classical one-dimensional partial summation. We suppose that the result is true up to $k-$dimensional sums, and we want to prove that it is true for $(k+1)-$dimensional sums.
\medskip
The $(k+1)-$dimensional sum $S_0$ may be written as 

$$S_0=\sum_{i\in I}\left|\sum_{{\bf m}\in {\cal P}_k}\sum_{n=1}^Na_i({\bf m}, n)\phi_i({\bf m},n)\right|$$
\medskip
In order to apply one dimensional partial summation to the sum in $n$, we set $A_i({\bf m},n)=\sum_{\nu=1}^na_i({\bf m},\nu)$ and $\psi_{i,n}({\bf m})=\phi_i({\bf m},n)-\phi_i({\bf m},n+1)$. We have

$$S_0\le \sum_i\left|\sum_{{\bf m}\in {\cal P}_k}A_i({\bf
m},N)\phi_i({\bf m},N)\right|+\sum_i\sum_{n=1}^{N-1}\left|\sum_{{\bf
m}\in {\cal P}_k}A_i({\bf m}, n)\psi_{i,n}({\bf m})\right|$$
\medskip
We apply the recurrence hypothesis to both terms :
$$S_0\le 2^k D\max_{{\cal P}'_k}\sum_i\left|\sum_{{\bf m}\in {\cal
P}'_k}\sum_{\nu=1}^Na_i({\bf m}, \nu)\right|+ 2^k\frac{D}{N}\max_{{\cal
P}'_k}\sum_i\sum_{n=1}^{N-1}\left|\sum_{{\bf m}\in {\cal
P}'_k}\sum_{\nu=1}^na_i({\bf m},\nu)\right|,$$
and the desired result follows.\qed

\medskip
{\twelvebf 2.3. Third derivative test and partial summation}
\medskip
The following lemma is not essential in the proof of Theorem 1, but it
gives rise to simplifications.
\proclaim Lemma 3. Let $M$ be a positive integer, and let $g$ and
$u:[1,M]\to \r $ be two functions, respectively ${\cal C}^3$ and
${\cal C}^1$, such that :
$$\mu\le |g'''(x)|\ll \mu \ {\rm and } \ u'(x)\ll\mu^{1/2},\ {\rm
for}\ 1\le x\le M,\eqdef{thirdcond}$$
where $\mu $ is a small positive number. We then have :
$$\sum_{m=1}^M\e(g(m)+u(m))\ll M\mu^{1/6}+\mu^{-1/3}.\eqdef{third}$$   
\par
{\it Proof} : If $u\equiv 0$, Lemma 3 is the third derivative test for
exponential sums ([7], Corollary 4.2). If $M\ll \mu^{-1/2}$, we can
eliminate without cost the term $u(m)$ by a (one dimensional) partial
summation, and \eqref{third} follows. Now, we suppose $M\gg
\mu^{-1/2}$. Then we divide the initial sum into $O(M\mu^{1/2})$ sums
of length $\mu^{-1/2}$ and we apply the previous case to each short sum. \qed 
\medskip
{\twelvebf 2.4. Double large sieve inequality}
\medskip
We consider the exponential sum :
$${ \widetilde S}=\sum_{0<|r|<R}\sum_{m=1}^M\left|\sum_{0<|q|<Q}\sum_{h=H}^{2H-1}\sum_{n=1}^Nb_r(q,h,n)\e(x_mP_1(r,q,h,n)+y_mP_2(r,q,h,n))\right|,\eqdef{sievesum}$$
with the following notations :  $R,M,Q,H,N$ are positive integers ;
$b_r(q,h,n)$ are complex numbers with modulus at most one ;  $P_1$ and $P_2$ are polynomials in four variables, $P_1$
with integer coefficients and $P_2$ with real coefficients ;
$(x_m)_{1\le m\le M}$ and $(y_m)_{1\le m\le M}$ are two families of
real numbers. 
\medskip
We suppose that the next two inequalities hold :
\medskip

$$\max_{1\le m,m'\le M}|y_{m}-y_{m'}|\le \mu \eqdef{mu}$$
$$\max_{r,q,h,n}|P_i(r,q,h,n)|\le X_i, \ (i=1,2)\eqdef{Xi}$$
the latter maximum being taken over all quadruples $(r,q,h,n)$ of
integers such that
$$0<|r|<R,\ 0<|q|<Q,\ H\le h<2H,\ 1\le n\le N.$$
We introduce the numbers ${\cal N}$ and ${\cal B}$ which correspond to
spacing problems
$$\n=\max_{1\le Q_1<Q}\#\left\{
\matrix{
(r,q_1,q_2,h_1,h_2,n_1,n_2)\in\z^7  \ \ {\rm which \ satisfy}\cr
 {\rm conditions}\ (2.11.a)...(2.11.d)\cr}
 \right\}\eqdef{N}$$
$$\left\{
\matrix{
0<|r|<R, |Q_1|\le|q_i|<\min(2Q_1,Q),H\le h_i<2H,1\le n_i\le N,{\rm
for\ }i=1,2 &(2.11.a)\cr
q_1q_2>0& (2.11.b)\cr
P_1(r,h_1,q_1,n_1)=P_1(r,h_2,q_2,n_2)& (2.11.c)\cr
| P_2(r,h_1,q_1,n_1)-P_2(r,h_2,q_2,n_2)|\le 1/\mu & (2.11.d)}
\right.
$$
and
$${\cal B}=\#\left\{
(m_1,m_2)\in\{1,...,M\}^2
\left|\ \matrix{\|x_{m_1}-x_{m_2}\|\le X_1^{-1}\hfill\cr
{\rm and} \hfill\cr
|y_{m_1}-y_{m_2}|\le X_2^{-1}\hfill\cr}\right.
\right\}\eqdef{B}$$
with the usual notation : $ \|x\|=\min_{m\in\z}|x-m|.$
\medskip
We can now state the double large sieve inequality in the particular
form which will be needed later.
\proclaim Lemma 4. With the above notations, we have
:
$${ \widetilde S}^2\ll R(1+X_1)(1+\mu X_2)\n{\cal B}(\log Q)^2.\eqdef{sieve}$$
\par
\medskip
\medskip
{\it Proof} : We set 
$${ \widetilde S}(r,Q_1)=\sum_{m=1}^M\left|\sum_{Q_1<|q|<\max(2Q_1,Q)}\sum_{h=H}^{2H-1}\sum_{n=1}^Nb_r(q,h,n)\e(x_mP_1(r,q,h,n)+y_mP_2(r,q,h,n))\right|,$$ 
so that we have
$${\widetilde S}\ll \max_{1\le Q_1<Q}\left(\sum_{0<|r|<R}{S}(r,Q_1)\right)\log
Q. \eqdef{Sr}$$ 
 We apply Bombieri and Iwaniec's double large sieve [1] (cf also [3], Lemma 7.5
or [4], Lemma 5.6.6) to each sum ${ \widetilde S}(r,Q_1)$ and we get :

$${\widetilde S}(r,Q_1)^2\ll(1+X_1)(1+\mu X_2)\n(r,Q_1){\cal B},\eqdef{dls}$$

where $\n(r,Q_1)$ is the number of $(q_1,q_2,h_1,h_2,n_1,n_2)\in\z^6$
which satisfy  conditions (2.11.a),..., (2.11.d).

Next we apply Cauchy's
inequality and we report \eqref{dls} into \eqref{Sr} to get \eqref{sieve} \qed
\medskip

{\twelvebf 2.5. Iwaniec and Mozzochi's arithmetic lemma}
\medskip
For a positive integer $n$, we set
$$\tau(n)=\sum_{d|n}1 \ {\rm and} \
\sigma(n)=\sum_{d|n}d.\eqdef{tausigma}$$
The following lemma is contained in the proof of theorem 14.1 of
[6]. A complete and independent proof may be found in Lemma 13.1.2 of
[4](see also [8]).
\medskip
\proclaim Lemma 5. Let $a,b,c$ be three non zero integers, with
$gcd(a,b,c)=1$ and $c>0$. Let $V\ge 1,\ \alpha<\beta $ be real
numbers. We denote by ${\cal V}$ the number of triplets $(u,v,w)$ with
non zero integers such that $gcd(u,v,w)=1$, $V\le v\le 2V$ and :
$$au+bv+cw=0 \ {\rm and} \ \alpha\le \frac{u}{v}\le \beta.\eqdef{hyperplan}$$ 
We then have :
$${\cal V}\ll
\tau(c)+(\beta-\alpha)V^2\frac{\sigma(c)}{c^2}\eqdef{IwM}$$
\qed
\par

\paraun{The Diophantine problem}

\medskip
{\twelvebf 3.1. Statement of the result}
\medskip
Let $R,H,Q,N$ be real numbers $\ge 1$ with $R\le H/2$ and let
$\delta $ be a positive number. We denote by ${\cal
N}(R,Q,H,N,\delta)$ the number of integer points
$(r,q_1,q_2,h_1,h_2,n_1,n_2)\in\z^7$ lying in the domain :

$$0<|r|<R,\ Q\le| q_i|<2Q,\ H\le h_i<2H,\ 1\le n_i\le N\ {\rm for \ }
i=1,2,\ {\rm and\ } q_1q_2>0,\eqdef{paramcond}$$
and satisfying the system :
$$\left.\matrix{
rn_1+h_1q_1=rn_2+h_2q_2 \hfill \cr
     |rn_1^2+2h_1q_1n_1+h_1q_1^2 -(rn_2^2+2h_2q_2n_2+h_2q_2^2 )| \le \delta HQ^2\hfill \cr
}\right\} \eqdef{syst1}$$
\medskip
\proclaim Theorem 2. The number of solutions of the diophantine system
\eqref{syst1}, lying in the domain \eqref{paramcond} satisfies the bound
$${\cal N}(R,H,Q,N,\delta)\ll_{\ep}(RNHQ)^{1+\ep}(1+\delta Q).\eqdef{moyenne}$$
\par
The rest of this section is devoted to the proof of Theorem 2.
\medskip

{\twelvebf 3.2. Reduction of the problem}
\medskip
We reduce the diophantine system \eqref{syst1} into a simpler one by
means of easy calculations. Let ${\cal J}_1(R,Q,H,\delta)$ be the
number of integer points $(r,q_1,q_2,h_1,h_2,d)\in\z^6$ lying in the domain

$$0<|r|<R,\ Q\le |q_i|<2Q,\ H\le h_i<2H,\ {\rm for \ }
i=1,2,\ q_1q_2>0\ {\rm and}\ 0<|d|\ll(1+\delta)Q \eqdef{paramnewcond}$$
satisfying the system
$$\left.\matrix{
rd+h_1q_1-h_2q_2=0 \hfill \cr
     |rd^2+2h_1q_1d+h_1q_1^2 -h_2q_2^2 | \le \delta HQ^2,\hfill \cr
}\right\} \eqdef{syst2}$$
just as the additional condition :
$${\rm gcd}(d,q_1,q_2)=1,\ {\rm gcd}(r,h_1,h_2)=1.\eqdef{parareduction}$$
\medskip
\proclaim Lemma 6. With the above notations, we have :
$${\cal N}(R,Q,H,N,\delta)\ll_{\ep}(RNHQ)^{1+\ep}(1+\delta
Q)+N\sum_{1\le j\le R}\sum_{1\le k\le Q}{\cal J}(R/j,Q/k,H/j,\delta)\eqdef{reduction}$$
\par

{\it Proof} : we set $n_1=n_2+d$ and we insert this in
\eqref{syst1}. We observe that the terms containing $n_2$ cancel out
each other and we obtain \eqref{syst2}.
On the other hand, the system

$$\left.\matrix{
rd+h_1q_1-h_2q_2=0 \hfill \cr
     |(h_1q_1+h_2q_2)d+h_1q_1^2 -h_2q_2^2 | \le \delta HQ^2,\hfill \cr
}\right\} \eqdef{syst3}$$
is equivalent to \eqref{syst2}. From it, since $q_1$ and $q_2$ have
the same sign, we  deduce that
$$d\ll (1+\delta)Q\eqdef{dcond}$$
Using only the first line of \eqref{syst1}, we see that the number of
solutions of \eqref{syst1} with $d=0$ (i.e. $n_1=n_2$) is
$O_{\ep}((RQHN)^{1+\ep}),$ so that we may suppose now $d\not = 0$.
Let $j$ and $k$ be two positive integers and let
${\cal J}(j,k)$ be the number of integer points
$(r,q_1,q_2,h_1,h_2,d)\in \z^6,$ lying in the domain
\eqref{paramnewcond}, satisfying  system \eqref{syst2} and the
additional condition :

$${\rm gcd}(d,q_1,q_2)=k,\ {\rm gcd}(r,h_1,h_2)=j.$$
The following bound is then obvious :

$${\cal N}(R,Q,H,N,\delta)\ll_{\ep}(RNHQ)^{1+\ep}+N\sum_{1\le j\le
R}\sum_{1\le k<2Q}{\cal J}(j,k)\eqdef{obvious}$$
and it may be transformed into
$${\cal
N}(R,Q,H,N,\delta)\ll_{\ep}(RNHQ)^{1+\ep}(1+\delta)+N\sum_{1\le j\le
R}\sum_{1\le k\le Q}{\cal
J}(j,k).\eqdef{improved}$$
Indeed, if we assume that $k={\rm gcd}(d,q_1,q_2)$ is greater than
$Q$, then we have $q_1=q_2=k$ and there are $O(1+\delta)$
possibilities for $d$, so that the total number of solutions of
\eqref{syst2} with ${\rm gcd}(d,q_1,q_2)>Q$, is $O(RQHN(1+\delta))$
(we have only to use the first line of \eqref{syst2}). We have thus proved
\eqref{improved}.
But, for $j$ and $k$ fixed, with $1\le j \le R$ and $1\le k\le Q$, we
may divide the first line of \eqref{syst2} by $jk$ and the second line
of \eqref{syst2} by $jk^2$. The real numbers $R/j,Q/k$ and $H/j$ are
$\ge 1$, with $R/j\le \frac{1}{2}H/j$ and we have
$${\cal J}(j,k)={\cal J}(R/j,Q/k,H/j,\delta).$$ 
Thus \eqref{improved} implies \eqref{reduction} and the proof of Lemma
6 is complete.\qed
\medskip
From Lemma 6, we deduce that Theorem 2 is a consequence of the
following lemma :
\proclaim Lemma 7. Let $R,Q,H,\delta$ be positive real numbers with
$R\ge 1,\ Q\ge 1$ and $H\ge 2R$. We have :
$${\cal
J}(R,Q,H,\delta)\ll_{\ep}(RQH)^{1+\ep}(1+Q\delta),\eqdef{sufficient}$$
\par
which remains to be proved.
\medskip

{\twelvebf 3.3 Proof of Lemma 7}
\medskip

a) First we treat the case $\delta\ge 1$. System \eqref{syst3} reduces
then to 
$$\left\{
\matrix{
h_2q_2=h_1q_1+rd\cr
d\ll \delta Q\cr
}
\right.$$
from which we deduce \eqref{sufficient} at once. From now on, we
suppose $0<\delta<1.$
 
\medskip
b) We fix the integers $r,h_1$ and $h_2$. In order to apply Lemma 5,
we transform system \eqref{syst3}. We use the first line of
\eqref{syst3} to express $d$ and we report this expression into the
second line ; we divide the so obtained inequality by $q_1^2h_2(h_2-r)$
and we get :

$$\left|\frac{q_2^2}{q_1^2}-\frac{h_1(h_1-r)}{h_2(h_2-r)}\right|\le 2\delta \frac{|r|}{H},$$
since $h_2-r\ge H/2.$ Finally, system \eqref{syst3} implies
$$\left\{
\matrix{
rd+h_1q_1-h_2q_2=0\cr
\displaystyle{\frac{q_2}{q_1}=\sqrt{\frac{h_1(h_1-r)}{h_2(h_2-r)}}+O\left(\delta \frac{|r|}{H}\right)}}\right.\eqdef{q1surq2}$$

By Lemma 5, the number of triplets $(d,q_1,q_2)$ solutions of
\eqref{q1surq2} is
$$\ll H^{\ep}\left(1+\frac{\delta  Q^2}{H}\right),$$
so that we have
$${\cal J}(R,Q,H,\delta)\ll_{\ep}H^{\ep}(RH^2+RHQ^2\delta),\eqdef{maj1}$$
and this proves \eqref{sufficient} in the case $H\ll Q$.
\medskip
c) It only remains to prove Lemma 7 in the following two cases
$$0<\delta <1, \ Q\ll H {\rm \ and}\ \delta \ll R/H\eqdef{firstcase}$$
and
$$0<\delta <1, \ Q\ll H {\rm \ and}\ \delta \gg R/H\eqdef{secondcase}$$
We could fix the integers $d,q_1,q_2$ with the aim of applying Lemma 5
to bound the number of triplets $(r,h_1,h_2)$ which satisfy
\eqref{syst3}, as in the previous case. But this direct method does
not yield \eqref{sufficient} and some extra work is needed.
First we want to prove that \eqref{syst3} implies the two systems :

$$\left\{
\matrix{
dr+q_1h_1-q_2h_2=0 \cr
\displaystyle{\frac{h_1}{h_2}=\frac{q_2(q_2-d)}{q_1(q_1+d)}+O(\delta)}& {\rm if}\
\eqref{firstcase}\ {\rm holds}\cr
2d+q_1-q_2\ll RQ/H\cr
}
\right.\eqdef{syst4}$$
and

$$\left\{
\matrix{
dr+q_1h_1-q_2h_2=0 \cr
\displaystyle{\frac{h_1}{h_2}=\frac{q_2}{q_1}+O(R/H)}& {\rm if}\
\eqref{secondcase}\ {\rm holds}\cr
2d+q_1-q_2\ll Q\delta\cr
}
\right.\eqdef{syst5}$$

For this, we recall that $q_1$ and $q_2$ are of the same sign, so that
we have either $|q_1+d|\ge Q$ or $|q_2-d|\ge Q$. For example, we
assume that $q_1,q_2$ and $d$ are of the same sign. From
\eqref{syst3}, we deduce that
$$\frac{h_1}{h_2}=\frac{q_2(q_2-d)}{q_1(q_1+d)}+O(\delta).\eqdef{Fh1surh2}$$
From the first line of \eqref{syst3} and the bound $d\ll Q$, we deduce
at once
$$\frac{h_1}{h_2}=\frac{q_2}{q_1}+O(R/H).\eqdef{Sh1surh2}$$

At last, from \eqref{Fh1surh2} and \eqref{Sh1surh2}, we deduce
$$2d+q_1-q_2\ll Q\delta+RQ/H,\eqdef{dq1q2}$$

so that the systems \eqref{syst4} and \eqref{syst5} derive from
\eqref{syst3} as claimed above.
\medskip
Now, we suppose that \eqref{firstcase} holds. We suppose furthermore
that $|d|$ has a fixed size $D$, that is $D\le |d|<2D$, with $D\ll
Q$. We then fix the integers $d,q_1$ and $q_2$ with only
$O(DQ+DQ^2R/H)$ possibilities, by \eqref{dq1q2}. By lemma 5, the number
of triplets $(r,h_1,h_2)$ which satisfy \eqref{syst4} is
$O_{\ep}(Q^{\ep}(1+H^2\delta/D))$, so that the total number of integer
points $(d,q_1,q_2,r,h_1,h_2)$ lying in the domain
\eqref{paramnewcond} and satisfying \eqref{syst4} is :
$$\ll_{\ep} Q^{\ep}\max_{1\le D\ll Q}\left(1+H^2\delta/D\right)\left(DQ+DQ^2R/H\right)$$

$$\ll_{\ep} Q^{\ep}\left(Q^2+\delta H^2Q+Q^3R/H+Q^2HR\delta \right)$$
which proves \eqref{sufficient} in this case.
 The proof of lemma 7 in case \eqref{secondcase} is completely similar
 and we have only to use \eqref{syst5} instead of \eqref{syst4}. The
 proofs of lemma 7 and of Theorem 2 are complete.\qed

\paraun{Proof of Theorem 1}
\medskip
We are now going to prove Theorem 1. We may suppose that hypothesis
\eqref{VdCcond} holds with $\lambda $ small enough.
We split up the proof into short steps.
\medskip
{\bf Step 0 : the size of $M$}
\medskip
For proving Theorem 1, we may suppose that 
$$M\asymp \lambda^{-8/13}\eqdef{Mcourt}$$
(where the notation $u\asymp v$ means that we have both $u\ll v$ and
$v\ll u$).
Indeed, we set $M_0=[\lambda^{-8/13}]$. If we have $M\ge M_0$, we
divide the sum $S_M$ into $O(M\lambda^{8/13})$ shorter sums and the
problem reduces to \eqref{Mcourt}.
Now we consider the case $\lambda^{-7/13}\ll M<M_0$. We perform a
$C^4$ continuation of $f$ by setting :
$${\widetilde f}(M+t)=\sum_{j=0}^4f^{(j)}(M)\frac{t^j}{j!}$$
for $t>0$. Then by Lemma 5.2.3 of [4], we have
$$S_M\ll\max_{0\le \theta \le 1}\left|\sum_{m=1}^{M_0}\e({\widetilde
f}(m)+\theta m)\right|\log M_0$$
and the problem reduces again to \eqref{Mcourt}.

\medskip
{\bf Step 1 : A-process}
\medskip
We start with the sum $S_M=\sum_{m=1}^M\e(f(m))$ and we apply Weyl and
van der Corput's A-process in the form that uses symmetrical
differences (cf [4], Lemma 5.6.2). We set $\Delta_hf(m)=f(m+h)-f(m-h)$
and choose a positive integer $H$ such that
$$H\asymp \lambda^{-2/13}\eqdef{choixH}$$
We then have :
$$S_M^2\ll\frac{M^2}{H}+\frac{M}{H}\left|\sum_{h=1}^{H-1}\left(1-\frac{h}{H}\right)\sum_{m=h+1}^{M-h}\e(\Delta_hf(m))\right|.\eqdef{lemA}$$
Next, we remove the factor $(1-h/H)$ by partial summation and we use
the following remark : given any complex numbers $a_1,a_2,...,a_H$, there
exists a positive integer $H_1\le H/2$ such that 
$$\sum_{h=1}^{H-1}a_h\ll \left(\max_{1\le h\le
H-1}|a_h|+\left|\sum_{h=H_1}^{2H_1-1}a_h\right|\right)\log
H. \eqdef{dyadic}$$
Taking $a_h=\sum_{m=h+1}^{M-h}\e(\Delta_hf(m))$, with $|a_h|\le M$, we
get
$$S_M^2\ll\frac{M^2}{H}\log H+\frac{M}{H}|S(H_1)|\log
H,\eqdef{lemAmod}$$
for some integer $H_1\le H/2$, where we have set
$$S(H_1)=\sum_{h=H_1}^{2H_1-1}\sum_{m=h+1}^{M-h}\e(\Delta_hf(m)).\eqdef{SH1}$$
By the third derivative test (lemma 3), we see that, if $H_1\ll
\lambda^{-1/7}$, we have $S_M\ll M\lambda^{1/13}\log M$, and the
theorem is proved. Thus it remains to prove that :
$$S(H_1)\ll M^{1+\ep} \ ,\  {\rm for \ } \lambda^{-1/7}\ll H_1\ll
\lambda^{-2/13}.\eqdef{H1grand}$$
\medskip
{\bf Step 2 : A $\times $  A-process}
\medskip
We choose two integers $R$ and $Q$ such that
$$R\asymp \lambda^{-1/13} \ {\rm and} \ Q\asymp\lambda^{-3/13}.\eqdef{choixRQ}$$
We apply Lemma 1 to get
$$S(H_1)^2\ll
\frac{MH_1}{QR}\sum_{|r|<R}\sum_{|q|<Q}\left(1-\frac{|r|}{R}\right)\left(1-\frac{|q|}{Q}\right)\sum_{h\in
J_1(r)}\sum_{m\in
J_2(h,q)}\e(\Delta_hf(m+q)-\Delta_{h+r}f(m)),\eqdef{lem1}$$
where $J_1(r)$ and $J_2(h,q)$ are intervals defined by
$$J_1(r)=[\max(H_1,H_1-r),\min(2H_1-1,2H_1-1-r)],$$
$$ {\rm and}\ J_2(h,q)=[\max(1+h,1+h-q),\min(M-h,M-h-q)]$$
In the sum in \eqref{lem1}, we want to remove all terms with $r=0$ or
$q=0$ to get
$$\matrix{\displaystyle{S(H_1)^2\ll M^2+
\frac{MH_1}{QR}\sum_{0<|r|<R}\sum_{0<|q|<Q}\left(1-\frac{|r|}{R}\right)\left(1-\frac{|q|}{Q}\right)}\hfill\cr
\hfill\displaystyle{\sum_{h\in
J_1(r)}\sum_{m\in
J_2(h,q)}\e(\Delta_hf(m+q)-\Delta_{h+r}f(m)).}\cr}\eqdef{rqnonnul}$$
In order to prove \eqref{rqnonnul}, we first notice that the terms in
the sum \eqref{lem1} corresponding to $r=q=0$ have a contribution 
$$\ll \frac{(MH_1)^2}{QR}\ll M^2.$$
The terms corresponding to $r=0$ and $q\not =0$ may be treated as
exponential sums on the variable $m$, by van der Corput's inequality
([3], Theorem 2.2). Their contribution is $\ll M^2\lambda^{1/13}\ll
M^2$.
The terms corresponding to $q=0$ and $r\not =0$ may be treated
similarly, but with Lemma 3. In order to see that the hypotheses are
satisfied, we make use of Taylor's formula to write the phase as :
$$\Delta_hf(m)-\Delta_{h+r}f(m)=g(m)+u(m),$$
with $g(m)=-2rf'(m)$ and 
$$u(m)=\int_0^h\frac{(h-t)^2}{2!}(f'''(m+t)+f'''(m-t))\d
t-\int_0^{h+r}\frac{(h+r-t)^2}{2!}(f'''(m+t)+f'''(m-t))\d t.$$
If we set $\mu=|r|\lambda$, so that $|g'''(m)|\asymp \mu$, we have
$u'(m)\ll H^3\lambda\ll\mu^{1/2}$. An application of lemma 3 shows
that the contribution of these terms is $\ll M^2\lambda^{1/13}$. We
have completed the proof of \eqref{rqnonnul}.
\medskip
{\bf Step 3 : shift}
\medskip
We want to apply the following obvious equality :
$$\sum_{1\le m\le M}a(m)=\frac{1}{N}\sum_{n=1}^N\sum_{1-n\le m\le
M-n}a(m+n),$$
where $N$ is any positive integer and where $(a(m))_{1\le m\le M}$ are
any given complex numbers. Here we choose
$$N\asymp \lambda^{-3/13}.\eqdef{choixN}$$
If $g(m)$ is any real valued function, we have
$$\sum_{m\in J_2(h,q)}\e(g(m))=\frac{1}{N}\sum_{n=1}^N\sum_{m\in J_3(h,q,n)}\e(g(m+n)),$$
where $J_2(h,q)$ is as in \eqref{rqnonnul}. Set now
$J_0=[H_1+Q,M-2H_1-Q-N]$; we have $J_0\subset J_3(h,q,n)\subset
[1,M]$, so that :
$$\sum_{m\in J_2(h,q)}\e(g(m))=\frac{1}{N}\sum_{n=1}^N\sum_{m\in J_0}\e(g(m+n))+O(H_1+Q+N).$$
Inserting the above equality in \eqref{rqnonnul}, we finally deduce :
$$\matrix{\displaystyle{S(H_1)^2\ll
M^2+\frac{MH_1}{QRN}\sum_{0<|r|<R}\sum_{m\in J_0}}\hfill\cr
\displaystyle{\left|\sum_{0<|q|<Q}\left(1-\frac{|q|}{Q}\right)\sum_{h=H_1}^{2H_1-1}\sum_{n=1}^N\e(\Delta_hf(m+n+q)-\Delta_{h+r}f(m+n))\right|.}\cr}\eqdef{shift}$$
\medskip
{\bf Step 4 : Taylor's formula and partial summation}
\medskip
We write $$f(m+y)=f(m)+f'(m)y+f''(m)\frac{y^2}{2!}+f'''(m)\frac{y^3}{3!}+v_m(y),$$
with
$$v_m(y)=\int_0^y\frac{(y-t)^3}{3!}f^{(4)}(m+t)\d t.$$
Then $v_m$ is a ${\cal C}^4$ function which satisfies :
$$v_m^{(j)}(y)\ll Q^{4-j}\lambda \ , \ {\rm for } \ o\le j\le 4 {\rm \
and \ } y\ll Q.\eqdef{vm}$$
We introduce the function
$$u_{m,r}(q,h,n)=v_m(n+q+h)-v_m(n+q-h)-v_m(n+h+r)+v_m(n-h-r)\eqdef{umr}$$
and the two polynomials
$$P_1(r,q,h,n)=qh-rn \ , \
P_2(r,q,h,n)=hq^2+2hqn-rn^2-rh^2-r^2h,\eqdef{p1p2}$$
so that we have

$$\matrix{\displaystyle{
\Delta_hf(m+n+q)-\Delta_{h+r}f(m+n)=-2rf'(m)+2f''(m)P_1(r,q,h,n)}\hfill\cr
\hfill\displaystyle{+f'''(m)P_2(r,q,h,n)+\frac{r^3}{3}f'''(m)+u_{m,r}(q,h,n).}\cr}\eqdef{develop}$$
We bring this formula into \eqref{shift}. Our aim now is to remove the
term $u_{m,r}(q,h,n)$ from the triple exponential sum by partial
summation. But for the function
$$(q,h,n)\to \e(u_{m,r}(q,h,n)),$$
 the bound \eqref{abelcond} holds with $D\ll 1$, for we have
 $Q^4\lambda \ll 1$. If we use coefficients in \eqref{abel} instead of
 the sets ${\cal P}'$, we have finally obtained :
$$S(H_1)^2\ll M^2+\frac{MH_1}{QRN}{\widetilde S}(H_1),$$
where we have set
$${\widetilde
S}(H_1)=\sum_{0<|r|<R}\sum_{m=1}^M\left|\sum_{0<|q|<Q}\sum_{h=H_1}^{2H_1-1}\sum_{n=1}^Nb_r(q,h,n)\e(x_mP_1(r,q,h,n)+y_mP_2(r,q,h,n))\right|,\eqdef{srond}$$
and
$$x_m=2f''(m),\ y_m=f'''(m), {\rm for} \ m=1,...,M,\eqdef{xmym}$$
and where $b_r(q,h,n)$ are complex numbers of modulus at most
one. Taking \eqref{H1grand} into account, we see that Theorem 1 will
be proved if we obtain the bound
$${\widetilde S}(H_1)\ll_{\ep}M^{1+\ep}\lambda^{-5/13}.\eqdef{suffisant}$$
\medskip
{\bf Step 5 : double large sieve}
\medskip
We want to apply Lemma 4 to the sum ${\widetilde S}(H_1)$. We set
$$\mu=M\lambda, \ X_1=QH, \ X_2= QHN.\eqdef{mux1x2}$$
We define ${\cal N}$ and ${\cal B}$ as in \eqref{N} and \eqref{B}. The
size of parameters $M,H,H_1,R,Q$ and $N$ (cf
\eqref{Mcourt},\eqref{choixH},\eqref{H1grand},\eqref{choixRQ} and
\eqref{choixN}) shows that the hypothesis of lemma 4 are satisfied and
that \eqref{sieve} implies :
$${\widetilde S}(H_1)^2\ll RX_1\mu X_2{\cal N}{\cal B}(\log Q)^2.\eqdef{sieveimply}$$
It only remains to bound ${\cal B}$ and ${\cal N}$.
\medskip
{\bf Step 6 : a bound for ${\cal B}$}
\medskip
In \eqref{B}, we set $m_2=m$ and $m_1=m+k.$ We then have :
$${\cal B}\ll\left\{(k,m)\left|\right. \ 1\le m\le m+k\le M \ {\rm and \
such \ that \ the\  properties\ \eqref{Dkf'''}\ and\ \eqref{Dkf''} \
are\ satisfied} \right\}$$
with
$$|\overline{\Delta}_kf'''(m)|\le X_2^{-1},\eqdef{Dkf'''}$$
$$\|2\overline{\Delta}_kf''(m)\|\le X_1^{-1},\eqdef{Dkf''}$$
where we have set $\overline{\Delta}_k\phi(x)=\phi(x+k)-\phi(x)$.
\medskip
The inequality \eqref{Dkf'''} yields a bound for $k$, say $0\le k\le
K$, with $K\asymp (\lambda X_2)^{-1}$, while the inequality
\eqref{Dkf''} may be treated with respect to  $m$, with  fixed $k$, by
the first derivative test for integer points close to a curve (cf [4], Lemma
3.1.2). We then obtain :
$${\cal B}\ll M+\sum_{k=1}^K\left(\frac{M}{X_1}+Mk\lambda
+\frac{1}{k\lambda X_1}+1\right).$$
 The final bound is
$${\cal B}\ll M\log M.\eqdef{boundforB}$$
\medskip
{\bf Step 7 : A bound for ${\cal N}$}
\medskip
We have to bound the number ${\cal N}$ of integral solutions of the
system
$$\left\{\matrix{h_1q_1-rn_1=h_2q_2-rn_2 \hfill \cr
(h_1q_1^2+2h_1q_1n_1-rn_1^2-rh_1^2-r^2h_1)-(h_2q_2^2+2h_2q_2n_2-rn_2^2-rh_2^2-r^2h_2)\ll
\frac{1}{M\lambda}.\hfill\cr}\right.\eqdef{systmlambda}$$
We failed in solving this problem in its right generality. If we were
able to prove the expected bound (under some suitable restrictions), then we should obtain conjecture~1 by the
same method (we should only need to change the size of the
parameters).
\smallskip
Presently, we reduce the system \eqref{systmlambda} to the simpler one
of Theorem 2 by transferring the terms $rh_i^2$ and $r^2h_i$ in the
error term. This is possible since we have imposed the condition
$RH_1^2\ll 1/M\lambda .$

 Furthermore, the hypothesis $R\le H_1/2$ of theorem 2 is satisfied
 when $\lambda $ is small enough, because of \eqref{H1grand} and
 \eqref{choixRQ}. From  Theorem 2, we deduce
$${\cal N}\ll_{\ep} M^{\ep}\lambda^{-9/13}.\eqdef{boundforN}$$
We take back \eqref{boundforB} and \eqref{boundforN} into
\eqref{sieveimply}. This gives \eqref{suffisant} and the proof is complete.\qed

\bigskip
\bigskip
\bigskip
\centerline{\gras References}
\medskip

\medskip
[1] {\pcap E. Bombieri et H. Iwaniec,} \og On the order of
$\zeta(\frac{1}{2}+it)$ \fg,{\it Ann. Scuola Norm. Sup. Pisa Cl. Sci} (4) 13 (1986) 
\medskip
[2] {\pcap E. Bombieri et H. Iwaniec,} \og Some mean value theorems for exponential
sums\fg,{\it Ann. Scuola Norm. Sup. Pisa Cl. Sci} (4) 13 (1986) 473-486
\medskip
[3] {\pcap S. W. Graham and G.Kolesnik, } {\it Van der Corput's method for exponential
sums,} London Math. Soc. Lecture Notes Series 126 (Cambridge
University Press, 1991)
\medskip
[4] {\pcap M. N. Huxley, } {\it Area, lattice points and exponential
sums, } Clarendon Press (Oxford, 1996)
\medskip
[5] {\pcap M. N. Huxley and G. Kolesnik, }\og Exponential sums with a
large second derivative\fg

(preprint)
\medskip
[6] {\pcap H. Iwaniec and C. J. Mozzochi, } \og On the divisor and circle
problems.\fg {\it J. Number Theory } 29 (1988), 60-93
\medskip
[7] {\pcap P. Sargos, } \og Points entiers au voisinage d'une courbe,
sommes trigonom\'etriques courtes et paires d'exposants\fg, {\it
Proc. London. Math. Soc. } (3) 70 (1995) 285-312
\medskip
[8] {\pcap N. Watt, } \og A problem on square roots of integers\fg,
{\it P\'eriodica Mathematica Hungarica} Vol. 21 (1) (1990) 55-64 
\medskip

\medskip
\medskip
\medskip
\bigskip
\bigskip
{\it O. Robert et P. Sargos
\smallskip
Institut Elie Cartan
\smallskip
Universit\'e Henri Poincar\'e - Nancy I
\smallskip
BP 239
\smallskip
54 506 Vandoeuvre-l\`es-Nancy Cedex
\smallskip
  Adresse  \'electronique :
\smallskip
 Robert @ iecn.u-nancy.fr
\smallskip
 Sargos @ iecn.u-nancy.fr}

\end